\magnification=\magstep1
\input amstex
\documentstyle{amsppt}

\define\defeq{\overset{\text{def}}\to=}

\def \isom {\buildrel \sim \over \rightarrow}
\def \isomleft {\buildrel \sim \over \leftarrow}

\def \Pri{\operatorname {Pri}}
\def \Fr{\operatorname {Fr}}
\def \Div{\operatorname {Div}}

\def \lim{\operatorname {lim}}
\def \sep{\operatorname {sep}}

\def \Gal{\operatorname {Gal}}
\def \pr{\operatorname {pr}}
\def \Ker{\operatorname {Ker}}

\def \acts\ trivially\ on{\operatorname {acts\ trivially\ on}}
\def \closed\ subgroups{\operatorname {\closed\ subgroups}}
\def \Hom{\operatorname {Hom}}

\def \ur{\operatorname {ur}}

\def \sup{\operatorname {sup}}
\def \ab{\operatorname {ab}}

\def \Aut{\operatorname {Aut}}

\def \pr{\operatorname {pr}}

\def \et{\operatorname {et}}
\def \Pic{\operatorname {Pic}}

\def \Sub{\operatorname {Sub}}
\def \characteristic{\operatorname {characteristic}}
\def \char {\operatorname {char}}
\def \cl {\operatorname {cl}}
\def \div {\operatorname {div}}
\def \deg {\operatorname {deg}}

\define\OSub{\operatorname{OSub}}

\define\Primes{\frak{Primes}}

\define\GL{\operatorname{GL}}
\define\PGL{\operatorname{PGL}}
\define\nemp{\neq\varnothing}

\NoRunningHeads
\NoBlackBoxes

\document

\head
A refined version of Grothendieck's anabelian conjecture for hyperbolic curves over finite fields
\endhead
\bigskip
\centerline {MOHAMED SA\"IDI and AKIO TAMAGAWA}

\bigskip
\definition{Abstract} In this paper we prove a refined version of 
a theorem by Tamagawa and Mochizuki on isomorphisms 
between (tame) arithmetic fundamental groups of hyperbolic curves over finite fields,
where one ``ignores'' the information provided by a ``small'' set of primes.

\enddefinition

\bigskip

\noindent
\S 0. Introduction

\noindent
\S 1. Review of the local theory 

\noindent
\S 2. Large and small sets of primes relative to a hyperbolic curve over a finite field

\noindent
\S 3. Review of Mochizuki's cuspidalization theory of proper hyperbolic curves over finite fields

\noindent
\S 4. Isomorphisms between geometrically pro-$\Sigma$ arithmetic fundamental groups

\noindent
\S 5. 
On the fundamental theorem of projective geometry

\subhead
\S 0. Introduction
\endsubhead
Let $k$ be a finite field of characteristic $p>0$ and 
$U$ a hyperbolic curve over $k$. Namely, $U=X\setminus S$, where 
$X$ is a proper, smooth, geometrically connected curve 
of genus $g$ over $k$ and $S\subset X$ is a divisor 
which is finite \'etale of degree $r$ over $k$, such that 
$2-2g-r<0$. 
We have the following commutative diagram of profinite groups:
$$
\matrix
1&\to& \pi_1(U\times _k{\bar k},*)&\to& \pi_1(U,*)&\to& G_k&\to& 1\\
&&&&&&\\
&&\downarrow&&\downarrow&&\Vert\\
&&&&&&\\
1&\to& \pi_1^t(U\times _k{\bar k},*)&\to& \pi_1^t(U,*)&\to& G_k&\to& 1\\
\endmatrix
$$
in which both rows are exact and all vertical arrows are surjective 
(and bijective for $r=0$). 
Here, $G_k$ is the absolute Galois group $\Gal (\bar k/k)$, $*$ 
means a suitable geometric 
point, and $\pi _1$ (resp. $\pi_1^t$) stands for the \'etale 
(resp. tame) fundamental group. 
The following result is fundamental
in the anabelian geometry of hyperbolic curves over finite fields.

\proclaim{Theorem A (Tamagawa, Mochizuki)} Let $U$, $V$ be 
hyperbolic curves over 
finite fields $k_U$, $k_V$, respectively. Let 
$$\alpha : \pi_1(U,*)\isom \pi_1(V,*)$$
be an isomorphism of profinite groups. 
Then $\alpha$ arises from a uniquely determined 
commutative diagram of schemes:
$$
\CD
\Tilde U@>{\sim}>> \Tilde V\\
@VVV   @VVV \\
U @>{\sim}>> V\\
\endCD
$$
in which the horizontal arrows are isomorphisms and the vertical arrows are the 
profinite \'etale universal coverings determined 
by the profinite groups $\pi_1(U,*)$, $\pi_1(V,*)$, respectively.
\endproclaim

Theorem A was proved by Tamagawa (cf. [Tamagawa], Theorem (4.3)) in the affine case 
(together with the variant where $\pi_1$ is replaced by $\pi_1^t$), 
and by Mochizuki (cf. [Mochizuki1], Theorem 3.2) in the proper case. 
It implies, in particular, that one can embed a suitable category of hyperbolic 
curves over finite fields into the category of profinite groups. It is essential in 
the anabelian philosophy of Grothendieck, as was formulated in [Grothendieck], 
to be able to determine the image of this functor. Recall that the full structure 
of the profinite group $\pi_1(U\times _k{\bar k},*)$ is unknown 
(for any single example of 
$U$ which is hyperbolic). Hence, a fortiori, the structure of  
$\pi_1(U,*)$ is unknown. 
(Even if we replace the fundamental groups 
$\pi_1(U\times _k{\bar k},*)$, $\pi_1(U,*)$ by the tame fundamental groups 
$\pi_1^t(U\times _k{\bar k},*)$, $\pi_1^t(U,\ast)$, respectively, 
the situation is just the same.) 
Thus, the problem of determining the image of the above functor 
seems to be quite difficult, at least for the moment. 
In this paper we investigate the following question:

\definition{Question 0.1} 
Is it possible to prove 
any result analogous to the above Theorem A 
where $\pi_1(U,*)$ is replaced by some (continuous) quotient of 
$\pi_1(U,*)$ whose structure is better understood?
\enddefinition

Let $\Primes$ be the set of all prime numbers.
Let $\Sigma=\Sigma _X\subset \Primes$ be a set of prime numbers 
containing at least one prime number different from 
the characteristic $p$. 
Let $\Cal C$ be the full class of finite groups whose cardinality 
is divisible only by primes in $\Sigma$.
Let $\Delta_U\overset \text {def}\to=\pi_1^t(U\times _k\overline k,\ast)^{\Sigma}$ be the 
maximal pro-$\Cal C$ quotient of $\pi_1^t(U\times _k\overline k,\ast)$. 
Here, if $\Sigma$ does not contain $p$, 
the structure of  $\Delta_U$  is well understood: 
$\Delta_U$ is isomorphic to the pro-$\Sigma$ completion
of a certain well-known finitely generated discrete group (i.e., either a free 
group or a surface group). 
Let $ \Pi_U\overset \text {def}\to=\pi_1^t(U,\ast)/\Ker(\pi_1^t(U\times_k\overline k,\ast)\twoheadrightarrow 
\pi_1^t(U\times_k\overline k,\ast)^{\Sigma})$  be the corresponding quotient of $\pi_1^t(U,\ast)$. 
We shall refer to $\Pi_U$ as 
the maximal geometrically pro-$\Sigma$ quotient of the tame fundamental group 
$\pi_1^t(U,\ast)$ 
or, in short, the geometrically pro-$\Sigma$ tame fundamental group of $U$. 
(When $\Sigma$ does not contain $p$, we may and shall refer to it as 
the maximal geometrically pro-$\Sigma$ quotient of the fundamental group 
$\pi_1(U,\ast)$ 
or, in short, the geometrically pro-$\Sigma$ fundamental group of $U$.)

\definition{Question 0.2} 
Is it possible to prove any result analogous to the above Theorem A 
where $\pi_1(U,*)$ is replaced by $\Pi _U$, for some non-empty set of prime numbers 
$\Sigma$ 
containing at least one prime number different from 
the characteristic $p$?
\enddefinition

The first set $\Sigma$ to consider is the set $\Sigma\defeq \Primes 
\setminus \{\characteristic =p\}$. In this case we shall refer to 
$\Pi_U$ as the maximal 
geometrically prime-to-characteristic quotient
of the fundamental group $\pi_1(U,\ast)$. We have the following result:

\proclaim{Theorem B (A Prime-to-$p$ Version of Grothendieck's 
Anabelian Conjecture for
Hyperbolic Curves over Finite Fields)} 
Let $U$, $V$ be 
hyperbolic curves over finite fields $k_U$, $k_V$, 
respectively. 
Let $\Sigma _U \overset \text {def}\to=\Primes\setminus \{\char (k_U)\}$,  
$\Sigma _V\overset \text {def}\to=\Primes\setminus \{\char (k_V)\}$,
and write $\Pi_U$, $\Pi_V$ for the geometrically pro-$\Sigma _U$ 
\'etale 
fundamental group of $U$, and the geometrically pro-$\Sigma _V$ 
\'etale 
fundamental group of $V$, respectively.
Let
$$\alpha:\Pi_U\isom \Pi_V$$
be an isomorphism of profinite groups. Then $\alpha$ arises from a uniquely 
determined commutative diagram of schemes:
$$
\CD
\Tilde U @>{\sim}>> \Tilde V \\
@VVV   @VVV \\
U @>{\sim}>> V \\
\endCD
$$
in which the horizontal arrows are isomorphisms and the vertical arrows are the 
profinite \'etale coverings corresponding to the groups $\Pi_U$, $\Pi_V$, 
respectively. 
\endproclaim

Theorem B was proved by Sa\"\i di and Tamagawa (cf. [Sa\"\i di-Tamagawa1], Corollary 3.10).
Our main result in this paper is the following refined version of the above Theorems A and B (cf. Theorem 4.22).

\proclaim {Theorem C (A Refined Version of  the Grothendieck Anabelian Conjecture for 
Proper Hyperbolic Curves over Finite Fields)}
Let $X$, $Y$ be proper hyperbolic curves  over finite fields $k_X$, $k_Y$ 
of characteristic $p_X$, $p_Y$, respectively. Let 
$\Sigma_X, \Sigma_Y \subset \Primes $ 
be 
sets of prime numbers 
and set $\Sigma _X'\overset \text {def}\to=\Primes \setminus \Sigma_X$,
$\Sigma _Y'\overset \text {def}\to=\Primes \setminus \Sigma_Y$. 
Assume that 
neither 
the $\Sigma_X'$-adic representation 
$\rho_{\Sigma_{X}'}:G_{k_X}\to \prod_{l\in \Sigma _X'} \GL (T_l(J_X))$ 
nor 
the $\Sigma_Y'$-adic representation 
$\rho_{\Sigma_{Y}'}:G_{k_Y}\to \prod_{l\in \Sigma _Y'} \GL (T_l(J_Y))$, 
arising from the Jacobian varieties $J_X$, $J_Y$ of $X$, $Y$, respectively, 
is injective. 
Write $\Pi_X$, $\Pi_Y$ 
for the geometrically pro-$\Sigma_X$ \'etale fundamental group of $X$ 
and the geometrically pro-$\Sigma_Y$ \'etale fundamental group of $Y$, respectively. Let
$$\alpha:\Pi_X\isom \Pi_Y$$
be an isomorphism of profinite groups. Then $\alpha$ arises from a uniquely determined commutative diagram of schemes:
$$
\CD
\Tilde X @>{\sim}>> \Tilde Y\\
@VVV   @VVV \\
X @>{\sim}>> Y \\
\endCD
$$
in which the horizontal arrows are isomorphisms and the vertical arrows are the 
profinite \'etale coverings corresponding to the groups $\Pi_X$, $\Pi_Y$, 
respectively. 
\endproclaim

Note that the extra assumptions on $\Sigma_X$ and $\Sigma _Y$ in Theorem C are satisfied if $\Sigma _X'$, 
$\Sigma _Y'$ are finite.
We show that 
sets of primes $\Sigma_X$ and $\Sigma _Y$ 
satisfying the conditions in Theorem C must be 
of (natural) density $\neq 0$, while 
given any $\epsilon >0$ there exist 
sets of primes $\Sigma_X$ and $\Sigma _Y$ of (natural) density $<\epsilon$ 
satisfying the conditions in Theorem C (cf. Remark 2.8).

Theorem C above implies a ``similar'' version for affine hyperbolic curves (cf. Theorem 4.23).

\proclaim {Theorem D (A Refined Version of  the Grothendieck Anabelian Conjecture for 
(Not Necessarily Proper) Hyperbolic Curves over Finite Fields)}
Let $U$, $V$ be (not necessarily proper) hyperbolic curves
over finite fields $k_U$, $k_V$ of characteristic $p_U$, $p_V$, respectively. Let 
$\Sigma_U, \Sigma_V\subset \Primes$,
be sets of prime numbers,  
and set $\Sigma _U'\overset \text {def}\to=\Primes \setminus \Sigma_U$,
$\Sigma _V'\overset \text {def}\to=\Primes \setminus \Sigma_V$. Write $\Pi_U$, $\Pi_V$ 
for the geometrically pro-$\Sigma_U$ tame fundamental group of $U$ 
and the geometrically pro-$\Sigma_V$ tame fundamental group of $V$, respectively.  
Let
$$\alpha:\Pi_U\isom \Pi_V$$
be an isomorphism of profinite groups. 
Assume that there exist open subgroups $\Pi_{U'}\subset \Pi_U$,  $\Pi_{V'}\subset \Pi_V$, 
which correspond to each other via 
$\alpha$, i.e., $\Pi_{V'}= \alpha (\Pi_{U'})$,
corresponding to \'etale coverings $U'\to U$, $V'\to V$, such that the smooth compactifications $X'$ of $U'$ and $Y'$ of $V'$ are hyperbolic, and that 
neither the $\Sigma_U'$-adic representation 
$\rho_{\Sigma_{U}'}:G_{k_{U'}}\to \prod_{l\in \Sigma _U'} \GL (T_l(J_{X'}))$ 
nor 
the $\Sigma_V'$-adic representation 
$\rho_{\Sigma_{V}'}:G_{k_{V'}}\to \prod_{l\in \Sigma _V'} \GL (T_l(J_{Y'}))$, 
arising from the Jacobian varieties $J_{X'}$, $J_{Y'}$ of $X'$, $Y'$, respectively, 
is injective. 
(Here, 
$k_{U'}$, $k_{V'}$ denote the fields of constants of $U'$, $V'$, respectively.) 
Then $\alpha$ arises from a uniquely determined commutative diagram of schemes:
$$
\CD
\Tilde U @>{\sim}>> \Tilde V  \\
@VVV   @VVV \\
U @>{\sim}>> V \\
\endCD
$$
in which the horizontal arrows are isomorphisms and the vertical arrows are the 
profinite \'etale coverings corresponding to the groups $\Pi_U$, $\Pi_V$, 
respectively. 
\endproclaim

In what follows we explain the steps/ideas of the proof 
of Theorem C. 
Starting from an isomorphism
$$\alpha:\Pi_X\isom \Pi_Y$$
between profinite groups, one can first, using well-known results on the group-theoretic characterization of decomposition groups in arithmetic fundamental groups
as 
in [Tamagawa] (the so-called local theory), establish a set-theoretic bijection
$$\phi:X^{\cl}\setminus E_X\isom Y^{\cl}\setminus E_Y$$
between the set of closed points of $X$, $Y$, outside some ``exceptional sets'' $E_X\subsetneq X^{\cl}$ and $E_Y\subsetneq Y^{\cl}$, respectively, such that $\alpha (D_x)=D_{\phi(x)}$
for $x\in X^{\cl}\setminus E_X$ where $D_x$, $D_{\phi(x)}$ denote the decomposition group of 
$x$, $\phi(x)$ in $\Pi_X$, $\Pi_Y$, respectively (which are only defined up to conjugation).
It is not difficult to prove $p\defeq p_X=p_Y$ and $\Sigma\defeq \Sigma_X=\Sigma_Y$.
As a technical step in the proof we resort to a specific auxiliary prime number 
$l
$ and consider the $\Bbb Z_l$-extensions $k_X^{l}$, $k_Y^{l}$, of $k_X$, $k_Y$, respectively.
Let $X^l\defeq X\times _{k_X}k_X^{l}$, and $Y^l\defeq Y\times _{k_Y}k_Y^{l}$.
Write $E_{X^l}\defeq E_X\times _{k_X} k_X^{l}$ (resp. $E_{Y^l}\defeq E_Y\times _{k_Y}k_Y^{l}$), 
$\Cal O_{E_{X^l}}$,  $\Cal O_{E_{Y^l}}$ 
for the rings of rational functions on $X^l$, $Y^l$ 
whose poles are disjoint from $E_{X^l}$, $E_{Y^l}$, respectively, 
and $\Cal O_{E_{X^l}}^{\times}$,  $\Cal O_{E_{Y^l}}^{\times}$ 
the multiplicative groups of $\Cal O_{E_{X^l}}$,  $\Cal O_{E_{Y^l}}$, 
respectively.
We have a natural set-theoretic bijection $\phi^l:(X^{l})^{\cl}\setminus E_{X^l}\isom (Y^{l})^{\cl}\setminus E_{Y^l}$.
Next, 
certain finite index subgroups $H_{X^l}^{\times}$, $H_{Y^l}^{\times}$ 
of $\Cal O_{E_{X^l}}^{\times}$, $\Cal O_{E_{Y^l}}^{\times}$ are naturally 
associated with $\alpha$ via Kummer theory, such that 
$\alpha:\Pi_X\isom \Pi_Y$ induces 
a commutative diagram:
$$
\CD
H_{X^l}^{\times}/((k_X^{l})^{\times}\{\Sigma'\}) @<{\rho'}<< H_{Y^l}^{\times}
/((k_Y^{l})^{\times}\{\Sigma'\})\\
@VVV        @VVV   \\
H_{X^l}^{\times}/(k_X^{l})^{\times} @<{\bar \rho}<< H_{Y^l}^{\times}/(k_Y^{l})^{\times}\\
\endCD
$$
in which the vertical arrows are the natural surjective homomorphisms and the horizontal arrows 
are natural isomorphisms induced by $\alpha$, 
where $\Sigma'\defeq \Primes\setminus \Sigma$, and $(k_X^{l})^{\times}\{\Sigma'\}$ 
(resp. $(k_Y^{l})^{\times}\{\Sigma'\}$) is the $\Sigma'$-primary part of the multiplicative group 
$(k_X^{l})^{\times}$ (resp.   $(k_Y^{l})^{\times}$ ). 

The isomorphism $\bar \rho:
H_{Y^l}^{\times}/(k_Y^{l})^{\times}\isom H_{X^l}^{\times}/(k_X^{l})^{\times}
$ 
between subgroups of groups of principal divisors supported outside exceptional sets has 
the property that it preserves the valuations of 
functions, 
with respect to the set-theoretic bijection 
$\phi^l:(X^{l})^{\cl}\setminus E_{X^l}\isom (Y^{l})^{\cl}\setminus E_{Y^l}$. 
We think of elements of 
$\Cal O_{E_{X^l}}^{\times}/((k_X^{l})^{\times}\{\Sigma'\})$ and  $\Cal O_{E_{Y^l}}^{\times}
/((k_Y^{l})^{\times}\{\Sigma'\})$ as "pseudo-functions", i.e., classes of 
rational functions with divisor supported outside exceptional sets, modulo $\Sigma'$-primary constants.
In particular, given a pseudo-function $f'\in \Cal O_{E_{X^l}}^{\times}/((k_X^{l})^{\times}\{\Sigma'\})$ (resp. $g'\in \Cal O_{E_{Y^l}}^{\times}/((k_Y^{l})^{\times}\{\Sigma'\})$),
and a closed point $x\in X^{\cl}\setminus E_X$ (resp. $y\in Y^{\cl}\setminus E_Y$) it makes sense to consider the $\Sigma$-value $f'(x)\in (k(x)^{\times})^{\Sigma}$ 
(resp. $g'(y)\in (k(y)^{\times})^{\Sigma}$) of $f'$ (resp. $g'$) (cf. Lemma 4.5 and the discussion before it). Here, $(k(x)^{\times})^{\Sigma}$, $(k(y)^{\times})^{\Sigma}$ denote 
the maximal 
$\Sigma$-primary quotient 
of the multiplicative group of the residue fields $k(x)$, $k(y)$, respectively. 
Then the isomorphism $\rho':H_{Y^l}^{\times}/((k_Y^{l})^{\times}\{\Sigma'\})\to 
H_{X^l}^{\times}/((k_X^{l})^{\times}\{\Sigma'\})$
has the property that it preserves the $\Sigma$-value of the pseudo-functions 
with respect to the set-theoretic bijection 
$\phi^l:(X^{l})^{\cl}\setminus E_{X^l}\isom (Y^{l})^{\cl}\setminus E_{Y^l}$.
Let $R_{X^l}\defeq\langle H_{X^l}^\times\rangle$, $R_{Y^l}\defeq\langle H_{Y^l}^\times\rangle$ 
denote the abelian subgroups of $\Cal O_{E_{X^l}}$, $\Cal O_{E_{Y^l}}$ generated by 
$H_{X^l}^\times$, $H_{Y^l}^\times$, respectively. In fact, 
$R_{X^l}$, $R_{Y^l}$ are subalgebras of $\Cal O_{E_{X^l}}$, $\Cal O_{E_{Y^l}}$ 
over $k_{X}^l$, $k_{Y}^l$, respectively, having the same fields of fractions, 
and $\Cal O_{E_{X^l}}$, $\Cal O_{E_{Y^l}}$ are the normalizations of 
$R_{X^l}$, $R_{Y^l}$, respectively. 
We think of the multiplicative groups 
$H_{X^l}^{\times}/(k_X^{l})^{\times}$, 
$H_{Y^l}^{\times}/(k_Y^{l})^{\times}$ 
as subsets of the projective spaces 
$\Bbb P(R_{X^l})$, $\Bbb P (R_{Y^l})$ 
associated to the infinite-dimensional 
$k_X^{l}$-vector space $R_{X^l}$, 
$k_Y^{l}$-vector space 
$R_{Y^l}$, respectively. 
Using again, in an essential way, the fact that the set $\Sigma$ satisfies the assumptions 
in Theorem C we show that the isomorphism 
$\bar \rho: H_{Y^l}^{\times}/(k_Y^{l})^{\times}\to H_{X^l}^{\times}/(k_X^{l})^{\times}$ 
viewed as a bijection between subsets
of the projective spaces $\Bbb P(R_{Y^l})$ and $\Bbb P (R_{X^l})$ preserves ``partial'' 
collineations in the following sense: given a line $\ell\subset
\Bbb P(R_{Y^l})$ such that $\ell\cap (H_{Y^l}^{\times}/(k_Y^{l})^{\times})\neq \varnothing$ 
then there exists a unique line $\ell'\subset 
\Bbb P(R_{X^l})$ such that $\ell'\cap (H_{X^l}^{\times}/(k_X^{l})^{\times})\neq \varnothing$ and
$\bar \rho(\ell\cap (H_{Y^l}^{\times}/(k_Y^{l})^{\times}))=\ell'\cap 
(H_{X^l}^{\times}/(k_X^{l})^{\times})$.
If 
$H_{X^l}^\times/(k_X^l)^\times=\Cal O_{E_{X^l}}^{\times}/(k_X^{l})^{\times}
$,
$H_{Y^l}^\times/(k_Y^l)^\times=\Cal O_{E_{Y^l}}^{\times}/(k_Y^{l})^{\times}
$, 
and $E_X=E_Y=\varnothing$, 
then 
$\bar \rho:K_{Y^l}^{\times}/(k_Y^{l})^{\times}\to K_{X^l}^{\times}/(k_X^{l})^{\times}$ is 
a bijection between points of 
the projective spaces $\Bbb P(K_{Y^l})$ and $\Bbb P(K_{X^l})$, which preserves collineations, 
where
$K_{X^l}$ (resp. $K_{Y^l}$) is the function field of $X^l$ (resp. $Y^l$). 
Thus, by the fundamental theorem of projective geometry, it arises from a unique
semi-linear isomorphism $(K_{X^l},+)\isom (K_{Y^l},+)$. Unfortunately, at this stage 
we are even not able to prove that the exceptional sets $E_X$ and $E_Y$ 
are finite. 
This causes a very serious difficulty. To overcome this difficulty, we prove, 
in 
$\S5$, 
a refined version of the fundamental 
theorem of projective geometry, which may be of interest independently of the topic of this paper 
(cf. Theorem 5.7), and which applies well in our situation in order 
to recover the ring structures of $\Cal O_{E_{X^l}}$, $\Cal O_{E_{Y^l}}$, respectively. 
More precisely, given a commutative field $k$ we define the notion of 
an admissible 
set $\Cal S$ of subsets of $\Bbb P^1(k)$ (cf. Definition 5.4) (roughly speaking these are sets consisting of ``small'' 
subsets of $\Bbb P^1(k)$). For a subset
$\Cal U\subset \Bbb P(V)$ of a projective space $\Bbb P(V)$ associated to a $k$-vector space $V$,
we define the notion of being $\Cal S$-ample where $\Cal S$ as above is 
admissible 
(cf. Definition 5.6).
Roughly speaking, being $\Cal S$-ample means that $\Cal U$ is ``sufficiently large'' in some sense (cf. loc. cit.).
Let $\Bbb L (V)$ be the set of lines in $\Bbb P(V)$,
and $\Bbb L (V)_{\Cal U}\defeq \{\ell\in \Bbb L (V)\ \vert\ \ell\cap \Cal U\neq \varnothing\}$. 
Our main result is the following (cf. Theorem 5.7).

\proclaim{Theorem E (A Refined Version of the Fundamental Theorem of Projective Geometry)}
Let $V_i$ be a $k_i$-vector space for $i=1,2$.
Assume that $\dim_{k_i}(V_i)\geq 3$ for $i=1,2$. 
Let $U_i$ be a subset of $\Bbb P(V_i)$ for $i=1,2$, and assume that 
$U_i$ is $\Cal S_i$-ample for some 
admissible 
set $\Cal S_i$ of subsets of $\Bbb P^1(k_i)$ for $i=1,2$. 
Let $\sigma: U_1\isom U_2$ and $\tau: \Bbb L(V_1)_{U_1}\isom \Bbb L(V_2)_{U_2}$ 
be bijections such that 
for each $\ell\in\Bbb L(V_1)_{U_1}$, one has 
$\tau(\ell)_{U_2}= \sigma(\ell_{U_1})$. 
Then, each such 
$(\sigma, \tau): (U_1, \Bbb L(V_1)_{U_1})\isom (U_2, \Bbb L(V_2)_{U_2})$ 
uniquely extends to a collineation 
$(\tilde\sigma, \tilde\tau): (\Bbb P(V_1), \Bbb L(V_1))\isom (\Bbb P(V_2), \Bbb L(V_2))$. 
Thus, $\tilde \sigma$ is a bijection between projective spaces which preserves collineation.
In particular, there exists an isomorphism $\mu:k_1\isom k_2$, and a $\mu$-semi-linear  
isomorphism of abelian groups $\lambda: (V_1,+)\isom (V_2,+)$ 
that induces $(\sigma, \tau): (U_1, \Bbb L(V_1)_{U_1})\isom 
(U_2, \Bbb L(V_2)_{U_2})$. Moreover,
such an isomorphism $(\mu,\lambda)$ is unique 
up to scalar multiplication. 
\endproclaim

Theorem E applies well in our case. More precisely, applying Theorem E to the above situation 
we deduce that there exists a unique isomorphism
$\tilde \rho: \Bbb P(R_{Y^l})\isom \Bbb P(R_{X^l})$
which extends the bijection
$\bar \rho: H_{Y^l}^{\times}/(k_Y^{l})^{\times}\to H_{X^l}^{\times}/(k_X^{l})^{\times}$ 
and $\tilde \rho$ preserves collineation. In particular, 
the bijection $\tilde \rho$ arises from a $\psi _0$-isomorphism 
$$\psi: ( R_{Y^l} ,+) \isom  (R_{X^l},+),$$
where  $\psi _0:k_{Y}^{l} \isom k_{X}^{l}$ is a field isomorphism. Namely, 
$\psi$ is an isomorphism of abelian groups which is 
semilinear with respect to 
$\psi _0$ 
in the sense that $\psi(ax)=\psi_0(a)\psi(x)$ for $a\in k_{Y}^{l}$ and $x\in R_{Y^l}$. 
Further, $\psi_0$ is uniquely determined and $\psi$ is uniquely determined 
up to scalar multiplication. Moreover, if we normalize the isomorphism 
$\psi: (R_{Y^l} ,+) \isom  (R_{X^l},+),$ by the condition $\psi(1)=1$, 
it becomes a ring isomorphism such that 
the diagram
$$
\CD
R_{X^l}  @<\psi<<R_{Y^l}   \\
@AAA                    @AAA   \\
k_{X}^{l} @<\psi_0<< k_{Y}^{l}\\ 
\endCD
$$
commutes. Further, $\psi$ induces a natural commutative diagram
$$
\CD
X^l @>\psi>>   Y^l\\
@VVV     @VVV  \\
X  @>\psi>>  Y \\
\endCD
$$
where the horizontal maps are scheme isomorphisms and the vertical maps are natural morphisms.
By passing to open subgroups of $\Pi_X$ and $\Pi_Y$ which correspond to each other via $\alpha$, one constructs 
the desired scheme isomorphism $\Tilde X\isom \Tilde Y$ which is compatible with the isomorphism $\psi:X\isom Y$.
Here, one has to overcome the difficulty that the assumptions on the set $\Sigma$ in Theorem C are not preserved
by passing to open subgroups: 
even if the representation
$\rho_{\Sigma',X}:G_{k_X}\to \prod_{l\in \Sigma'} 
\GL (T_l(J_X))
$ is not injective and 
$\Pi_{X'}\subset \Pi_X$ is an open subgroup, 
the representation
$\rho_{\Sigma',X'}:G_{k_{X'}}\to \prod_{l\in \Sigma '} 
\GL (T_l(J_{X'}))
$ might be injective. We overcome this problem by introducing certain (weaker but more technical) conditions 
which are preserved by passing to open subgroups. 

In $\S1$, we review the main results of the local theory mainly from [Sa\"\i di-Tamagawa1], 
and how various invariants of the curve 
$X$ can be recovered group-theoretically from $\Pi_X$. In $\S2$, we define and discuss the notion of large 
set of primes relative to a hyperbolic curve over a finite field. In $\S3$, we review the main results of 
Mochizuki's theory
of cuspidalization of \'etale fundamental groups of proper hyperbolic curves, which plays an essential 
role in this paper. 
In $\S4$, we prove our main results: Theorem C and Theorem D. 
In 
$\S5$, 
we prove the refined version of the fundamental theorem of projective geometry: Theorem E.

\definition{Remark 0.3} 
(i) A function field version of the main results of the present paper 
is given in [Sa\"\i di-Tamagawa3]. (See also 
[Sa\"\i di-Tamagawa1] and [Sa\"\i di-Tamagawa2] for 
the special case $\Sigma=\Primes\setminus\{p\}$.)

\noindent 
(ii) At the moment of writing this paper, we do not know (even in the function field case) 
if pro-$l$ versions of the above results hold, 
namely if the above Theorems C and D 
hold (under a certain Frobenius-preserving assumption) 
in the case where $\Sigma =\{l\}$ consists of a single prime $l$ which is different form $p$.
\enddefinition

\subhead
\S 1. Review of the local theory
\endsubhead
In this section we briefly review the main results in [Sa\"\i di-Tamagawa1], $\S 1$
concerning the local theory in arithmetic 
fundamental groups of hyperbolic curves over finite fields.
Let $X$ be a proper, smooth, geometrically connected 
curve over a finite field $k=k_X$ of characteristic 
$p=p_X>0$. Write $K=K_X$ for the function field 
of $X$. 

Let $S$ be a (possibly empty) finite set of closed 
points of $X$, and set $U=U_S\defeq X\setminus S$. We assume that 
$U$ is hyperbolic. 

Fix a separable closure $K^{\sep}=K^{\sep}_X$ of $K$, 
and write $\overline{k}=\overline{k_X}$ 
for the algebraic closure of $k$ in $K^{\sep}$. 
Write 
$$G_{K}\defeq\Gal(K^{\sep}/K),$$
$$G_{k}\defeq\Gal(\overline{k}/k)$$ 
for the absolute Galois groups of $K$ and $k$, respectively. 

The tame fundamental group $\pi_1^t(U)$ 
with respect to the base point defined by $K^{\sep}$ 
(where ``tame'' is with respect to the complement of $U$ in $X$)
can be naturally identified with a quotient of 
$G_{K}$. 
Write $\Gal(K^t_U/K)$ for 
this quotient. 
(In case $S=\varnothing$, we also write $K^{\ur}_U$ for $K^t_U$.) 
It is easy to see that 
$K^t_U$ contains 
$K\overline{k}$. 

Let $\Sigma=\Sigma_X$ be a set of prime numbers that contains at least
one prime number different from $p$. Write 
$$\Sigma ^{\dag}\defeq \Sigma \setminus  \{p\}.$$ 
Thus, $\Sigma ^{\dag}\neq \varnothing$ by our assumption.
Denote by $\hat \Bbb Z^{\Sigma^\dag}$ the maximal pro-$\Sigma ^{\dag}$ 
quotient of $\hat \Bbb Z$. 
Set $\Sigma'=\Sigma'_X=\Primes\setminus\Sigma_X$. We say that 
$\Sigma$ is cofinite if $\sharp(\Sigma')<\infty$. 

We define ${\tilde K}_U$ to be the 
maximal pro-$\Sigma$ 
subextension of $K\overline{k}$ in $K^{t}_U$. 
Now, set 
$$\Pi_U=\Gal({\tilde K}_U/K),$$ 
which is a quotient of 
$\pi_1^t(U)=\Gal(K^{t}_U/K)$. 
This fits into the exact sequence 
$$1\to \Delta_U \to \Pi_U \overset{\pr_U}\to{\to} G_{k}\to 1.$$ 
Here, $\Delta_U$ is the maximal 
pro-$\Sigma$ quotient of 
$\pi_1^t({\overline U})$, where, 
for a $k$-scheme $Z$, we set $\overline Z \defeq Z \times_k \overline k$. 

Define ${\tilde X}_U$ to be the integral closure 
of $X$ in ${\tilde K}_U$. 
Define 
${\tilde U}$ to be the integral closure 
of $U$ in ${\tilde K}_U$, 
which can be naturally identified with the inverse image 
(as an open subscheme) of $U$ in ${\tilde X}_U$. Define 
${\tilde S}_U$ to be the inverse image (as a set) 
of $S$ in ${\tilde X}_U$.

For a scheme $Z$, write $Z^{\cl}$ for the set of closed points of $Z$. 
Then we have 
$$X^{\cl}=U^{\cl}\coprod S,$$ 
$$({\tilde X}_U)^{\cl}={\tilde U}^{\cl}\coprod {\tilde S}_U.$$
Moreover, 
$({\tilde X}_U)^{\cl}$ admits a 
natural action of $\Pi_U$, and 
the corresponding quotient can be naturally identified with 
$X^{\cl}$. 

For each $\tilde x \in ({\tilde X}_U)^{\cl}$, we define the 
decomposition group $D_{\tilde x}\subset \Pi_U$ 
(respectively, the inertia group $I_{\tilde x}\subset D_{\tilde x}$) 
to be the stabilizer at $\tilde x$ of the 
natural action of $\Pi_U$ on $({\tilde X}_U)^{\cl}$ 
(respectively, the kernel of the natural action of $D_{\tilde x}$ 
on $k({\tilde x})=\overline{k(x)}=\overline{k}$, where $x$ is the image of
$\tilde x$ in $X^{\cl}$). These 
groups fit into the following commutative diagram in which 
both rows are exact:
$$
\matrix
1 &\to& I_{\tilde x} &\to& D_{\tilde x} &\to& G_{k(x)} &\to& 1\\
&&&&&&&&\\
&& \cap && \cap && \cap && \\
&&&&&&&&\\
1 &\to& \Delta_U &\to& \Pi_U &\to& G_{k} &\to& 1
\endmatrix
$$
Moreover, 
$I_{\tilde x}=\{1\}$ (respectively,  
$I_{\tilde x}$ is (non-canonically) isomorphic to 
${\hat {\Bbb Z}}^{\Sigma^\dag}$), 
if ${\tilde x}\in {\tilde U}^{\cl}$ (respectively, 
${\tilde x}\in {\tilde S}_U$). 
Since $I_{\tilde x}$ is normal in $D_{\tilde x}$, 
$D_{\tilde x}$ acts on $I_{\tilde x}$ by conjugation. 
Since $I_{\tilde x}$ is abelian, this action factors through 
$D_{\tilde x}\to G_{k(x)}$ and induces a natural 
action of $G_{k(x)}$ on $I_{\tilde x}$. 

%
%
%

Let $G$ be a profinite group. Then, define $\Sub(G)$ 
(respectively, $\OSub(G)$) 
to be the set of closed (respectively, open) subgroups of $G$. 

By conjugation, $G$ acts on $\Sub(G)$. More generally, 
let $H$ and $K$ be closed subgroups of $G$ such that 
$K$ normalizes $H$. Then, by conjugation, $K$ acts on 
$\Sub(H)$. We denote by $\Sub(H)_K$ the quotient 
$\Sub(H)/K$ by this action. In particular, 
$\Sub(G)_G$ is the set of conjugacy classes of 
closed subgroups of $G$. 

For any closed subgroups $H,K$ of $G$ with $K\subset H$, 
we have 
a natural inclusion 
$\Sub(K)\subset \Sub(H)$, as well as 
a natural map
$\Sub(H)\to\Sub(K)$, $J\mapsto J\cap K$. By using this latter 
natural map, we define
$$\overline{\Sub}(G)\defeq\underset{H\in\OSub(G)}\to{\varinjlim}
\Sub(H).$$
Observe that $\overline{\Sub}(G)$ can be identified with 
the set of commensurate classes of closed subgroups of $G$. 
(Closed subgroups $J_1$ and $J_2$ of $G$ are called commensurate 
(to each other), if $J_1\cap J_2$ is open both in $J_1$ and in $J_2$.) 

With these notations, we obtain natural maps
$$D=D[U]:({\tilde X}_U)^{\cl}\to\Sub(\Pi_U), {\tilde x}\mapsto D_{\tilde x},$$
$$I=I[U]:({\tilde X}_U)^{\cl}\to\Sub(\Delta_U)\subset
\Sub(\Pi_U), {\tilde x}\mapsto I_{\tilde x},$$
which fit into the commutative diagram
$$\CD
({\tilde X}_U)^{\cl} @>{D}>> \Sub(\Pi_U) \\
\| @. @VVV \\
({\tilde X}_U)^{\cl} @>{I}>> \Sub(\Delta_U)
\endCD
$$
where the vertical arrow stands for the natural map 
$\Sub(\Pi_U)\to \Sub(\Delta_U)$, 
$J\mapsto J\cap \Delta_U$. 
By composition with the natural map 
$\Sub(\Pi_U)\to\overline{\Sub}(\Pi_U)$, $D,I$ yield 
$$\overline{D}=\overline{D}[U]:({\tilde X}_U)^{\cl}\to\overline{\Sub}(\Pi_U),$$
$$\overline{I}=\overline{I}[U]:({\tilde X}_U)^{\cl}\to\overline{\Sub}
(\Delta_U)\subset \overline{\Sub}(\Pi_U).$$
Note that, unlike the case of $D,I$, 
the maps $\overline{D},\overline{I}$ are essentially unchanged 
if we replace $U$ by any covering corresponding to an open 
subgroup of $\Pi_U$.

Since the maps $D,I$ 
are $\Pi_U$-equivariant, they 
induce natural maps 
$$D_{\Pi_U}=D[U]_{\Pi_U}:X^{\cl}\to\Sub(\Pi_U)_{\Pi_U},$$
$$I_{\Pi_U}=I[U]_{\Pi_U}:X^{\cl}\to\Sub(\Delta_U)_{\Pi_U}
\subset\Sub(\Pi_U)_{\Pi_U},$$
respectively. 

\definition{Definition 1.1} 
Let $f:A\to B$ be a map of sets. 

\noindent
(i) We define $\mu_f: B\to {\Bbb Z}\cup\{\infty\}$ by 
$\mu_f(b)=\sharp(f^{-1}(b))$. (Thus, 
$f$ is injective (respectively, surjective) if 
$\mu_f(b)\leq 1$ (respectively, $\mu_f(b)\geq 1$) for 
any $b\in B$. We also have  
$f(A)=\{b\in B\mid \mu_f(b)\geq 1\}$.)

\noindent
(ii) We say that $f$ is quasi-finite, 
if $\mu_f(b)<\infty$ for any $b\in B$. 

\noindent
(iii) We say that an element $a$ of $A$ is an exceptional element 
of $f$ (in $A$), if $\mu_f(f(a))>1$. We refer to the set of 
exceptional elements as the exceptional set. 

\noindent
(iv) We say that a pair $(a_1,a_2)$ of elements of $A$ is an exceptional 
pair of $f$ (in $A$), if $a_1\neq a_2$ and $f(a_1)=f(a_2)$ hold. 

\noindent
(v) We say that $f$ is almost injective (in the strong sense), 
if the exceptional set of $f$ is finite. 
(Observe that almost injectivity implies quasi-finiteness.) 
\enddefinition

\definition{Definition 1.2}
Denote by $E_{\tilde U}$ the exceptional set 
of $\overline{D}$ in $({\tilde X}_U)^{\cl}$. 
\enddefinition


%

\proclaim{Proposition 1.3} 
%
%
%
Let $\overline \rho$ denote the natural morphism 
${\tilde X}_U\to {\overline X}$. Then, 
for each ${\overline x}\in {\overline X}^{\cl}$, 
$\overline{D}|_{{\overline \rho}^{-1}({\overline x})}$ is injective. 
%
%
%
%
%
\endproclaim

\demo{Proof}  Cf. [Sa\"\i di-Tamagawa1], Proposition 1.8(iii). 
\qed
\enddemo

\definition{Definition 1.4}
We define $E_U$ to be the image of $E_{\tilde U}$ in $X^{\cl}$. 
(This can be identified with 
$E_{\tilde U}/\Pi_U$. 
) 
\enddefinition

Next, we shall explain how various 
invariants and structures of $U$ can be recovered 
group-theoretically (or $\varphi$-group-theoretically) from $\Pi_U$, 
in the following sense. 

\definition{Definition 1.5}
(i) We say that $\Pi=(\Pi, \Delta, \varphi_\Pi)$ is a 
$\varphi$-(profinite) group, 
if $\Pi$ is a profinite group, $\Delta$ is a closed normal subgroup 
of $\Pi$ and $\varphi_\Pi$ is an element of $\Pi/\Delta$. 

\noindent
(ii) An isomorphism from a $\varphi$-group $\Pi=(\Pi, \Delta, \varphi_\Pi)$ 
to another $\varphi$-group $\Pi'=(\Pi', \Delta', \varphi_{\Pi'})$ is 
an isomorphism $\Pi\overset{\sim}\to{\to}\Pi'$ as profinite groups that 
induces $\Delta \overset{\sim}\to{\to} \Delta'$, 
hence also $\Pi/\Delta \overset{\sim}\to{\to} \Pi'/\Delta'$, such that 
the last isomorphism sends $\varphi_\Pi$ to $\varphi_{\Pi'}$. 
\enddefinition

{}From now on, we regard $\Pi_U$ as a $\varphi$-group by 
$\Pi_U=(\Pi_U, \Delta_U,\varphi_k)$, 
where $\varphi_k$ stands for the $\sharp(k)$-th power Frobenius element 
in $G_k=\Pi_U/\Delta_U$. We shall say that an isomorphism 
$\alpha: \Pi_U\isom \Pi_{U'}$ as profinite groups 
is Frobenius-preserving, if $\alpha$ is an isomorphism 
as $\varphi$-groups. 

\definition{Definition 1.6}
(i) Given an invariant $F(U)$ 
(e.g., a number, a set of numbers, etc.) 
that depends on the isomorphism class 
(as a scheme) of a hyperbolic curve $U$ over a finite field, we say 
that $F(U)$ can be recovered group-theoretically (respectively, 
$\varphi$-group-theoretically) from $\Pi_U$, 
if any isomorphism (respectively, any Frobenius-preserving isomorphism) 
$\Pi_U\isom \Pi_V$ 
implies $F(U)=F(V)$ for two such curves $U,V$. 

\noindent
(ii) Given an additional structure ${\Cal F}(U)$ (e.g., 
a family of subgroups, quotients, elements, etc.) 
on the profinite group $\Pi_U$ that depends 
functorially on a hyperbolic curve $U$ over a finite field 
(in the sense that, for any isomorphism (as schemes) 
between two such curves $U,V$, 
any isomorphism $\Pi_U\isom\Pi_V$ 
induced by this isomorphism 
$U\isom V$ (which is unique up to composition with inner automorphisms) 
preserves the structures ${\Cal F}(U)$ and ${\Cal F}(V)$, 
we say that ${\Cal F}(U)$ can be recovered group-theoretically 
(respectively, $\varphi$-group-theoretically) from $\Pi_U$, 
if any isomorphism (respectively, any Frobenius-preserving isomorphism) 
$\Pi_U\isom\Pi_V$
between two such curves 
$U,V$ preserves the structures ${\Cal F}(U)$ and ${\Cal F}(V)$. 
\enddefinition

\proclaim{Proposition 1.7}
The following invariants and structures can be 
recovered group-theoretically from $\Pi_U$: 

\noindent
{\rm (i)} The subgroup $\Delta_U$ of $\Pi_U$, hence 
the quotient $G_k=\Pi_U/\Delta_U$. 

\noindent
{\rm (ii)} The subsets $\Sigma$ and $\Sigma^\dag$ of $\Primes$. 
\noindent

%
%
%
%
%
%
%
%
\endproclaim

\demo{Proof} 
Cf. [Sa\"\i di-Tamagawa1], Proposition 1.15(i)(ii).
\qed
\enddemo

%
%
%

%
%

Finally, we shall explain that the set of decomposition groups in $\Pi_U$ 
can be recovered group-theoretically from $\Pi_U$. First, we shall treat 
decomposition groups at points of ${\tilde S}_U$. 

\proclaim{Theorem 1.8} 
{\rm (i)} The set of inertia groups at points of ${\tilde S}_U$ 
(i.e., the image of the injective 
map $I|_{{\tilde S}_U}: {\tilde S}_U \to \Sub(\Delta_U) 
\subset \Sub(\Pi_U)$) can be recovered 
$\varphi$-group-theoretically from $\Pi_U$. 

\noindent
{\rm (ii)} The set of decomposition groups at points of ${\tilde S}_U$ 
(i.e., the image of the injective 
map $D|_{{\tilde S}_U}: {\tilde S}_U \to \Sub(\Pi_U)$) 
can be recovered $\varphi$-group-theoretically from $\Pi_U$. 
\endproclaim

\demo{Proof} 
Cf. [Sa\"\i di-Tamagawa1], Theorem 1.18.
\qed
\enddemo

Next, we shall consider decomposition groups at points of ${\tilde U}^{\cl}$. 
This is done along the lines of [Tamagawa], \S 2, but slightly more subtle 
than the case of [Tamagawa], due to the existence of the exceptional 
set $E_{\tilde U}$.

\proclaim{Theorem 1.9} The following hold.

\noindent
{\rm (i)}\ The set of decomposition groups at points of ${\tilde U}^{\cl}$ 
(respectively, ${\tilde U}^{\cl}\setminus  E_{\tilde U}$, respectively, $E_{\tilde U}$) 
(i.e., the image of the map 
$D|_{{\tilde U}^{\cl}}: {\tilde U}^{\cl} \to \Sub(\Pi_U)$ 
(respectively, 
$D|_{{\tilde U}^{\cl}\setminus E_{\tilde U}}: {\tilde U}^{\cl}\setminus E_{\tilde U} \to \Sub(\Pi_U)$, 
respectively, 
$D|_{E_{\tilde U}}: E_{\tilde U} \to \Sub(\Pi_U)$))
can be recovered $\varphi$-group-theoretically from $\Pi_U$.

\noindent
{\rm (ii)}\ The set of decomposition groups at points of $({\tilde X}_U)^{\cl}$ 
(i.e., the image of the map 
$D: ({\tilde X}_U)^{\cl} \to \Sub(\Pi_U)$) 
can be recovered $\varphi$-group-theoretically from $\Pi_U$. 
\endproclaim

\demo{Proof} 
Cf. [Sa\"\i di-Tamagawa1], Theorem 1.24 and Corollary 1.25.
\qed
\enddemo

\subhead 
\S 2. Large and small sets of primes relative to a hyperbolic curve over a finite field
\endsubhead
Throughout this section, let 
$\Sigma \subset \Primes$ 
be a set of prime numbers, 
and set 
$\Sigma '\defeq \Primes \setminus \Sigma$. Let $k$ be a finite field of characteristic $p>0$ 
and set $\Sigma^{\dag}=\Sigma\setminus\{p\}$. 
Write 
$$\hat \Bbb Z^{\Sigma^{\dag}}\defeq \prod_{l\in \Sigma^{\dag}} \Bbb Z_l.$$ 
For a prime number $l\in \Primes \setminus \{p\}$ let 
$$\chi_l:G_k\to\Bbb Z_l ^{\times}$$ 
be the $l$-adic cyclotomic character,
and define the $\Sigma$-part of the cyclotomic character by: 
$$\chi_\Sigma\defeq (\chi_l)_{l\in\Sigma^{\dag}}: G_k\to 
(\hat \Bbb Z^{\Sigma^{\dag}})^{\times}=\prod_{l\in\Sigma^{\dag}}\Bbb Z_l^{\times}.$$ 
Thus, we have 
$$\bar k^{\Ker(\chi_\Sigma)}=
k_{\Sigma}\defeq k(\zeta _{l^j}\ \vert \ l\in \Sigma^{\dag}, j\in 
\Bbb Z_{\ge 0}).$$
For a prime number $l\in \Primes$,
let $G_{k,l}\subset G_k$ be the pro-$l$-Sylow subgroup of $G_k$. 
(Recall that $G_k\simeq\hat \Bbb Z$ and $G_{k,l}\simeq\Bbb Z_l$.) 
Next, we recall the notion of $k$-small and $k$-large set of primes. (Cf. [Sa\"\i di-Tamagawa3], $\S3$ 
for more details).

%
%
%
%
%
%
%
%

\definition 
{Definition 2.1} ($k$-Small/$k$-Large Set of Primes) 
%
(i) We say that the set $\Sigma$ is $k$-small if the $\Sigma$-part $\chi_{\Sigma
}$ 
of the cyclotomic character is not injective.

\noindent
(ii) We say that the set $\Sigma$ is $k$-large if the set $\Sigma '$ is $k$-small. 
\enddefinition

%
%
%
%
%
%
%
%
%
%

Note that a $k$-large set of primes is not $k$-small, by [Sa\"\i di-Tamagawa3], Proposition 3.3. 

The following results are slight generalizations of results in [Sa\"\i di-Tamagawa3], \S 3. 
Let $X$ be a proper, smooth, geometrically connected curve over $k$, 
$f,g:X\to \Bbb P^1_k$ nonconstant $k$-morphisms, and 
$F$ a proper subfield of $\bar k$ containing $k$. 
Write $
X(F)^{\cl}
\subset X^{\cl}$ for the image of $X(F)$ in $X^{\cl}$. 

\definition {Definition 2.2}
%
%
%
%
%
%
%
We say that the pair $(f,g)$ has property 
$P_\Sigma$
(respectively, 
$Q_{F,\Sigma}$, 
$P_\Sigma$ and 
$\overline Q_{F,\Sigma}$)
if the following holds:

\medskip
\noindent
$P_\Sigma(f,g)$: 
$\exists a,b \in k^{\times}\{\Sigma'\}$, such that $f=a+bg$.

\smallskip\noindent
$Q_{F,\Sigma}(f,g)$: 
$\forall' x\in X^{\cl}\setminus 
X(F)^{\cl}
$, 
$\exists a_x,b_x\in k(x)^{\times}\{\Sigma'\}$, such that $f(x)=a_x+b_x g(x)$. 

\smallskip
\noindent
$\overline P_\Sigma(f,g)$: 
$\exists a,b \in \bar k^{\times}\{\Sigma'\}$, such that $f=a+bg$.

\smallskip\noindent
$\overline Q_{F, \Sigma}(f,g)$: 
$\forall' x\in X^{\cl}\setminus 
X(F)^{\cl}
$, 
$\exists a_x,b_x\in \bar k^{\times}\{\Sigma'\}$, such that $f(x)=a_x+b_x g(x)$. 

\medskip\noindent
Here the symbol $\forall '$ means ``for all but finitely many''. 
\enddefinition

\proclaim {Proposition 2.3} 
{\rm (i)} We have the following implications: 
$$
\matrix
P_\Sigma(f,g) &\iff & \overline P_\Sigma(f,g) \\
\Downarrow && \Downarrow \\
Q_{F,\Sigma}(f,g) &\implies &\overline Q_{F, \Sigma}(f,g) 
\endmatrix
$$
\noindent
{\rm (ii)} 
Assume that $\Sigma$ is $k$-large. 
Then we have the following implication: 
$$\overline Q_{F,\Sigma}(f,g)\implies \overline P_{\Sigma}(f,g).$$
\endproclaim

\demo{Proof} 
(i) The equivalence in the first row is given in [Sa\"\i di-Tamagawa3], 
Definition/Proposition 3.5. 
The remaining implications are immediate. 

\noindent
(ii) Similar to the proof of [Sa\"\i di-Tamagawa3], Proposition 3.11. Indeed, as in the proof of loc. cit.,
if property $\overline P_{\Sigma}(f,g)$ does not hold one deduces that 
there exists a non-empty open subscheme $V\subset X$ such that 
for every $x\in V^{\cl}\setminus 
X(F)^{\cl}
$, one has
$k(x)\subset K$ where $K/k$ is a subextension of $\bar k/k$ such that $\bar k/K$ is infinite. In particular,
for every $x\in V^{\cl}$,  
one has $k(x)\subset K$ or $k(x)\subset F$. 
This is not possible: let $\phi:V\to \Bbb A^1_k$ be 
a finite $k$-morphism, $a\in \bar k \setminus  K\cup F$, and $x\in \phi^{-1}(a)\subset V^{\cl}$, 
then $k(a)\subset k(x)\subset K\cup F$, which is absurd.
\qed
\enddemo

\proclaim {Definition/Proposition 2.4} For a pair $(f,g)$ as in the above discussion, a positive integer $m$,
and a set of prime numbers $\Sigma\subset \Primes$. We define the following properties:

\noindent
$P_{\Sigma}^{(m)}(f,g): \exists a,c \in k^{\times}\{\Sigma '\}$, such that $f=a(1+cg)^m$.

\noindent
$\overline P_{\Sigma}^{(m)}(f,g): \exists a,c \in \bar k^{\times}\{\Sigma '\}$, such that $f=a(1+cg)^m$.

\noindent
$\overline Q_{F,\Sigma}^{(m)} (f,g): \forall' 
x\in X^{\cl}\setminus 
X(F)^{\cl}
$, 
$\exists a_x, c_x\in \bar k^{\times}\{\Sigma '\}$, such that  $f(x)=a_x(1+c_xg(x))^m$.

Then:  

\noindent 
{\rm (i)} The implications 
$$P_{\Sigma}^{(m)}(f,g) \iff \overline P_{\Sigma}^{(m)}(f,g) \implies \overline Q_{F,\Sigma}^{(m)} (f,g)$$
hold.

\noindent
{\rm (ii)}  If $\Sigma$ is $k$-large, 
then the implication
$$\overline Q_{F,\Sigma}^{(m)} (f,g)\implies \overline P_{\Sigma}^{(m)}(f,g)$$
holds.
 \endproclaim

\demo{Proof}
Similar to the proof of Proposition 2.3.
\qed
\enddemo

The following is the first application of the $k$-largeness property 
to the (geometrically pro-$\Sigma$, tame) fundamental groups of 
hyperbolic curves over 
$k$. 

\proclaim
{Proposition 2.5} 
Let $U$ be a hyperbolic curve over $k$, 
$X$ the smooth compactification of $U$, 
$g$ the genus of $X$, 
$r$ the cardinality of $X_{\bar k}\setminus U_{\bar k}$, 
and $\Pi_U$ the geometrically pro-$\Sigma$ 
tame fundamental group of $U$ (cf. \S 1). 
Assume that $\Sigma$ is $k$-large. 
Then the following invariants and structures can be recovered group-theoretically from 
$\Pi_U$ (cf. Definition 1.6 for the meaning of being recovered group-theoretically). 

\noindent
{\rm (i)} The prime number $p$. 

\noindent
{\rm (ii)} The $\sharp(k)$-th power Frobenius element 
$\varphi_{k}\in G_k$. 

\noindent
{\rm (iii)} The cardinality $q\defeq \sharp(k)$ 
(or, equivalently, the isomorphism class of the finite field $k$). 
\endproclaim

\demo {Proof} 
First, consider the natural character 
$$\rho^{\det}: 
G_k\to\Aut(\bigwedge^{\max}_{\hat \Bbb Z^{\Sigma^{\dag}}}(\Delta_X^{\ab})^{\Sigma^{\dag}})
=(\hat\Bbb Z^{\Sigma^{\dag}})^{\times},$$
which can be group-theoretically recovered, by Proposition 1.7(i)(ii). 
As in [Tamagawa], Proposition 3.4 and its proof, $\rho^{\det}$ 
coincides with $\lambda\cdot(\chi_{\Sigma
})^a$, where $a=g$ 
(resp. $a=g+r-1$) for $r=0$ (resp. $r>0$), and $\lambda$ is 
a certain character with values in $\{\pm 1\}$. 
(Note that $\lambda=1$ when $r=0$.) 
It follows from the 
hyperbolicity assumption $2-2g-r<0$ that $a>0$. In particular, we have 
$(\rho^{\det})^2=(\chi_{\Sigma
})^{2a}$. 

\noindent
(i) For each $N\in\Bbb Z_{>0}$, 
let $k_N/k$ denote the unique finite subextension of 
$\bar k/k$ of degree $N$. 
Then $G_{k_N}=(G_k)^N\subset G_k$ can be recovered group-theoretically. 
Consider the coinvariant quotient $\hat \Bbb Z^{\Sigma^{\dag}}_{(\rho^{\det})^2(G_{k_N})}$ 
and define $w_{X,N}$ to be its cardinality, which is a group-theoretic invariant. 
As $(\rho^{\det})^2=(\chi_{\Sigma
})^{2a}$, this invariant is computed as: 
$$w_{X,N}= (q^{2aN}-1)_{\Sigma}=\frac{q^{2aN}-1}{(q^{2aN}-1)_{\Sigma'}},$$
where, for a positive integer $n$, 
$n=n_{\Sigma}n_{\Sigma'}$ stands for the unique decomposition where 
every prime divisor of $n_{\Sigma}$ (resp. $n_{\Sigma'}$) 
belongs to $\Sigma$ (resp. $\Sigma'$). 
Set 
$$M_0\defeq \inf\Cal M,\  \Cal M\defeq\{M\in\Bbb R_{>0}\mid
\exists C>0, \forall N\in\Bbb Z_{>0}, 
w_{X,N}\leq CM^N\},$$
which is also group-theoretic. 

We claim that $M_0=q^{2a}$. Indeed, 
since $w_{X,N}\leq q^{2aN}-1\leq q^{2aN}$, 
we have $q^{2a}\in\Cal M$. On the other hand, set 
$F_0\defeq\bar k^{\Ker((\rho^{\det})^2)}$. 
As $\Sigma$ is $k$-large, 
$F_0$ is a proper subfield of $\bar k$. 
Take a prime $l$ so that $k(l)\defeq F_0\cap k^l$ is 
finite, where $k^l$ denotes the unique 
$\Bbb Z_l$-extension of $k$ (cf. Proposition 2.13 below). 
Write $[k(l):k]=l^{n_0}$. 
Then for each $n\geq 0$, we have 
$$
w_{X,l^n}= (q^{2al^n}-1)_{\Sigma}=\frac{q^{2al^n}-1}{(q^{2al^n}-1)_{\Sigma'}}
\geq 
\frac{q^{2al^n}-1}{(q^{2al^{n_0}}-1)_{\Sigma'}}
\geq 
\frac{q^{2al^n}-1}{q^{2al^{n_0}}-1}
\geq \frac{q^{2al^n}}{q^{2al^{n_0}}}. 
$$
Thus, if $(0<)M<q^{2a}$, we have $M\notin \Cal M$. The claim now follows. 

Now, $p$ can be recovered as the unique prime divisor of $M_0=q^{2a}$. 

\noindent 
(ii) Similar to [Tamagawa], Proposition 3.4(i)(ii). 

\noindent
(iii) Similar to [Tamagawa], Proposition 3.4(iii). 
\qed
\enddemo

The notion of $k$-small/$k$-large set of primes can be naturally generalized as follows, 
by replacing the multiplicative group $\Bbb G_{m,k}$ by an abelian variety. 
Let $A$ be an abelian variety over $k$, and $T(A)=\prod_{l\in\Primes}T_l(A)$ 
the (full) Tate module of $A\times _k{\bar k}$. 
Let 
$T_{\Sigma}(A)=\prod_{l\in\Sigma}T_l(A)$ 
be 
the maximal pro-$\Sigma$ quotient
of $T(A)$. Recall that one has a natural Galois representation 
$\rho _{A,\Sigma}:G_k \to \Aut (T_{\Sigma}(A))$.

\definition 
{Definition 2.6} ($A$-Small/$A$-Large Set of Primes) 
Let $A$ be an abelian variety over $k$. 

\noindent
(i) We say that the set $\Sigma$ is $A$-small if 
the Galois representation $\rho _{A,\Sigma}:
G_k \to \Aut (T_{\Sigma}(A))$ is not injective.

\noindent
(ii) We say that the set $\Sigma$ is $A$-large if the set $\Sigma '$ is $A$-small. 
\enddefinition

\proclaim{Lemma 2.7}
Let $A$ be an abelian variety over $k$ of dimension $g>0$ and 
let $\Sigma \subset \Primes$ be a set of prime numbers.  
If $\Sigma$ is $A$-large, then $\Sigma$ is $k$-large. 
\endproclaim

\demo{Proof}
It is well-known that 
the $2g$-th exterior 
power of the representation $\rho _{A,\Sigma}$
coincides 
with the $g$-th power of $\chi_\Sigma$. 
Hence, an open subgroup of $\Ker(\rho_{A,\Sigma})$ (of index $\mid g$) 
is contained in $\Ker(\chi_{\Sigma
})$, from which the assertion follows. 
\qed\enddemo

\definition{Remark 2.8} 
Let $A$ be an abelian variety over $k$ of dimension $g>0$ and 
let $\Sigma \subset \Primes$ be a set of prime numbers.  

\noindent
(i) By Lemma 2.7 and [Sa\"\i di-Tamagawa3], Remark 3.4.1, if 
a set of primes $\Sigma\subset\Primes$ is $A$-large, then $\Sigma$ is of 
(natural) density $\neq 1$. 

\noindent
(ii) On the other hand, for any given $\epsilon >0$, there exists a set of primes 
$\Sigma\subset \Primes$ such that $\Sigma$ is $A$-large and that 
the (natural) density of $\Sigma$ 
is $<\epsilon$. 
Indeed, take a prime number $r\neq p$ 
satisfying $\frac {2g(2g+1)}{2\epsilon}<r-1$. 
Let $\Sigma\defeq 
\cup _{k=1}^{2g}\{l\in \Primes\ \vert\ l^k\equiv 1 \pmod{r}\}\cup \{r,p\}$. 
Observe that the condition $l^k\equiv 1 \pmod{r}$ is equivalent to saying that 
$l\bmod r\in \mu_k(\Bbb F_r)$, hence 
the (natural) density of $\Sigma$ is $\le \frac {\sum _{k=1}^{2g} k}{r-1}
\leq\frac {2g(2g+1)}{2(r-1)}<\epsilon$. We claim
that $\Sigma'\defeq \Primes\setminus \Sigma$ is $A$-small. Indeed, $l\in \Sigma'$ implies that $l\neq r,p$ and $r$ does 
not divide
$\sharp \GL _{2g}(\Bbb F_l)=(l^{2g}-1)(l^{2g-1}-1)\cdots(l^2-1)(l-1)l^{\frac {2g(2g-1)}{2}}$. Consider the following homomorphism
$G_k\to \Aut(T_{\Sigma'} (A))\to \Aut(T_{l} (A))\to \Aut(A[l])$ where the far left map is the 
natural $\Sigma'$-adic representation and the right maps are the natural ones. Note that the kernel of the natural
surjective map $\Aut(T_{l} (A))\to \Aut(A[l])$ is pro-$l$. In particular, the image of the $r$-Sylow subgroup of 
$G_k$ in $\Aut(A[l])$ is trivial, hence the image of the $r$-Sylow subgroup of 
$G_k$ in $\Aut(T_{\Sigma'}(A))$ is trivial and $\Sigma'$ is $A$-small.
\enddefinition

Next, let $X$ be a proper, smooth, and geometrically connected hyperbolic curve over the finite field $k$. 
We apply the notations in \S 1 to $U=X$. In particular, 
$\Delta_X$ and $\Pi_X$ denote the maximal pro-$\Sigma$ quotient of 
the geometric fundamental group $\pi_1(X\times_k\bar k)$ and the maximal 
geometrically pro-$\Sigma$ quotient of the arithmetic fundamental group 
$\pi_1(X)$, respectively. For the definition of the exceptional set $E_X\subset X^{\cl}$, 
see Definition 1.4. 
Further, let $J_X$ denote the Jacobian variety of $X$. 

\definition
{Definition 2.9} 
%
(i) We denote by $F_X$ the compositum of $k(x)$ in $\bar k$ 
for all $x\in E_X$. (Note that $F_X$ depends on $\Sigma$, as 
so does $E_X$.) 

\noindent
(ii) Let $F$ be a proper subfield of $\bar k$ containing $k$: $k\subset F\subsetneq \bar k$. 
We say that $X$ is almost $\Sigma$-separated with respect to $F$ if 
$E_X\subset 
X(F)^{\cl}
$ or, equivalently, if $F_X\subset F$. 

\noindent
(iii) We say that $X$ is almost $\Sigma$-separated if it is almost 
$\Sigma$-separated with respect to some proper subfield of $\bar k$ 
containing $k$ or, equivalently, if $F_X\subsetneq\bar k$.
\enddefinition




Let $k$, $k'$ be 
finite fields and $X$, $X'$ 
proper, smooth, geometrically connected curves over $k$, $k'$, respectively. Let 
$f: X'\to X$ be a finite, generically \'etale morphism (as schemes), 
which induces a finite separable extension  
$k'(X')/k(X)$ of function fields and a finite extension $k'/k$ of 
constant fields. (In particular, we may identify $\bar k'=\bar k$.) 
Let $L'/k(X)$ denote the Galois closure of $k'(X')/k(X)$. 

\definition
{Definition 2.10} 
(i) We say that $f$ is a $\Sigma$-covering if 
the cardinality of the finite group 
$\Gal(L'\bar k/k(X)\bar k)$ is divisible only by primes in $\Sigma$.

\noindent
(ii) We say that $f$ is tame-Galois if $k'(X')/k(X)$ is a Galois extension 
(i.e., $L'=k'(X')$) and is at most tamely ramified everywhere on $X$. 
\enddefinition

\proclaim{Proposition 2.11} Assume that $f: X'\to X$ is a $\Sigma$-covering. 

\noindent
{\rm (i)} Assume that $f$ is 
\'etale. Then 
we have $E_{X'}=f^{-1}(E_X)$. Further, 
$X$ is almost $\Sigma$-separated if and only if so is $X'$. More precisely, 
if $X'$ is almost $\Sigma$-separated with respect to $F'$, then 
$X$ is almost $\Sigma$-separated with respect to $F'$; and, if 
$X$ is almost $\Sigma$-separated with respect to $F$, then 
$X'$ is almost $\Sigma$-separated with respect to some finite extension 
of $Fk'$. 

\noindent
{\rm (ii)} Assume that $f$ is tame-Galois. 
Then we have $E_{X'}\subset f^{-1}(E_X)\cup S'$, where 
$S'\subset (X')^{\cl}$ is the non-\'etale locus of $f$. 
Further, if $X$ is almost 
$\Sigma$-separated, then so is $X'$. More precisely, 
if $X$ is almost $\Sigma$-separated with respect to $F$, then 
$X'$ is almost $\Sigma$-separated with respect to some finite extension 
of $Fk'$. 
\endproclaim

\demo{Proof}
(i) When $f$ is an \'etale $\Sigma$-covering, we have the following 
commutative diagram: 
$$\matrix
(\tilde X')^{\cl}&\overset{\overline{D}}\to\to&\overline{\Sub}(\Pi_{X'})\\
\Vert&&\Vert\\
(\tilde X)^{\cl}&\overset{\overline{D}}\to\to&\overline{\Sub}(\Pi_X) 
\endmatrix
$$
{}from which we get $E_{X'}=f^{-1}(E_X)$. It is clear that, 
if $X'$ is almost $\Sigma$-separated with respect to $F'$, then 
$X$ is almost $\Sigma$-separated with respect to $F'$. Next, assume that 
$X$ is almost $\Sigma$-separated with respect to $F$. 
Define $F'$ to be the finite extension of $F$ corresponding to 
the open subgroup $(G_F)^{d!}\subset G_F$ (of index $\mid d!$), where 
$d$ is the degree of $f$. Then $F'\supset Fk'$ 
(as $[Fk':F]\mid [k':k]\mid [k'(X'):k(X)]=d$), 
and $X'$ is almost $\Sigma$-separated with respect to $F'$. 

\noindent
(ii) 
Set $G\defeq\Gal(k'(X')/k(X))$ and $\Delta_G\defeq \Gal(k'(X')\bar k/k(X)\bar k)$. 
Let $\tilde X\to X$ (resp. $\tilde X'\to X'$) be the profinite covering corresponding to 
$\Pi_X$ (resp. $\Pi_{X'}$). 
Note that ${\tilde X}'\to X$ is a profinite Galois covering with group $\Pi_{X',G}$
which sits naturally in the following exact sequences 
$1\to \Pi_{X'}\to \Pi_{X',G}\to G\to 1$, 
$1\to \Delta _{X',G}\to \Pi_{X',G} \to G_{k}\to 1$, 
where $\Delta _{X',G}$ is defined so that the latter sequence is exact and 
sits naturally in the following exact sequence 
$1\to \Delta _{X'} \to \Delta _{X',G} \to \Delta_G\to 1$. 
Note that if we view $X$ as an orbicurve, 
being the stack-theoretic quotient of $X'$ by the action of the finite group $G$, then 
$\Pi_{X',G}$ is nothing but the geometrically pro-$\Sigma$ 
\'etale fundamental group of the orbicurve $X$. 

Now, let $x'_1\in E_{X'}\subset (X')^{\cl}$. Then there exists 
$\tilde x'_1,\tilde x'_2\in (\tilde {X'})^{\cl}$, $\tilde x'_1\neq \tilde x'_2$, 
such that $\tilde x'_1$ is above $x_1'$ and that 
$D_{\tilde x'_1}, D_{\tilde x'_2}$ 
are commensurate in $\Pi_{X'}$. 
Let $\tilde x_1$, $\tilde x_2\in \tilde X^{\cl}$, 
$x'_1,x'_2\in (X')^{\cl}$ and 
$x_1$, $x_2\in \tilde X^{\cl}$, 
the images of $\tilde x'_1,\tilde x'_2$, respectively. 
Then $D_{\tilde x_1},D_{\tilde x_2}\subset \Pi_X$ are commensurate to each other, 
hence either $\tilde x_1=\tilde x_2$ or $\tilde x_1,\tilde x_2\in E_{\tilde X}$. 
In the latter case, we have $x_1'\in f^{-1}(E_X)$, as desired. In 
the former case, 
in particular, the images of 
$\tilde x'_1$, $\tilde x'_2$ in $\overline X^{\cl}$ are equal, hence 
there exists $\sigma \in \Delta _{X',G}$
such that $\sigma\cdot\tilde x'_1=\tilde x'_2$. 
Write $Z\defeq D_{\tilde x'_1}\cap D_{\tilde x'_2}\subset \Pi_{X'}\subset \Pi_{X',G}$ for simplicity.
First, we follow the proof of [Sa\"\i di-Tamagawa1], Lemma 1.7 in order to deduce that $\sigma$ is torsion.
More precisely, let $Z_0\defeq Z\cap \sigma Z \sigma ^{-1}$. Then as in loc. cit. 
we deduce that $\sigma$ commutes with any element of $Z_0$, i.e., 
$Z_0\subset Z_{\Pi_{X',G}}(\langle \sigma\rangle )$. 
(Here, given a profinite group $G$ and a closed subgroup $H$, 
we write $Z_G(H)$ for the centralizer of $H$ in $G$.) 
Moreover, arguing by contradiction, suppose $\sigma \neq 1$ and let $\overline N$ be 
a sufficiently small 
open characteristic subgroup of $\Delta _{X',G}$ such that 
$\sigma \notin \overline N$ and set $\overline H\defeq \langle \overline N,\sigma\rangle \subset \Delta _{X',G}$. 
Then, as in loc. cit. $Z_0$ normalizes $\overline H$ and  the image of $\sigma$
in $\overline H^{\ab}$ is nontrivial and fixed by the action of $Z_0$. 
By observing the Frobenius weights in the action of $Z_0$ we deduces that 
$\langle \sigma\rangle \cap \Delta _{X'}=\{1\}$ and $\sigma$ has finite order.

If $\sigma=1$, then $\tilde x'_1=\tilde x'_2$ and we are done. 
Suppose that $\sigma\neq 1$. Then it follows from 
[Mochizuki2], Lemma 4.1(iii) that
there exists a unique closed point $\tilde y'\in {\tilde X}'$ such that 
$\langle \sigma\rangle \subset I_{\tilde y'}$, where $I_{\tilde y'}$ 
is the inertia subgroup at $\tilde y'$ (necessarily
finite). We claim that $Z_{\Pi_{X',G}}(\langle \sigma\rangle )$ is commensurate to the 
decomposition group $D_{\tilde y'}$ at $\tilde y'$. Indeed, first there exists an open subgroup
$D_{\tilde y'}^o\subset D_{\tilde y'}$ of $D_{\tilde y'}$ such that $D_{\tilde y'}^o\subset 
Z_{\Pi_{X',G}}(\langle \sigma\rangle )$, as follows easily from the fact that $D_{\tilde y'}$ 
acts by inner conjugation on its normal subgroup $I_{\tilde y'}$ and the group of 
automorphisms of $I_{\tilde y'}$ is finite. 
Second, we have a natural exact sequence $1\to Z_{\Delta_{X',G}}(\langle \sigma\rangle )
\to Z_{\Pi_{X',G}}(\langle \sigma\rangle )\to G_{k}$.
Now, our second claim is that  $Z_{\Delta_{X',G}}(\langle \sigma\rangle )$ is finite. 
Indeed, after possibly replacing the orbicurve $X$ by a suitable \'etale cover, corresponding to an
open subgroup of $\Pi_{X',G}$, we can assume that $\langle \sigma\rangle =I_{\tilde y'}$. The assertion 
then follows from the fact that $I_{\tilde y'}$ is normally terminal in
$\Delta_{X',G}$ (cf. [Mochizuki2], Lemma 4.1(i)).  This implies that the subgroups $Z_0$ and 
$D_{\tilde y'}^o$ of $Z_{\Pi_{X',G}}(\langle \sigma\rangle )$ are commensurate (cf. above exact 
sequence), hence $D_{\tilde x_1'}$ and $D_{\tilde y'}$ are commensurate in $\Pi_{X',G}$, since 
$Z_0$ is open in $D_{\tilde x_1'}$, and $D_{\tilde y'}^o$ is open in $D_{\tilde y'}$. 
In particular, $D_{\tilde x_1}$ and $D_{\tilde y}$ are commensurate in $\Pi_{X}$, 
where $\tilde y$ is the image of $\tilde y'$ in $\tilde X^{\cl}$. 
Hence, either $\tilde x_1=\tilde y$ or $\tilde x_1\in E_{\tilde X}$. 
In the former case, we have $x_1=y$, hence $x_1'\in S'$, as desired, 
In the latter case, we have $x_1\in E_X$, hence $x_1'\in f^{-1}(E_{X'})$, as desired. 
\qed\enddemo

The following application of the $J_X$-largeness property 
is crucial in later sections. 

\proclaim
{Proposition 2.12}
If $\Sigma$ is $J_X$-large, then $X$ is almost $\Sigma$-separated. More precisely, then 
$X$ is almost $\Sigma$-separated with respect to some finite extension ($\neq\bar k$) 
of $\bar k^{\Ker(\rho_{J_X,\Sigma'})}$. 
\endproclaim

\demo{Proof}
Let $F_0\defeq \bar k^{\Ker (\rho _{J_X,\Sigma'})}$. 
Then $k\subset F_0\subsetneq \bar k$. 
Let $F$ denote the finite extension of $F_0$ corresponding to the 
open subgroup $(G_{F_0})^{e}\subset G_{F_0}$ 
(of index $\leq  e$), 
where $e=2$ (resp. $e=1$) 
if $X$ is hyperelliptic (resp. otherwise). 
Now, we claim that the field $F$ satisfies the above property. 

Consider the morphism
$\delta :X\times X\to J_X$, $(P,Q)\mapsto \cl (P-Q)$. Then 
$\delta\vert _{X\times X\setminus  \iota (X)}$ is quasi-finite 
with geometric fibers of cardinality $\leq e$, where $\iota : X\to X\times X$ is the diagonal map 
(cf. [Sa\"\i di-Tamagawa1], Claim 1.9(i) and its proof). 
Thus, we have a quasi-finite map $X(\bar k)\times X(\bar k)\setminus  \iota (X)(\bar k)\to J_X(\bar k)$ 
with fibers of cardinality $\leq e$. 

Let $\tilde x,\tilde x'\in \tilde X^{\cl}$ 
be an exceptional pair of the map $\overline D$ and 
$\bar x,\bar x'\in\overline X^{\cl}=X(\bar k)$ the images of 
$\tilde x,\tilde x'$, respectively. 
By Proposition 1.3, we have $\bar x\neq \bar x'$. 
Then 
the image of $(\bar x,\bar x')\in X(\bar k)\times X(\bar k)\setminus\iota(X)(\bar k) $ 
via the composed map $X(\bar k)\times X(\bar k)\setminus  \iota (X)(\bar k)
\overset\delta \to{\to} J_X(\bar k)
\twoheadrightarrow J_X(\bar k)/(J_X(\bar k)\{\Sigma'\})$ is trivial (cf. [Sa\"\i di-Tamagawa1], proof 
of Proposition 1.8(vi)). Thus, the image of $(\bar x,\bar x')$ in $J_X(\bar k)$ lies in the subgroup 
$J_X(\bar k)\{\Sigma'\}$ which is contained in $J_X(F_0)$ by the choice of $F_0$. 
It follows from this that $(\bar x,\bar x')\in X(F)\times X(F)$ since the above map 
$\delta\vert _{X\times X \setminus \iota (X)}$ is quasi-finite with geometric 
fibers of cardinality $\leq e$. 
\qed
\enddemo

%

\proclaim{Proposition 2.13} 
%
{\rm (i)} If $X$ is almost $\Sigma$-separated, then there exists a prime number 
$l
$, such that, for every finite extension $k'/k^{l}$ and every finite extension 
$F/F_X$, the field $k'\cap F$ is finite, where $k^l$ denotes the unique 
$\Bbb Z_l$-extension of $k$. (We shall refer to such a prime number $l$ as being 
$(X,\Sigma)$-admissible.) 
In particular, the set 
$E_X\cap 
X(k')^{\cl}
$ is finite.

\noindent
{\rm (ii)} If, moreover, $\Sigma$ is $J_X$-large (resp. $\Sigma$ is $J_X$-large and 
$J_X$ has positive $p$-rank), then the prime number $l$ in (i) can be 
chosen to be in $\Sigma\cup\{p\}$ (resp. $\Sigma$). 
\endproclaim

\demo{Proof}
(i) As $1\neq G_{F_X} \subset G_k=\hat\Bbb Z=\prod_{l\in\Primes}\Bbb Z_l$, there exists 
an $l\in\Primes$ such that the image of $G_{F_X}$ in $\Bbb Z_l$ is nontrivial, or, equivalently, 
open in $\Bbb Z_l$. Such an $l$ satisfies the desired property. 

\noindent
(ii) By Proposition 2.12, there exists an open subgroup $1\neq H\subset \Ker (\rho _{J_X,\Sigma'})$ 
that is contained in $G_{F_X}$. So, as in the proof of (i), take any $l\in\Primes$ such that 
the image of $H$ under the projection $G_k=\hat\Bbb Z\twoheadrightarrow\Bbb Z_l$ is nontrivial. 
Then the assertions in (i) hold for this $l$. 

Now, assume that $\Sigma$ is $J_X$-large (resp. $\Sigma$ is $J_X$-large and 
$J_X$ has positive $p$-rank), and suppose $l\in \Sigma'\setminus\{p\}$ 
(resp. $l\in\Sigma'$). Then the image of the $l$-adic representation 
$\rho_{J_X,\{l\}}: G_k\to \GL (T_l(J_X))$ is infinite and almost pro-$l$, hence the image of 
$\Ker (\rho _{J_X,\Sigma'})\subset \Ker (\rho _{J_X,\{l\}})$ in $\Bbb Z_l$ is trivial. 
Thus, we must have $l\in\Sigma\cup\{p\}$ (resp. $l\in\Sigma$) automatically. 
\qed
\enddemo

%
%
%

\subhead
\S 3. Review of Mochizuki's cuspidalization theory of proper hyperbolic curves over finite fields
\endsubhead
In this \S, we review the main results of Mochizuki's theory 
of cuspidalizations of arithmetic fundamental groups of proper hyperbolic curves over finite fields, 
developed in [Mochizuki1], which plays an important role in this paper. 
We maintain the notations of \S 1 and further assume 
$X=U$. (Thus, the finite set $S$ in \S 1 is empty, and, in this \S, 
we save the symbol $S$ for another finite set of closed points of $X$.) 
Accordingly, $X$ is a proper hyperbolic curve over a finite field 
$k=k_X$.

Recall that $\Delta _X$ stands for the maximal 
pro-$\Sigma$ quotient of $\pi _1(\overline X)$, that 
$\Pi _X$ stands for 
$\pi _1(X)/\Ker (\pi _1(\overline X)\twoheadrightarrow \Delta_X)$, and that 
they fit into the following exact sequence: 
$$1\to \Delta _X \to \Pi _X@>\pr_X>> G_k\to 1.$$

Similarly, if we write $X\times X\overset \text {def}\to=X\times _kX$, 
then we obtain (by considering the maximal pro-$\Sigma$ quotient $\Delta _{X\times X}$ 
of $\pi _1(\overline {X\times X}$)) an exact sequence:
$$1\to \Delta _{X\times X} \to \Pi _{X\times X}\to G_k\to 1,$$
where $\Pi _{X\times X}$ (respectively, $\Delta _{X\times X}$) may be 
identified with $\Pi _X\times _{G_k}\Pi_X$ (respectively, $\Delta _X\times 
\Delta_X$).

\definition{Definition 3.1} (cf. [Mochizuki1], Definition 1.1(i).) 
Let $H$ be a profinite 
group equipped with a homomorphism $H\to \Pi _X$.
Then we shall refer to the kernel $I_H$ of $H\to 
\Pi _X$ as the cuspidal subgroup of $H$ (relative to $H\to \Pi _X$).
We shall refer to an inner automorphism of $H$ by an element of $I_H$ as a
cuspidally inner automorphism.
We shall say that $H$ is cuspidally abelian (respectively, 
cuspidally pro-$\Sigma ^{*}$, where $\Sigma ^{*}$ is a set of prime numbers)
(relative to $H\to \Pi _X$) if $I_H$ is abelian (respectively, if $I_H$ is a 
pro-$\Sigma ^{*}$ group). If $H$ is cuspidally abelian, then observe that 
$H/I_H$ acts naturally (by conjugation) on $I_H$. We shall say that $H$ is 
cuspidally central (relative to $H\to \Pi _X$) if this action of 
$H/I_H$ on $I_H$ is trivial. Also, we shall use the same terminology 
for  $H\to \Pi _X$ when $\Pi _X$ is replaced by $\Delta _X$, 
$\Pi _{X\times X}$, or $\Delta _{X\times X}$.
\enddefinition

For a finite subset $S\subset X^{\cl}$ write $U_S\overset \text {def}
\to=X\setminus S$. Let $\Delta _{U_S}$ be the maximal cuspidally
(relative to the natural map to $\Delta _X$) 
pro-$\Sigma ^{\dag}$ quotient of the maximal pro-$\Sigma$ quotient of 
the tame fundamental group of $\overline {U_{S}}$ (where ``tame'' is 
with respect to the complement of $U_S$ in $X$), and let $\Pi 
_{U_S}$ be the corresponding quotient $\pi _1(U_S)/\Ker (\pi _1
(\overline {U_S})\twoheadrightarrow \Delta _{U_S}$) of $\pi_1(U_S)$. 
Thus, we have an exact sequence:  
$$1\to \Delta _{U_S}\to \Pi _{U_S} \to G_k\to 1,$$
which fits into the following commutative diagram: 
$$\matrix
1&\to &\Delta _{U_S}&\to &\Pi _{U_S} &\to &G_k&\to& 1\\
&&&&&&&&\\
&& \downarrow && \downarrow && \Vert && \\
&&&&&&&&\\
1&\to &\Delta _{X}&\to &\Pi _{X} &\to& G_k&\to& 1
\endmatrix
$$

Further, let $\iota :X \to X\times X$ be the diagonal
morphism, and write
$$U_{X\times X}\overset \text {def}\to=X\times X\setminus \iota (X).$$
We shall denote by $\Delta _{U_{X\times X}}$ the maximal cuspidally 
(relative to the natural map to $\Delta _{X\times X}$) 
pro-$\Sigma^{\dag}$ quotient of the maximal pro-$\Sigma$ quotient of 
the tame fundamental group of 
$(U_{X\times X})_{\bar k}$ (where ``tame'' is 
with respect to the divisor $\iota (X)\subset X\times X$), and by $\Pi 
_{U_{X\times X}}$ the corresponding quotient $\pi _1(U_{X\times X})/\Ker (\pi _1
(\overline {U_{X\times X}})\twoheadrightarrow \Delta _{U_{X\times X}}$) of 
$\pi _1(U_{X\times X})$.
Thus, we have an exact sequence:
$$1\to \Delta _{U_{X\times X}}\to \Pi _{U_{X\times X}}\to G_k\to 1,$$
which fits into 
the following commutative diagram: 
$$\matrix
1&\to &\Delta _{U_{X\times X}}&\to &\Pi _{U_{X\times X}}&\to& G_k&\to& 1\\
&&&&&&&&\\
&&\downarrow&&\downarrow&&\Vert&&\\
&&&&&&&&\\
1&\to &\Delta _{X\times X} &\to &\Pi _{X\times X}&\to &G_k&\to& 1 
\endmatrix
$$

Finally, set
$$M_X\overset \text {def}\to=\Hom _ {\hat \Bbb Z^{\Sigma^{\dag}}}(H^2(\Delta_X,
\hat \Bbb Z^{\Sigma^{\dag}}),\hat \Bbb Z^{\Sigma^{\dag}}).$$
Thus, $M_X$ is a free $\hat \Bbb Z^{\Sigma^{\dag}}$-module of rank $1$, 
and $M_X$
is isomorphic to $\hat \Bbb Z^{\Sigma^{\dag}}(1)$ as a $G_k$-module 
(where the ``$(1)$'' denotes a ``Tate twist'', 
i.e., $G_k$ acts on $\hat \Bbb Z^{\Sigma^{\dag}}
(1)$ via the cyclotomic character) (cf. [Mochizuki1], the discussion 
following Proposition 1.1). 

For the rest of this \S, let $X$, $Y$ be proper, 
hyperbolic curves over finite fields $k_X$, $k_Y$ of characteristic
$p_X$, $p_Y$, respectively. Let $\Sigma_X$ (respectively, $\Sigma_Y$)
be a set of prime numbers that contains at least one prime number
different from $p_X$ (respectively, $p_Y$).
Write $\Delta _X$ (respectively, $\Delta _Y$) 
for the maximal pro-$\Sigma _X$ quotient of $\pi _1(\overline X)$ 
(respectively, the maximal  pro-$\Sigma_Y$ quotient of $\pi _1
(\overline Y)$), and $\Pi _X$ (respectively, $\Pi _Y$) for the 
quotient $\pi_1(X)/\Ker (\pi_1(\overline X)\twoheadrightarrow \Delta_X)$
of $\pi _1(X)$ (respectively, the quotient 
$\pi_1(Y)/\Ker (\pi_1(\overline Y)\twoheadrightarrow \Delta_Y)$ of $\pi _1(Y)$). 

Let 
$$\alpha :\Pi _X\isom \Pi _Y$$
be an isomorphism of profinite groups.

The following is one of the main results of Mochizuki's theory 
(cf. [Mochizuki1], Theorem 1.1(iii)).

\proclaim{Theorem 3.2} (Reconstruction of Maximal Cuspidally Abelian Extensions)
Let $\iota _X:X\to X\times X$ (respectively, 
$\iota _Y:Y\to Y\times Y$) be the diagonal morphism, and write
$U_{X\times X}\overset \text {def}
\to=X\times X\setminus  \iota (X)$ (respectively, $U_{Y\times Y}
\overset \text {def}\to=Y\times Y\setminus  \iota (Y)$). Denote by 
$\Pi _{U_{X\times X}}\twoheadrightarrow\Pi _{U_{X\times X}}^{\text{\rm c-ab}}$, 
$\Pi _{U_{Y\times Y}}\twoheadrightarrow\Pi _{U_{Y\times Y}}^{\text{\rm c-ab}}$ 
the maximal cuspidally (relative to 
the natural surjections  $\Pi _{U_{X\times X}}\twoheadrightarrow \Pi_
{X\times X}$, $\Pi _{U_{Y\times Y}}\twoheadrightarrow \Pi_{Y\times Y}$, 
respectively)
abelian quotients. Then there is a commutative diagram:
$$
\CD
\Pi _{U_{X\times X}}^{\text{\rm c-ab}} @>{\alpha^{\text{\rm c-ab}}}>> \Pi _{U_{Y\times Y}}^{\text{\rm 
c-ab}} \\
@VVV                @VVV \\
\Pi_{X\times X}   @>{\alpha\times \alpha}>>  \Pi _{Y\times Y}\\
\endCD
$$
where $\alpha^{\text{\rm c-ab}}$ is an isomorphism which is well-defined up to 
cuspidally inner automorphism (i.e., an inner automorphism of $\Pi _ 
{U_{Y\times Y}}^{\text{\rm c-ab}}$ by an element of the cuspidal subgroup 
$\Ker (\Pi _{U_{Y\times Y}}^{\text{\rm c-ab}}\twoheadrightarrow 
\Pi _{Y\times Y})$). Moreover, the correspondence 
$$\alpha \mapsto \alpha ^{\text{\rm c-ab}}$$
is functorial (up to cuspidally inner automorphism) with respect to $\alpha$.
\endproclaim

\demo {Proof}
See [Mochizuki1], Theorem 1.1(iii). (See also [Sa\" \i di-Tamagawa1], Theorem 2.2.) 
\qed
\enddemo

Let $\tilde x\in \tilde X^{\cl}$ and $x$ the image of $\tilde x$ in $X$.
In this and next \S\S, we sometimes refer to the decomposition group
$D_{\tilde x}$ as the decomposition group of $\Pi_X$ at $x$, and denote it simply
by $D_x$. Thus, $D_x$ is well-defined only up to conjugation by an element of $\Pi_X$. 

For the rest of this \S, we shall assume 
that $\alpha$ is Frobenius-preserving (cf. Definition 1.5).
(Note that this assumption is automatically satisfied 
in the case where 
$\Sigma_X$ is $k_X$-large and $\Sigma_Y$ is $k_Y$-large, cf. Proposition 2.5(ii).) 
Thus, by Theorem 1.9, 
one deduces naturally from $\alpha$ a bijection
$$\phi:X^{\cl}\setminus E_X\isom Y^{\cl}\setminus E_Y$$ 
such that 
$$\alpha (D_x)=D_{\phi (x)}$$ 
holds (up to conjugation) for any $x\in X^{\cl}\setminus  E_X$. 
Here, considering the images of $D_x$ and $D_{\phi(x)}$ in 
$G_{k_X}=\Pi_X/\Delta_X\isom \Pi_Y/\Delta_Y=G_{k_Y}$ (cf. Proposition 1.7(i)), 
we obtain $[k(x):k_{X}]=[k(\phi(x)):k_Y]$, hence 
$\sharp(k(x))=\sharp(k(\phi(x)))$ by Proposition 2.5(iii). 

As an important consequence of Theorem 3.2 we deduce the following:

\proclaim{Corollary 3.3} With the above assumptions,
let $S\subset X^{\cl}\setminus E_X$,  
$T\subset Y^{\cl}\setminus E_Y$ be finite subsets that correspond 
to each other via $\phi$. Then $\alpha$, $\alpha ^{\text{\rm c-ab}}$ induce 
isomorphisms (well-defined up to cuspidally inner automorphisms, i.e.,
inner automorphisms by elements of $\Ker (\Pi _{V_T}^{\text{\rm c-ab}} \to \Pi_Y)$)
$$\alpha_{S,T}^{\text{\rm c-ab}}: 
\Pi_{U_S}^{\text{\rm c-ab}}\isom \Pi _{V_T}^{\text{\rm c-ab}}$$
lying over $\alpha$, where $U_S\defeq X\setminus S$, $V_T\defeq Y\setminus T$, 
and $\Pi_{U_S}\twoheadrightarrow \Pi_{U_S}^{\text{\rm c-ab}}$,
$\Pi_{V_T}\twoheadrightarrow \Pi_{V_T}^{\text{\rm c-ab}}$, are the maximal cuspidally 
abelian quotients (relative to the maps $\Pi_{U_S}\twoheadrightarrow \Pi_X$,
$\Pi_{V_T}\twoheadrightarrow \Pi_Y$, respectively).
These isomorphisms are functorial with respect to $\alpha$, $S$, 
$T$, as well as with respect to passing to connected finite \'etale coverings
of $X$, $Y$, which arise from open subgroups of $\Pi_X$, $\Pi_Y$, 
in the following sense:
Let $\xi:X'\to X$ (respectively, $\eta : Y'\to Y$) be a finite \'etale covering
which arises from the open subgroup $\Pi_{X'}\subset \Pi_X$ (respectively,  
$\Pi_{Y'}\subset \Pi_Y$), such that $\alpha (\Pi_{X'})=\Pi_{Y'}$; 
set $U'_{S'}\defeq X'\setminus S'$, $V'_{T'}\defeq Y'\setminus T'$, 
$S' \defeq \xi ^{-1} (S)$, $T' \defeq \eta^{-1} (T)$; 
and denote by $\alpha'$ the isomorphism 
$\Pi_{X'}\isom\Pi_{Y'}$ induced by $\alpha$. 
Then we have the following commutative diagram:
$$
\CD
\Pi_{U'_{S'}}^{\text{\rm c-ab}} @>{(\alpha')_{S',T'}^{\text{\rm c-ab}}}>> 
\Pi_{V'_{T'}}^{\text{\rm c-ab}}  \\
@VVV                            @VVV        \\
\Pi_{U_S}^{\text{\rm c-ab}} @>{\alpha_{S,T}^{\text{\rm c-ab}}}>> \Pi_{V_{T}}^{\text{\rm c-ab}}  \\
\endCD
$$
where the vertical arrows are the natural maps. 
\endproclaim

\demo {Proof} 
The proof of [Mochizuki1], Theorem 2.1(i) 
(where $E_X=E_Y=\varnothing$ is assumed) 
works as it is. See also [Sa\"\i di-Tamagawa1], Corollary 2.3. 
\qed
\enddemo

Next, let 
$$1\to M_X\to \Cal D\to \Pi_{X\times X}\to 1$$
be a fundamental extension, i.e., an extension whose corresponding 
extension class in $H^2_{\et}(X\times X,M_X)$ (via the natural identification
$H^2(\Pi_{X\times X}, M_X)\isom H^2_{\et}(X\times X,M_X)$ 
(cf. [Mochizuki1], Proposition 1.1)) coincides with the first 
(\'etale) Chern class of the diagonal $\iota (X)$ (cf. [Mochizuki1], Proposition 1.5). 
Let $x,y\in X(k)$ and 
write $D_x,D_y\subset \Pi_X$ for the associated 
decomposition groups (which are well-defined up to conjugation). Set
$$\Cal D_x\overset \text {def}\to=\Cal D|D_x\times _{G_k}\Pi_X,\ \ \ 
\Cal D_{x,y}\overset \text {def}\to=\Cal D|D_x\times _{G_k}D_y.$$
Thus, $\Cal D_x$ (respectively, $\Cal D_{x,y}$) is an extension of $\Pi_X$
(respectively, $G_k$) by $M_X$. Similarly, if $D=\sum _{i}m_i.x_i,\  
E=\sum_jn_j.y_j$ are divisors on $X$ supported 
on $k$-rational points, then set 
$$\Cal D_D\overset \text {def}\to=\sum _i m_i.\Cal D_{x_i},\  
\Cal D_{D,E}\overset \text {def}\to=\sum _{i,j} m_i.n_j.\Cal D_{x_i,y_j}$$
where the sums are to be understood as sums of extensions of $\Pi _X$, 
$G_k$, respectively, 
by $M_X$, i.e., the sums are induced by the additive structure of $M_X$.

For a finite subset $S\subset X(k)$, we shall write
$$\Cal D_S\overset \text {def}\to=\underset {x\in S}\to \prod \Cal D_x$$
where the product is to be understood as a fiber product over $\Pi_X$.
Thus, $\Cal D_S$ is an extension of $\Pi_X$ by a product of copies of $M_X$ 
indexed by the points of $S$. We shall refer to  $\Cal D_S$ as the 
$S$-cuspidalization of $\Pi_X$.
Observe that if $T\subset X(k)$
is a finite subset containing $S$, then we obtain a natural  
projection morphism $\Cal D_T\to\Cal D_S$. More generally, for a 
finite subset $S\subset X^{\cl}$ which does not 
necessarily consist of $k$-rational
points, one can still construct the object $\Cal D_S$ by passing to a finite 
extension $k_S$ of $k$ over which the points of $S$ are rational, performing the
above construction over $k_S$, and then descending to $k$. 
(See [Mochizuki1], Remark 5 for more details.) 

\proclaim {Proposition 3.4} 
(Maximal Geometrically Cuspidally Central Quotients)

\noindent
{\rm (i)}\ For $S\subset X^{\cl}$ a finite subset, the $S$-cuspidalization 
$\Cal D_S$ of $\Pi_X$ may be identified with the 
quotient $\Pi_{U_S}\twoheadrightarrow \Pi_{U_S}^{\text{\rm c-cn}}
\overset \text {def}\to=\Pi_{U_S}/ \Ker (\Delta _{U_S}\twoheadrightarrow \Delta _{U_S}^{\text{\rm 
c-cn}})$
of $\Pi_{U_S}$, where $\Delta _{U_S}^{\text{\rm c-cn}}$ is the maximal cuspidally
central quotient of $\Delta _{U_S}$ relative to the natural map 
$\Delta _{U_S} \twoheadrightarrow \Delta _X$.

\noindent
{\rm (ii)}\ The fundamental extension $\Cal D$ may be identified with the 
quotient $\Pi_{U_{X\times X}}\twoheadrightarrow 
\Pi_{U_{X\times X}}^{\text{\rm c-cn}}\overset \text {def}\to=\Pi_{U_{X\times X}}/\Ker 
(\Delta _{U_{X\times X}}\twoheadrightarrow   \Delta _{U_{X\times X}}
^{\text{\rm c-cn}})$ of $\Pi_{U_{X\times X}}$, where $\Delta _{U_{X\times X}}^{\text{\rm c-cn}}$ is
the maximal cuspidally central quotient of $\Delta _{U_{X\times X}}$
relative to the natural map $\Delta _{U_{X\times X}}\twoheadrightarrow \Delta_{X\times X}$.
\endproclaim

\demo 
{Proof}\ See [Mochizuki1], Proposition 1.6(iii)(iv). 
(Precisely speaking, Proposition 1.6(iii) loc. cit. only treats the 
special case where $S\subset X(k)$ holds. However, the proof for 
the general case is easily reduced to this special case by passing 
to a finite extension of $k$. cf. Remark 5, loc. cit.) See also 
[Sa\"\i di-Tamagawa1], Proposition 2.4. 
\qed
\enddemo

\definition {Remark 3.5}
Let $\Cal D$ (respectively, $\Cal E$) be a fundamental extension 
of $X$ (respectively, $Y$). 
The isomorphism $\alpha:\Pi_X\isom \Pi_Y$ induces an isomorphism:
$$\Cal D\isom \Cal E$$
up to cyclotomically inner automorphisms (i.e., inner automorphisms 
by elements of $M_X, M_Y$) and the actions of $(k_X^{\times})^{\Sigma_X^{\dag}}, 
(k_Y^{\times})^{\Sigma_Y^{\dag}}$, where $(k_X^{\times})^{\Sigma_X^{\dag}}$ 
(respectively, $(k_Y^{\times})^{\Sigma_Y^{\dag}}$) is the maximal $\Sigma_X^{\dag}$-
(respectively, $\Sigma_Y^{\dag}$-) quotient of $k_X^{\times}$ (respectively, 
$k_Y^{\times}$) (cf. [Mochizuki1], Proposition 1.4(ii)). 
Moreover, let $S\subset X^{\cl}\setminus E_X$ and $T\subset Y^{\cl}\setminus E_Y$ be 
as in Corollary 3.3 and write $\Cal D_S$ (respectively, $\Cal E_T$) 
for the $S$-cuspidalization of $\Pi_X$ 
(respectively, the $T$-cuspidalization of $\Pi_Y$). 
Then the isomorphism $\Cal D\isom \Cal E$ induces an isomorphism 
$$\Cal D_S\isom \Cal E_T$$
lying over $\alpha$. 

On the other hand, 
let $\Pi_{U_S}\twoheadrightarrow \Pi_{U_S}^{\text{\rm c-cn}}$ and
$\Pi_{V_T}\twoheadrightarrow \Pi_{V_T}^{\text{\rm c-cn}}$ be 
the maximal geometrically cuspidally central quotients 
(here, $U_S\overset \text {def}\to=X\setminus S$, 
$V_T\overset \text {def}\to=Y\setminus T$) (cf. Proposition 3.4). 
Note that the isomorphism 
$\alpha_{S,T}^{\text{\rm c-ab}}: 
\Pi_{U_S}^{\text{\rm c-ab}}\isom \Pi _{V_T}^{\text{\rm c-ab}}$ 
in Corollary 3.3 naturally induces an isomorphism
$$\Pi_{U_S}^{\text{\rm c-cn}}\isom \Pi _{V_T}^{\text{\rm c-cn}}$$ 
lying over $\alpha$, 
which is well-defined up to cuspidally inner automorphism. 
Now, by Proposition 3.4(i), $\Pi_{U_S}^{\text{\rm c-cn}}$ 
(respectively, $\Pi_{V_T}^{\text{\rm c-cn}}$) may be identified with 
$\Cal D_S$ (respectively, $\Cal E_T$). 
Thus, we deduce another isomorphism 
$$\Cal D_S\isom \Cal E_T$$
lying over $\alpha$. 

Now, the above two isomorphisms between $\Cal D_S$ and $\Cal E_T$ 
coincide with each other up to cyclotomically inner automorphisms 
and the actions of $(k_X^{\times})^{\Sigma_X^{\dag}}, 
(k_Y^{\times})^{\Sigma_Y^{\dag}}$.
\enddefinition

Another main result of Mochizuki's theory is the following, 
which allows us to recover $\varphi$-group-theoretically the maximal cuspidally 
pro-$l$ extension of $\Pi_X$, in the case where 
the set of cusps consists of a single rational point.

\proclaim {Theorem 3.6} (Reconstruction of Maximal Cuspidally Pro-$l$ Extensions)

\noindent
Let $x_*\in X(k_X)$, $y_*\in Y(k_Y)$, and 
set $S\overset \text {def}\to=\{x_*\}$, $T\overset \text {def}\to=\{y_*\}$,
$U_S\overset \text {def}\to=X\setminus S$, $V_T\overset \text {def}\to=Y\setminus T$. 
Assume that the Frobenius-preserving isomorphism $\alpha:\Pi_X\isom\Pi_Y$ 
maps 
the decomposition group of $x_*$ in $\Pi_X$ (which is well-defined up to conjugation) 
to the decomposition group of $y_*$ in $\Pi_Y$ (which is well-defined up to 
conjugation). Then, for 
each prime $l\in \Sigma ^{\dag}$ (thus, $l\neq p$), 
there exists a commutative diagram:
$$
\CD
\Pi_{U_S}^{\text{\rm c-$l$}} @>\alpha^{\text{\rm c-$l$}}>>  
\Pi_{V_T}^{\text{\rm c-$l$}}\\
@VVV                   @VVV   \\
\Pi_X        @>\alpha>> \Pi_Y \\
\endCD
$$ 
in which $\Pi_{U_S}\twoheadrightarrow  \Pi_{U_S}^{\text{\rm c-$l$}}$, $\Pi_{V_T}
\twoheadrightarrow  \Pi_{V_T}^{\text{\rm c-$l$}}$ are the maximal cuspidally pro-$l$ 
quotients (relative to the maps $\Pi_{U_S}\twoheadrightarrow    \Pi_X$, 
$\Pi_{V_T}\twoheadrightarrow  \Pi_Y$, respectively), the vertical arrows are the 
natural surjections, and $\alpha^{\text{\rm c-$l$}}$
is an isomorphism well-defined up to composition with a cuspidally inner 
automorphism (i.e., an inner automorphism by an element of 
$\Ker (\Pi_{V_T}^{\text{\rm c-$l$}} \to \Pi_Y$)), which is compatible relative to the natural 
surjections
$$\Pi_{U_S}^{\text{\rm c-$l$}}\twoheadrightarrow \Pi_{U_S}^{\text{\rm c-ab},l}
,\ \ \ \ \Pi_{V_T}^{\text{\rm c-$l$}}\twoheadrightarrow \Pi_{V_T}^{\text{\rm c-ab},l}$$
where the subscript ``$\text{\rm c-ab},l$'' denotes the maximal cuspidally pro-$l$ 
abelian quotient, with the isomorphism
$$\alpha_{S,T}^{\text{\rm c-ab}}: 
\Pi_{U_S}^{\text{\rm c-ab}}\isom \Pi _{V_T}^{\text{\rm c-ab}}$$
in Corollary 3.3. Moreover, $\alpha^{\text{\rm c-$l$}}$ is compatible, up to cuspidally
inner automorphisms, with the decomposition groups of $x_*$, $y_*$ in 
$\Pi_{U_S}^{\text{\rm c-$l$}}$, $\Pi_{V_T}^{\text{\rm c-$l$}}$.
\qed
\endproclaim

\demo
{Proof}\ See [Mochizuki1], Theorem 3.1. (See also [Sa\"\i di-Tamagawa1], Theorem 2.6.) 
\qed
\enddemo

If $n$ is 
an integer all of whose prime factors belong to $\Sigma ^{\dag}$, 
then we have the Kummer exact sequence
$$1\to \mu_n\to \Bbb G_m\to \Bbb G_m\to 1,$$
where $\Bbb G_m\to \Bbb G_m$ is the $n$-th power map. We shall identify 
$\mu_n$ with $M_X/nM_X$ according to the identification in [Mochizuki1], 
the discussion at the beginning of \S 2. 

Consider a subset
$$E\subset X^{\cl}.$$
(We will set $E=E_X$ eventually, but $E$ is arbitrary for the present.) 
Let $S\subset X^{\cl}\setminus E$ be a finite set. If we consider the above 
Kummer exact sequence on 
the \'etale site of $U_S \overset \text {def}\to=X\setminus S$ 
and pass to the inverse limit with respect to $n$, then we obtain a  
natural homomorphism
$$\Gamma (U_S,\Cal O_{U_S}^{\times})\to H^1(\Pi _{U_S},M_X)$$
(cf. loc. cit.). (Note that here it suffices to consider the group cohomology of 
$\Pi _{U_S}$ (i.e., as opposed to the \'etale cohomology of $U_S$),
since the extraction of $n$-th roots of an element of $\Gamma (U_S,\Cal O_
{U_S}^{\times})$ yields finite \'etale coverings of $U_S$ that correspond to 
open subgroups of $\Pi _{U_S}$.)  The above homomorphism induces a natural
injective homomorphism
$$\Gamma (U_S,\Cal O_{U_S}^{\times})/(k^{\times}\{\Sigma'\})\to H^1(\Pi _{U_S},M_X)$$
where $k^{\times}\{\Sigma'\}$ stands for the $\Sigma'$-primary part of  the multiplicative group $k^{\times}$
(since the abelian group $\Gamma (U_S,\Cal O_{U_S}
^{\times})/(k^{\times}\{\Sigma'\})$ is finitely generated and free of $\Sigma'$-primary torsion, hence injects into
its pro-$\Sigma$ completion). 
In particular, by allowing $S$ to vary among
all finite subsets of $X^{\cl}\setminus E$, we obtain a natural injective
homomorphism
$$\Cal O_{E}^{\times}/(k^{\times}\{\Sigma '\})\to \underset {S}\to {\varinjlim}\ H^1(\Pi _{U_S},M_X),$$  
where 
$$\Cal O_{E}^{\times}\overset \text {def}\to=\{f\in K_X^{\times}\mid
\sup (\div (f))\cap E=\varnothing\}$$ 
is the multiplicative group of the units
in the
ring $\Cal O_{E}$ of functions on $X$ 
which are regular at all points of $E$. 
(Here, $K_X$ denotes the function field of $X$.)

\proclaim {Proposition 3.7} (Kummer Classes of Functions) 
Suppose that $S\subset X^{\cl}\setminus E$ is a finite subset. 
Write 
$$\Delta _{U_S}\twoheadrightarrow \Delta _{U_S}^{\text{\rm c-ab}}\twoheadrightarrow 
\Delta _{U_S}^{\text{\rm c-cn}}$$
for the maximal cuspidally abelian and 
the maximal cuspidally central quotients,
respectively, relative to the map $\Delta _{U_S}\twoheadrightarrow \Delta _X$,
and
$$\Pi _{U_S}\twoheadrightarrow \Pi _{U_S}^{\text{\rm c-ab}}\twoheadrightarrow 
\Pi _{U_S}^{\text{\rm c-cn}}$$
for the corresponding quotients of $\Pi _{U_S}$ (i.e., $\Pi_{U_{S}}^{\text{\rm c-ab}}
\defeq \Pi_{U_{S}}/\Ker (\Delta _{U_{S}}\twoheadrightarrow 
\Delta_{U_{S}}^{\text{\rm c-ab}})$,
$\Pi_{U_{S}}^{\text{\rm c-cn}}
\defeq \Pi_{U_{S}}/\Ker (\Delta _{U_{S}}\twoheadrightarrow 
\Delta_{U_{S}}^{\text{\rm c-cn}})$). 
If $x\in X^{\cl}$, then we shall write
$$D_x[U_S]\subset \Pi_{U_S}$$
for the decomposition group at $x$ in $\Pi_{U_S}$ (which is well-defined 
up to conjugation), and $I_x[U_S]  \overset \text {def}\to=D_x[U_S]\cap 
\Delta _{U_S}$ for the inertia subgroup of $D_x[U_S]$. Thus, when 
$x\in S$ we have a natural isomorphism of $M_X$ with 
$I_x[U_S]$ (cf. [Mochizuki1], Proposition 1.6(ii)(iii)). 
Then: 

\noindent
{\rm (i)} The natural surjections above induce the following isomorphisms:
$$ H^1(\Pi _{U_S}^{\text{\rm c-cn}},M_X)\isom H^1(\Pi _{U_S}^{\text{\rm c-ab}},M_X)\isom 
H^1(\Pi _{U_S},M_X).$$
In particular, we obtain the following natural injective homomorphisms: 
$$\Gamma (U_S,\Cal O_{U_S}^{\times})/(k^{\times}\{\Sigma'\})\hookrightarrow  H^1(\Pi _{U_S}^{\text{\rm c-cn}},M_X)\isom 
H^1(\Pi _{U_S}^{\text{\rm c-ab}},M_X)\isom
H^1(\Pi _{U_S},M_X),$$
$$\Cal O_{E}^{\times}/(k^{\times}\{\Sigma'\})\hookrightarrow \underset {S}\to {\varinjlim}
\ H^1(\Pi _{U_S}^{\text{\rm c-cn}},M_X)\isom \underset {S}\to 
{\varinjlim}\ H^1(\Pi _{U_S}^{\text{\rm c-ab}},M_X)\isom \underset 
{S}\to {\varinjlim}\ H^1(\Pi _{U_S},M_X),$$
where $S$ varies among all finite subsets of $X\setminus E$.

\noindent
{\rm (ii)}
Restricting 
cohomology classes of $\Pi _{U_S}$ to the various $I_x[U_S]$ for $x\in S$ 
yields a natural exact sequence:
$$1\to (k^{\times})^{\Sigma^{\dag}}\to H^1(\Pi _{U_S},M_X) \to 
(\underset {s\in S}\to \oplus 
\hat \Bbb Z^{\Sigma^{\dag}})$$
(where we identify $\Hom _{\hat \Bbb Z^{\Sigma^{\dag}}}(I_x[U_S],M_X)$ with $\hat 
\Bbb Z^{\Sigma^{\dag}}$, and $(k^{\times})^{\Sigma^{\dag}}$ is the maximal pro-$\Sigma^{\dag}$ quotient of 
the multiplicative group $k^{\times}$). Moreover, the image (via the natural homomorphism given in
{\rm (i)}) of $\Gamma (U_S,\Cal O_{U_S}^{\times})/(k^{\times}\{\Sigma'\})$ in  $H^1(\Pi _{U_S},M_X)$ is 
equal to the inverse image in  $H^1(\Pi _{U_S},M_X)$ of the submodule of
$$(\underset {s\in S}\to \oplus \Bbb Z)\subset (\underset {s\in S}\to \oplus 
\hat \Bbb Z^{\Sigma^{\dag}})$$
determined by the principal divisors (with support in $S$). A similar 
statement holds when $\Pi _{U_S}$ is replaced by $\Pi _{U_S}^{\text{\rm c-cn}}$ 
or $\Pi _{U_S}^{\text{\rm c-ab}}$.

\noindent
{\rm (iii)} 
If $f\in \Gamma (U_S,\Cal O_{U_S}^{\times})$, write $f'$ for its image
in $\Gamma (U_S,\Cal O_{U_S}^{\times})/(k^{\times}\{\Sigma'\})$. Write
$$\kappa _{f'}^{\text{\rm c-cn}}\in H^1(\Pi _{U_S}^{\text{\rm c-cn}},M_X),\ \ 
\kappa _{f'}^{\text{\rm c-ab}}\in H^1(\Pi _{U_S}^{\text{\rm c-ab}},M_X),\ \ 
\kappa _{f'}\in H^1(\Pi _{U_S},M_X)$$
for the associated Kummer classes. If $x\in (X^{\cl}\setminus E)\setminus S$, 
then $D_x[U_S]$ maps, via the natural surjection $\Pi _{U_S}\twoheadrightarrow G_k$, 
isomorphically onto the open 
subgroup $G_{k(x)}\subset G_k$ (where $k(x)$ is the residue field of $X$ at 
$x$). Moreover, the images of the pulled back classes
$$
\align
\kappa _{f'}^{\text{\rm c-cn}}|_{D_x[U_S]}=\kappa _{f'}^{\text{\rm c-ab}}|_{D_x[U_S]}
=\kappa _{f'}|_{D_x[U_S]}
\in H^1(D_x[U_S],M_X)&\simeq  H^1(G_{k(x)},M_X)\\
&\simeq (k(x)^{\times})^{\Sigma^{\dag}}
\endalign
$$
in $(k(x)^{\times})^{\Sigma^{\dag}}$ are equal to 
the image in $(k(x)^{\times})^{\Sigma^{\dag}}$ of the value 
$f(x)\in k(x)^{\times}$ of $f$ at $x$. 
\endproclaim

\demo {Proof} See [Sa\"\i di-Tamagawa1], Proposition 3.1. (See also [Mochizuki1], Proposition 2.1.) 
\qed
\enddemo

\definition{Remark 3.8} (cf. [Mochizuki1], Remark 12.) 
In the situation of 
Proposition 3.7(iii), assume $x\in X(k)$ and $S\subset X(k)$ for simplicity.
If we think of the extension $\Pi _{U_S}^{\text{\rm c-cn}}$ 
of $\Pi_X$ as being given by the extension $\Cal D_S$, where $\Cal D$ is a
fundamental extension of $\Pi _{X\times X}$ (cf. Proposition 3.4(i)), 
then it follows that the 
image of $D_x[U_S]$ in $\Pi _{U_S}^{\text{\rm c-cn}}$ may be thought of 
as the image of $D_x[U_S]$ in $\Cal D_S$. 
This image of $D_x[U_S]$ in $\Cal D_S$
amounts to a section of $\Cal D_S\twoheadrightarrow \Pi_X
\twoheadrightarrow G_k$ lying over the section $s_x:G_k\to \Pi_X$ 
corresponding to the rational point $x$ (which is well-defined up to
conjugation). Since $\Cal D_S$ is defined as a 
certain fiber product, this section is equivalent to a collection 
of sections (regarded as ``cyclotomically outer homomorphisms'', i.e., well-defined
up to composition with an inner automorphism of $\Cal D_{y,x}$ by an element
of $\Ker (\Cal D_{y,x}\twoheadrightarrow G_k)$)
$$\gamma _{y,x}:G_k\to \Cal D_{y,x},$$
where $y$ ranges over all points of $S$. Namely, from this point 
of view, Proposition 3.7(iii) may be regarded as saying that the image
in $(k(x)^{\times})^{\Sigma^{\dag}}
=(k^{\times})^{\Sigma^{\dag}}$ of the value $f(x)$ of the function 
$f\in \Gamma (U_S,\Cal O_{U_S}^{\times})$ at $x\in X(k)$ may be 
computed from its Kummer class, as
soon as one knows the sections $\gamma _{y,x}:G_k\to \Cal D_{y,x}$ 
for $y\in S$. Observe that $\gamma _{y,x}$ depends only on $x$, $y$, and 
not on the choice of $S$.
\enddefinition

\definition {Definition 3.9} (cf. [Mochizuki1], Definition 2.1.)
For $x,y\in X(k)$ with $x\neq y$, we shall 
refer to the above section (regarded as a cyclotomically outer homomorphism)
$$\gamma _{y,x}:G_k\to \Cal D_{y,x}$$
as the Green's trivialization of $\Cal D$ at $(y,x)$. If $D$ is a 
divisor on $X$ supported on $k$-rational points $\neq x$, 
then multiplication of the various  Green's trivializations for 
the points in the support of $D$ yields a section (regarded as a 
cyclotomically outer homomorphism)
$$\gamma _{D,x}:G_k\to \Cal D_{D,x}$$
which we shall refer to as the Green's trivialization of $\Cal D$ at 
$(D,x)$.
\enddefinition

\definition {Definition 3.10} (cf. [Mochizuki1], Definition 2.2.) 
Let the notations and the assumptions as in Corollary 3.3. 

\noindent
{\rm (i)} Write $\Cal D$ (respectively, $\Cal E$) for the  
fundamental extension of $\Pi _{X\times X}$ (respectively, 
$\Pi _{Y\times Y}$) that arises as the 
quotient of $\Pi _{U_{X\times X}}^{\text{\rm c-ab}}$ (respectively, $\Pi _{U_
{Y\times Y}}^{\text{\rm c-ab}}$) by the kernel of the maximal cuspidally central
quotient $\Delta _{U_{X\times X}}^{\text{\rm c-ab}}\twoheadrightarrow 
\Delta _{U_{X\times X}}^{\text{\rm c-cn}}$ (respectively,  $\Delta _{U_{Y\times Y}}
^{\text{\rm c-ab}}\twoheadrightarrow \Delta _{U_{Y\times Y}}^{\text{\rm c-cn}}$) 
(cf. Proposition 3.4(ii)). The isomorphism 
$\alpha ^{\text{\rm c-ab}}$ induces naturally an isomorphism:
$$\alpha ^{\text{\rm c-cn}}:\Cal D\isom \Cal E$$
We shall say that $\alpha$ is $(S,T)$-locally Green-compatible 
outside exceptional sets if, for
every pair of points $(x_1,x_2)\in X(k_X)\times X(k_X)$ corresponding via 
$\phi$ to a pair of points $(y_1,y_2)\in Y(k_Y)
\times Y(k_Y)$, such that $x_1\in (X^{\cl}\setminus E_X)\setminus S$, 
$y_1\in (Y^{\cl}\setminus E_Y)\setminus T$, $x_2\in S$, $y_2\in T$, the isomorphism
$$\Cal D_{x_1,x_2}\isom \Cal E_{y_1,y_2}$$
(obtained by restricting $\alpha ^{\text{\rm c-cn}}$ to the various decomposition groups) 
is compatible with the 
Green's trivializations. We shall say that $\alpha$ is 
$(S,T)$-locally bi-principally
Green-compatible outside exceptional sets 
if, for every point $x\in X(k_X)\cap S$ and 
every 
principal divisor
$P$ supported on $k_X$-rational points $\neq x$ 
contained in $X^{\cl}\setminus E_X$ 
corresponding via $\phi$ to a pair $(y,Q)$ (so $y\in Y(k_Y)\cap T$) 
with $Q$ principal, the isomorphism
$$\Cal D_{P,x}\isom \Cal E_{Q,y}$$
obtained from $\alpha^{\text{\rm c-cn}}$
is compatible with the Green's trivializations.

\noindent
{\rm (ii)} We shall say that $\alpha$ is totally globally 
Green-compatible (respectively, 
totally globally bi-principally Green-compatible) 
outside exceptional sets 
if, for all pair of connected finite \'etale coverings $\xi: X'\to X$, 
$\eta: Y'\to Y$ that arise from open subgroups $\Pi_{X'}\subset \Pi_X$, 
$\Pi _{Y'}\subset \Pi_Y$, 
corresponding to each other via $\alpha$, then for any subset $S\subset X^{\cl}\setminus E_X$ 
that corresponds, via $\phi$, to $T\subset Y\setminus E_Y$ the isomorphism
$$\Pi_{X'}\isom \Pi_{Y'}$$
induced by $\alpha$ is $(S',T')$-locally Green-compatible (respectively, 
$(S',T')$-locally bi-principally Green-compatible) 
outside exceptional sets, 
where $S'\defeq\xi ^{-1}(S)\subset X^{'\cl}$, 
$T'\defeq\eta^{-1} (T)  \subset Y^{'\cl}$ 
are the inverse images of $S$, $T$, respectively.
\enddefinition

\definition{Remark 3.11}
In [Sa\"\i di-Tamagawa1], Definition 3.4, we adopted a slightly different notion of 
being ``$(S,T)$-locally (or totally globally) principally Green-compatible outside exceptional sets'', 
where the divisor $Q$ of $Y$ appearing in (i) above is not assumed to be principal. 
Here we adopt the above notion of 
being ``$(S,T)$-locally (or totally globally) bi-principally Green-compatible 
outside exceptional sets'', since it is more natural in our settings (although both notions work). 
\enddefinition

\proclaim{Proposition 3.12} (Total Global Green-Compatibility Outside Exceptional Sets)
In the situation of Theorem 3.2, assume further that 
$\alpha$ is Frobenius-preserving. 
Then the isomorphism $\alpha$ is totally globally, and totally globally bi-principally,  Green-compatible 
outside exceptional sets. 
\endproclaim

\demo
{Proof} Similar to the proof of Proposition 3.8 in [Sa\"\i di-Tamagawa1]. 
(See also [Mochizuki1], Corollary 3.1.) 
\qed
\enddemo

\subhead
\S 4. Isomorphisms between geometrically pro-$\Sigma$ arithmetic fundamental groups
\endsubhead

We maintain the notations of \S 3. 

\definition {Definition/Remark 4.1} 
Let $J=J_X$ be the Jacobian variety of $X$. 
Let $\Div _{X\setminus  E}^0$ be the group of degree zero divisors 
on $X$ which are supported on points in $X\setminus E$. 
Write $D_{X\setminus  E}$ for the kernel of the natural homomorphism 
$\Div _{X\setminus  E}^0\to J(k)^{\Sigma}$. 
Here, $J(k)^{\Sigma}$ stands for the maximal pro-$\Sigma$ quotient 
$J(k)/(J(k)\{\Sigma'\})$ of $J(k)$, where, 
for an abelian group $M$, $M\{\Sigma'\}$ stands 
for the subgroup of torsion elements $a$ of $M$ such 
that every prime divisor of the order of $a$ 
belongs to $\Sigma'$. 
Then $D_{X\setminus  E}$ sits naturally in the following exact sequence:
$$0\to \Pri  _{X\setminus  E}\to D_{X\setminus  E}\to J(k)\{\Sigma'\} \to 0,$$
where $\Pri _{X\setminus  E}\defeq\Cal O_{E}^{\times}/k_X^{\times}$ 
stands for the group of principal divisors supported in $X\setminus E$.
Further, let $\Cal D_{X\setminus E}$ be the inverse image of $D_{X\setminus E}$ in 
${\varinjlim}\ H^1(\Pi _{U_S}^{\text{\rm c-ab}},M_X)$ (cf. Proposition 3.7(ii)). 
Then
$\Cal D_{X\setminus E}$ sits naturally in the following exact sequence
$$0\to \Cal O_{E}^{\times}/(k_X^{\times}\{\Sigma'\})  \to \Cal D_{X\setminus  E}\to J(k)\{\Sigma'\} \to 0.$$
\enddefinition

Now, let $X$, $Y$ be proper 
hyperbolic curves over finite fields $k_X$, $k_Y$ of characteristic
$p_X$, $p_Y$, respectively, and define $K_X$, $K_Y$ to be the 
function fields of $X$, $Y$, respectively. 
Let $\Sigma_X$ (respectively, $\Sigma_Y$)
be a set of prime numbers that contains at least one prime number
different from $p_X$ (respectively, $p_Y$).
Write $\Delta _X$ (respectively, $\Delta _Y$) 
for the maximal pro-$\Sigma _X$ quotient of $\pi _1(\overline X)$ 
(respectively, the maximal  pro-$\Sigma_Y$ quotient of $\pi _1
(\overline Y)$), and $\Pi _X$ (respectively, $\Pi _Y$) for the 
quotient $\pi_1(X)/\Ker (\pi_1(\overline X)\twoheadrightarrow \Delta_X)$
of $\pi _1(X)$ (respectively, the quotient 
$\pi_1(Y)/\Ker (\pi_1(\overline Y)\twoheadrightarrow \Delta_Y)$ of $\pi _1(Y)$). 

Let 
$$\alpha :\Pi _X\isom \Pi _Y$$
be an isomorphism of profinite groups and write 
$\Sigma\defeq\Sigma_X=\Sigma_Y$ (cf. Proposition 1.7(ii)). 

\proclaim {Theorem 4.2} (Reconstruction of Pseudo-Functions) 
Assume that $\alpha$ is Frobenius-preserving. 
Then:

\noindent
{\rm (i)} The bijection $\phi:X^{\cl}\setminus E_X\isom Y^{\cl}\setminus E_Y$
induced by $\alpha$ (where $E_X$ and $E_Y$ are the exceptional sets), 
together with the isomorphisms  in Corollary 3.3,
induce natural bijections 
$\bar \rho : D_{Y\setminus  E_Y}\isom D_{X\setminus  E_X},$
$\rho':\Cal D_{Y\setminus  E_Y}\isom \Cal D_{X\setminus  E_X},$
which fit into the following 
commutative diagrams
$$
\CD
0 @>>> \Pri  _{X\setminus  E_X}=\Cal O_{E_X}^{\times}/k_X^{\times} @>>> D_{X\setminus  E_X}@>>> J_X(k_X)\{\Sigma'\} @>>> 0 \\
@.    @.    @A{\bar \rho}AA   @.\\
0@>>> \Pri  _{Y\setminus  E_Y}=\Cal O_{E_Y}^{\times}/k_Y^{\times} @>>> D_{Y\setminus  E_Y}@>>> J_Y(k_Y)\{\Sigma'\} @>>> 0 \\
\endCD
$$
and
$$
\CD
0 @>>> \Cal O_{E_X}^{\times}/(k_X^{\times}\{\Sigma'\}) @>>> \Cal D_{X\setminus  E_X}@>>> J_X(k_X)\{\Sigma'\} @>>> 0 \\
@.    @.    @A{\rho'}AA   @.\\
0@>>> \Cal O_{E_Y}^{\times}/(k_Y^{\times}\{\Sigma'\}) @>>> \Cal D_{Y\setminus  E_Y}@>>> J_Y(k_Y)\{\Sigma'\} @>>> 0 \\
\endCD
$$
Moreover, the following diagram commutes
$$
\CD
\Cal D_{X\setminus  E_X} @<{\rho'}<<  \Cal D_{Y\setminus  E_Y} \\
@VVV           @VVV  \\
D_{X\setminus  E_X} @<{\bar \rho}<<  D_{Y\setminus  E_Y} \\
\endCD
$$
where the vertical maps are the natural ones.

\noindent
{\rm (ii)} The bijection $\bar \rho$ in (i) induces a natural isomorphism
$$\bar \rho:\overline H_Y\isom \overline H_X.$$
Here, 
$\overline H_X\defeq 
\Ker(\Pri_{X\setminus E_X} \overset{\bar \rho^{-1}}\to{\to} D_{Y\setminus E_Y} \to J_Y(k_Y)\{\Sigma'\})$ 
and 
$\overline H_Y\defeq 
\Ker(\Pri_{Y\setminus E_Y} \overset{\bar \rho}\to{\to} D_{X\setminus E_X} \to J_X(k_X)\{\Sigma'\})$ 
are 
finite index subgroups of $\Cal O_{E_X}^{\times}/k_X^{\times}$ 
and $\Cal O_{E_Y}^{\times}/k_Y^{\times}$, respectively. 
(More precisely, the quotients 
$(\Cal O_{E_X}^{\times}/k_X^{\times})/\overline H_X$ and 
$(\Cal O_{E_Y}^{\times}/k_Y^{\times})/\overline H_Y$ 
are $\Sigma'$-primary finite abelian groups that are embeddable into 
$J_Y(k_Y)\{\Sigma'\}$ and 
$J_X(k_X)\{\Sigma'\}$, respectively.) 
Moreover, let $H_X'$ (resp. $H_Y'$) be the inverse image of $\overline H_X$ (resp. $\overline H_Y$)
in $\Cal O_{E_X}^{\times}/(k_X^{\times})\{\Sigma'\}$ (resp. $\Cal O_{E_Y}^{\times}/(k_Y^{\times}\{\Sigma '\})$). 
Then the isomorphism $\rho'$ in (i)
induces a natural isomorphism $\rho': H_Y'\isom H_X'$ which fits into the following commutative diagram
$$
\CD
H_X' @<\rho'<<   H_Y' \\
@VVV      @VVV   \\
\overline H_X  @<\bar \rho<< \overline H_Y   \\
\endCD
\tag 4.1
$$
where the vertical maps are the natural surjections. 

In particular, $\rho'$ induces a natural isomorphism 
$$\tau: (k_Y^{\times})^{\Sigma}=k_Y^{\times}/(k_Y^{\times}\{\Sigma'\})\isom 
k_X^{\times}/(k_X^{\times}\{\Sigma'\})=(k_X^{\times})^{\Sigma}.$$

\noindent
{\rm (iii)} The diagram in (ii) is functorial  in $X$, $Y$, in the following sense: if $\xi :X'\to X$ is 
a finite \'etale covering, arising from an open subgroup $\Pi_{X'}\subset \Pi_X$, which corresponds 
to a finite \'etale covering $Y'\to Y$ via $\alpha$ (thus, $\Pi_{Y'}\defeq \alpha (\Pi _{X'})$), then $\alpha$ induces
natural isomorphisms 
${\gamma}' :H_{Y'}'\isom  H_{X'}' $ and $\overline { \gamma}:\overline H_{Y'}\isom \overline H_{X'}$,
which fit into the following commutative diagrams 

$$
\CD
H_{X'}' @<{\gamma}'<<   H_{Y'}' \\
@AAA      @AAA  \\
H_X'  @<\rho'<< H_Y'   \\
\endCD
$$
and 
$$
\CD
\overline H_{X'} @<\overline {\gamma}<<  \overline  H_{Y'} \\
@AAA     @AAA  \\
\overline H_X  @<\overline \rho<< \overline H_Y   \\
\endCD
$$
where the vertical maps are natural homomorphisms induced by the natural injective homomorphisms
$\Cal O_{E_X}^{\times}/k_X^{\times} \{\Sigma '\}\hookrightarrow \Cal O_{E_{X'}}^{\times}/k_{X'}^{\times} \{\Sigma'\}$,
$\Cal O_{E_Y}^{\times}/k_Y^{\times} \{\Sigma'\}\hookrightarrow \Cal O_{E_{Y'}}^{\times}/k_{Y'}^{\times} \{\Sigma'\}$.
$\Cal O_{E_X}^{\times}/k_X^{\times}\hookrightarrow \Cal O_{E_{X'}}^{\times}/k_{X'}^{\times}$, and
$\Cal O_{E_Y}^{\times}/k_Y^{\times}\hookrightarrow \Cal O_{E_{Y'}}^{\times}/k_{Y'}^{\times}$. 

In particular, 
$\alpha$ induces a natural isomorphism 
$$\tau: (\bar k_Y^{\times})^{\Sigma}=\bar k_Y^{\times}/(\bar k_Y^{\times}\{\Sigma'\})\isom 
\bar k_X^{\times}/(\bar k_X^{\times}\{\Sigma'\})=(\bar k_X^{\times})^{\Sigma}$$
which extends $\tau:
(k_Y^{\times})^{\Sigma}\isom (k_X^{\times})^{\Sigma}$ in (ii). 
\endproclaim

\demo {Proof}
Similar to the proof of [Sa\"\i di-Tamagawa1], Theorem 3.6. (See also [Sa\"\i di-Tamagawa3], Lemma 4.5 
for a similar statement in the birational case.) Here, the commutativity of the diagrams in (iii) follows 
basically from the functoriality of Kummer theory, together with the commutativity of the diagram 
$$
\CD
M_{X'} @<<<  M_{Y'} \\
@AAA     @AAA  \\
M_X  @<<< M_Y   \\
\endCD
$$
where the horizontal arrows are isomorphisms induced by $\alpha$ 
(or, more precisely, by $\alpha^{-1}$, via the functoriality of $H^2$) 
and the vertical arrows are 
isomorphisms defined geometrically via the identification of $M_X$ and $M_{X'}$ (resp. $M_Y$ and $M_{Y'}$) 
with the Tate module of $\Bbb G_m$ over $\bar k_X=\bar k_{X'}$ (resp. 
$\bar k_Y=\bar k_{Y'}$). The commutativity of this last diagram follows from the fact that 
the above (geometrically defined) isomorphism $M_X\isom M_{X'}$ 
(resp. $M_Y\isom M_{Y'}$) is identified with the composite of 
the $(X'_{\bar k_X}:X_{\bar k_X})$-multiplication map $M_X\to M_X$ 
(resp. the $(Y'_{\bar k_Y}:Y_{\bar k_Y})$-multiplication map $M_Y\to M_Y$) 
and the inverse of the natural isomorphism 
$M_{X'}\isom (X'_{\bar k_X}:X_{\bar k_X})M_X\subset M_X$ (resp. 
$M_{Y'}\isom (Y'_{\bar k_Y}:Y_{\bar k_Y})M_Y\subset M_Y$) induced by 
the inclusion $M_{X'}\to M_X$ (resp. 
$M_{Y'}\to M_Y$) 
arising from the functoriality of $H^2$, 
and the fact that 
$(X'_{\bar k_X}:X_{\bar k_X})
=(\Delta_{X}:\Delta_{X'})
=(\Delta_{Y}:\Delta_{Y'})
=(Y'_{\bar k_Y}:Y_{\bar k_Y})$. 
\qed
\enddemo

\definition{Definition 4.3}
We say that $\alpha$ is pseudo-constants-additive if $\tau: (\bar k_Y^{\times})^{\Sigma}\isom 
(\bar k_X^{\times})^{\Sigma}$ in Theorem 4.2(iii) satisfies the following: 
For $\eta\in \overline k_{Y}^{\times}$ and $\zeta\in \overline k_{X}^{\times}$, 
if 
$$\text{$1+\eta\neq0$ and $\tau (\eta')=\zeta'$},$$  
then there exist $\alpha, \beta \in \overline k_{X}^{\times}\{\Sigma'\}$, such that 
$$\text{$\alpha+\beta\zeta \neq 0$ and $\tau ((1+\eta)')=(\alpha+\beta \zeta)'$}.$$ 
Here, for an element $\xi$ of $\bar k_X^{\times}$ (resp. $\bar k_Y^{\times}$), 
$\xi'$ denotes its image in 
$(\bar k_X^{\times})^{\Sigma}$ (resp. $(\bar k_Y^{\times})^{\Sigma}$). 
\enddefinition


\proclaim
{Theorem 4.4} (Pseudo-Constants-Additive Isomorphisms)
Assume 
that $\Sigma$ is $k_X$-large and $k_Y$-large, 
that at least one of $X$ and $Y$ is almost $\Sigma$-separated, and 
that $\alpha$ is pseudo-constants-additive. 
Then $\alpha$ arises from a uniquely 
determined commutative diagram of schemes: 
$$
\CD
\Tilde X @>{\sim}>> \Tilde Y \\
@VVV   @VVV \\
X @>{\sim}>> Y \\
\endCD
$$
in which the horizontal arrows are isomorphisms and the vertical arrows are the 
profinite \'etale coverings determined by the groups $\Pi_X$, $\Pi_Y$.
\endproclaim

The rest of this section is devoted to the proof of Theorem 4.4 and deducing corollaries.

First, since $\Sigma$ is $k_X$-large and $k_Y$-large, $\alpha$ is Frobenius-preserving 
by Proposition 2.5(ii). Thus, we may apply Theorem 1.9, Theorem 3.12 and Theorem 4.2 
to $\alpha$. 
In particular, $\alpha$ is totally globally bi-principally Green-compatible 
outside exceptional sets by Theorem 3.12. 
Also, we see that both $X$ and $Y$ are almost 
$\Sigma$-separated. Indeed, let $F_X$ (resp. $F_Y$) denote the compositum of 
$k(x)$ (resp. $k(y)$) in $\bar k_X$ (resp. $\bar k_Y$) for all $x\in E_X$ 
(resp. $y\in E_Y$). Then it follows from Theorem 1.9(i) that 
$G_{F_X}\subset G_{k_X}$ corresponds to $G_{F_Y}\subset G_{k_Y}$ 
via the natural isomorphism $G_{k_X}\isom G_{k_Y}$ induced by $\alpha$ 
(cf. Proposition 1.7(i)). Thus, $F_X\subsetneq \bar k_X$ (i.e., $G_{F_X}\neq\{1\}$) 
is equivalent to $F_Y\subsetneq \bar k_Y$ (i.e., $G_{F_Y}\neq\{1\}$). 

Let $x\in X^{\cl}\setminus E_X$ and $y\defeq \phi(x)\in Y^{\cl}\setminus E_Y$. 
Then,  
as $\alpha$ preserves the decomposition groups at $x$, $y$, respectively,  
$\alpha$ induces naturally an isomorphism
$$\tau _{x,y}:(k(y)^{\times})^{\Sigma}\isom  (k(x)^{\times})^{\Sigma}$$ 
(cf. Proposition 3.7(iii))), 
which fits into a commutative diagram
$$
\CD
(k(x)^{\times})^{\Sigma} @<\tau_{x,y}<< (k(y)^{\times})^{\Sigma} \\
@AAA             @AAA   \\
(k_X^{\times})^{\Sigma}  @<\tau<<   (k_Y^{\times})^{\Sigma}   \\
\endCD
$$
where the vertical maps are the natural homomorphisms.

Next, we shall 
think of elements of $\Cal O_{E_{X}}^{\times}/(k_{X})^{\times}$, $\Cal O_{E_{Y}}^{\times}/(k_{Y})^{\times}$
as principal divisors of rational functions on $X$, $Y$, respectively, and denote them 
$\bar f, \bar g, \dots$, where $f$, $g$ are rational functions on $X$, $Y$, 
whose supports of divisors are disjoint from $E_{X}$, $E_{Y}$, 
respectively.  We will denote the elements of 
$\Cal O_{E_{X}}^{\times}/(k_{X}^{\times}\{\Sigma'\})$,  $\Cal O_{E_{Y}}^{\times}/
(k_{Y}^{\times}\{\Sigma '\})$, by $f', g',\dots$, 
where $f$, $g$, are rational functions on $X$, $Y$, whose supports of divisors 
are disjoint from $E_{X}$, $E_{Y}$, respectively,
and refer to them as ``pseudo-functions''$\defeq$classes of elements of $\Cal O_{E_{X}}^{\times}$,  
$\Cal O_{E_{Y}}^{\times}$,
modulo $k_{X}^{\times}\{\Sigma'\}$,  $k_{Y}^{\times}\{\Sigma '\}$, respectively. 
For each $x\in X^{\cl}\setminus E_X$ (resp. $y\defeq \phi(x)\in Y^{\cl}\setminus E_Y$), 
we denote by 
$v_x: \Cal O_{E_{X}}^{\times}/(k_{X})^{\times}\to \Bbb Z$ 
(resp. $v_y: \Cal O_{E_{Y}}^{\times}/(k_{Y})^{\times}\to \Bbb Z$) 
the function induced by the (normalized, additive) valuation 
$v_x: K_X^{\times} \to \Bbb Z$ at $x$ 
(resp. $v_y: K_Y^{\times} \to \Bbb Z$ at $y$). Further, 
we denote by 
$\deg: \Cal O_{E_{X}}^{\times}/(k_{X})^{\times}\to \Bbb Z_{\geq 0}$ 
(resp. $\deg: \Cal O_{E_{Y}}^{\times}/(k_{Y})^{\times}\to \Bbb Z_{\geq 0}$) 
the function induced by the degree function 
$\deg: K_X^{\times} \to\Bbb Z_{\geq 0}$ 
(resp. $\deg: K_Y^{\times}\to \Bbb Z_{\geq 0}$) 
that sends $f\in K_X^{\times}$ (resp. $f\in K_Y^{\times}$) 
to the degree of the pole divisor of $f$. 

\proclaim {Lemma 4.5}(Recovering the Valuations and the $\Sigma$-Values of Pseudo-Functions) 
Consider the commutative diagram (4.1). 
Let $x\in X^{\cl}\setminus E_X$ and $y\defeq \phi(x)\in Y^{\cl}\setminus E_Y$. 
Then the following implications hold:

\noindent 
{\rm (i)} For $\bar f\in \overline H_{Y}$ and $\bar g\in\overline H_{X}$:
$$\bar \rho (\bar f)=\bar g \Longrightarrow v_{y}(\bar f)=v_{x}(\bar g).$$ 
In particular, in terms of divisors, if: 
$$\bar f=y_1+y_2+\dots+y_n-y_1'-\dots-y_{n'}',$$ 
then: 
$$\bar g= x_1+x_2+\dots+x_n-x_1'-\dots-x_{n'}',$$ 
where $y_i\defeq \phi (x_i)$ 
(resp. $y_{i'}'\defeq \phi (x_{i'}')$) for $i\in \{1,\dots,n\}$ 
(resp. $i'\in \{1,\dots,n'\}$). In other words, 
the isomorphism $\bar \rho: \overline H_Y\isom \overline H_X$ 
preserves the valuations of the classes of functions in $\overline H_X, \overline H_Y$
with respect to the bijection $\phi:X^{\cl}\setminus E_X\isom Y^{\cl}\setminus E_Y$ between points.
Further, 
the isomorphism 
$\bar \rho$ preserves the degrees of the classes of functions in $\overline H_X, \overline H_Y$. 

\noindent
{\rm (ii)} For $f'\in H_{Y}'$ and $g'\in H_{X}'$:
$$\text{$v_{y}(\bar f)=0$ and $\rho(f')=g'$}
\Longrightarrow \text{$v_{x}(\bar g)=0$ and $\tau _{x,y} (f'(y))=g'(x)$}$$
where 
$$y=\phi (x)\ and\ \tau _{x,y}:(k(y)^{\times})^{\Sigma} \isom (k(x)^{\times})^{\Sigma}$$ 
is the natural identification above. In other words, 
the isomorphism 
$\rho': H_Y'\isom H_X'$ preserves the $\Sigma$-values of the pseudo-functions 
in $H_X', H_Y'$
with respect to the bijection $\phi:X^{\cl}\setminus E_X\isom Y^{\cl}\setminus E_Y$ between points.
\endproclaim

\demo{Proof} 
Assertion (i) follows from Proposition 3.7(i)(ii), together with 
the fact that $[k(x):k_{X}]=[k(\phi(x)):k_Y]$ for each $x\in X^{\cl}\setminus E_X$ 
(cf. discussion before Corollary 3.3), and assertion (ii) follows from
Proposition 3.7(iii), together with the fact that $\alpha$ is totally globally 
bi-principally Green-compatible outside exceptional sets.
\qed
\enddemo

\proclaim
{Lemma 4.6} 
Let $x\in X^{\cl}\setminus E_X$ and $y\defeq \phi(x)\in Y^{\cl}\setminus E_Y$. 
The natural identification 
$$\tau _{x,y}:(k(y)^{\times})^{\Sigma} \isom (k(x)^{\times})^{\Sigma},$$
induced by $\alpha$, satisfies the following property: 
For $\eta\in k(y)^{\times}$ and $\zeta\in k(x)^{\times}$, if 
$$\text{$1+\eta\neq0$ and $\tau_{x,y} (\eta')=\zeta'$},$$  
then there exist $\alpha, \beta \in \overline k_X^{\times}\{\Sigma'\}$, such that 
$$\text{$\alpha+\beta\zeta \neq 0$ and $\tau_{x,y} ((1+\eta)')=(\alpha+\beta \zeta)'$}.$$ 
\endproclaim

\demo{Proof}
Similar to the proof of [Sa\"\i di-Tamagawa3], Lemma 4.10. 
\qed
\enddemo

Let $l\in \Primes$ be a prime number which is both 
$(X,\Sigma)$-admissible and $(Y,\Sigma)$-admissible, 
i.e.,
for every  finite extension $k'$ of $k_X$ (resp. $k_Y$),  
$(k')^{l}\cap F_X$ (resp. $(k')^{l}\cap F_Y$) is finite, 
where $(k')^{l}$ is the maximal pro-$l$ extension of $k'$, and, in particular, 
$E_X\cap 
X((k')^{l})^{\cl}
$ (resp. $E_Y\cap 
Y((k')^{l})^{\cl}
$) is finite. 
Let $X^l$, $Y^l$ be the normalization of $X$, $Y$ in $K_X k_X^{l}$, $K_Yk_Y^{l}$, respectively.
Set $E_{X^l}\defeq E_X\times _{k_X} k_X^{l}$, $E_{Y^l}\defeq E_Y\times _{k_Y}k_Y^{l}$, 
and write 
$\Cal O_{E_{X^l}}^{\times}$,  $\Cal O_{E_{Y^l}}^{\times}$ for the group of multiplicative functions on $X^l$, $Y^l$
whose supports of divisors are disjoint from $E_{X^l}$, $E_{Y^l}$, respectively.
Define $D_{X^l\setminus E_{X^l}}$, $\Cal D_{X^l\setminus E_{X^l}}$ (resp. $D_{Y^l\setminus E_{Y^l}}$,
$\Cal D_{Y^l\setminus E_{Y^l}}$) in a similar way as in Definition 4.1. Thus, we have natural exact sequences
$$0\to \Cal O_{E_{X^l}}^{\times}/(k_X^{l})^{\times}  \to 
D_{X^l\setminus  E_{X^l}}\to J_X(k_X^{l})\{\Sigma'\} \to 0$$
and 
$$0\to \Cal O_{E_{X^l}}^{\times}/(k_X^{l})^{\times}\{\Sigma'\})  \to 
\Cal D_{X^l\setminus  E_{X^l}}\to J_X(k_X^{l})\{\Sigma'\} \to 0$$ 
(resp. similar sequences for $D_{Y^l\setminus E_{Y^l}}$ and $\Cal D_{Y^l\setminus  E_{Y^l}}$). 

The isomorphism $\alpha:\Pi_X\isom \Pi_Y$ induces natural isomorphisms
$\bar \rho:D_{Y^l\setminus E_{Y^l}}\isom D_{X^l\setminus E_{X^l}}$, 
and $\rho':\Cal D_{Y^l\setminus E_{Y^l}}\isom \Cal D_{X^l\setminus E_{X^l}}$ 
(by passing to finite subextensions of 
$k_X^{l}/k_X$ and $k_Y^{l}/k_Y$ corresponding to each other by $\alpha$,
cf. Theorem 4.2(i)), which fit into the following commutative diagram
$$
\CD
\Cal D_{X^l\setminus E_{X^l}} @<\rho'<< \Cal D_{Y^l\setminus E_{Y^l}} \\
@VVV                 @VVV \\
D_{X^l\setminus E_{X^l}} @<\bar \rho<< D_{Y^l\setminus E_{Y^l}} \\
\endCD
$$
where the vertical maps are the natural ones. 
The above isomorphism $\bar \rho$ implies the existence
of 
subgroups $\overline H_{X^l}$
(resp. $\overline H_{Y^l}$) of $\Cal O_{E_{X^l}}^{\times}/(k_X^{l})^{\times}$ 
(resp. $\Cal O_{E_{Y^l}}^{\times}/(k_Y^{l})^{\times}$),
which are functorial in $X$ and $Y$, 
and a natural isomorphism
$$\bar \rho:\overline H_{Y^l}\isom \overline H_{X^l},$$
which lies above the isomorphism $\bar \rho:\overline H_Y\isom \overline H_X$
in Theorem 4.2(ii). 
(Use similar arguments as those in the proof of Theorem 4.2(ii).) 
Moreover, let $H_{X^l}'$ (resp. $H_{Y^l}'$) be the inverse image of $\overline H_{X^l}$ (resp. $\overline H_{Y^l}$)
in $\Cal O_{E_{X^l}}^{\times}/((k_X^{l})^{\times}\{\Sigma'\})$ 
(resp. $\Cal O_{E_{Y^l}}^{\times}/((k_Y^{l})^{\times}\{\Sigma'\})$). Then there exists a 
natural isomorphism $\rho': H_{Y^l}'\isom H_{X^l}'$ (cf. Theorem 4.2(ii)), and $H_{X^l}'$ (resp. $H_{Y^l}' $) is 
a 
subgroup of 
$\Cal O_{E_{X^l}}^{\times}/((k_X^{l})^{\times}\{\Sigma'\})$ 
(resp. $\Cal O_{E_{Y^l}}^{\times}/((k_Y^{l})^{\times}\{\Sigma'\})$).
Further, let
$$H_{X^l}^{\times} \defeq \{f\in \Cal O_{E_{X^l}}^{\times}\  \vert \  \bar f\in \overline H_{X^l}\}\ 
\text{(resp. 
$H_{Y^l}^{\times} \defeq \{f\in \Cal O_{E_{Y^l}}^{\times}\  \vert \  \bar f\in \overline H_{Y^l}\}$
)}.$$
Then  $H_{X^l}^{\times}$ (resp. $H_{Y^l}^{\times}$) is a 
subgroup of $\Cal O_{E_{X^l}}^{\times}$ (resp. $\Cal O_{E_{Y^l}}^{\times}$), and the 
quotient
$\Cal O_{E_{X^l}}^{\times}/H_{X^l}^{\times}$ 
(resp. $\Cal O_{E_{Y^l}}^{\times}/H_{Y^l}^{\times}$) 
is embeddable into $J_Y(k_Y^{l})\{\Sigma'\}$ 
(resp. $J_X(k_X^{l})\{\Sigma'\}$), hence 
is a $\Sigma'$-primary torsion abelian group. 
(More precisely, let $r$ be a prime number. Then the $r$-primary part 
of $\Cal O_{E_{X^l}}^{\times}/H_{X^l}^{\times}$ 
(resp. $\Cal O_{E_{Y^l}}^{\times}/H_{Y^l}^{\times}$) 
is trivial (resp. finite, resp. embeddable into 
a finite direct sum of $\Bbb Q_r/\Bbb Z_r$), if $r\in\Sigma$ (resp. 
$r\in\Sigma'\setminus\{l\}$, resp. $r$ is any prime.) 
We have a commutative diagram
$$
\CD
H_{X^l}' @<\rho'<<   H_{Y^l}' \\
@VVV      @VVV   \\
\overline H_{X^l}  @<\bar \rho<< \overline H_{Y^l}   \\
\endCD
$$
where the vertical maps are the natural surjections (cf. Theorem 4.2(ii)), and lies
above the commutative diagram (4.1) 
in Theorem 4.2(ii). 

We state the following results (Remark 4.7, Lemma 4.8, Corollary 4.10 and Lemmas 4.11-4.13) 
only for $X$, but similar statments also hold for $Y$. 

\definition{Remark 4.7}
If $\Sigma$ is $J_Y$-large, then $F_Y$ is contained in a finite extension ($\neq\bar k_Y$) of 
$\bar k_Y^{\Ker(\rho_{J_Y,\Sigma'})}$ by Proposition 2.12. Thus, any $l$ such that 
the image of the pro-$l$ Sylow subgroup $G_{k_Y,l}$ of $G_{k_Y}$ under 
$\rho_{J_Y,\Sigma'}$ is finite 
is $(Y,\Sigma)$-admissible. 
(Further, by Proposition 2.13(ii), such an $l$ is automatically in $\Sigma\cup\{p\}$.) 
Let us take such an $l$. Then 
it follows that 
$(\rho_{J_Y,\Sigma'}(G_{k_Y}):\rho_{J_Y,\Sigma'}(G_{k_Y^{l}}))<\infty$, hence 
$\sharp(J_Y(k_Y^{l})\{\Sigma'\})<\infty$. Thus, in this case, 
the 
quotient
$\Cal O_{E_{X^l}}^{\times}/H_{X^l}^{\times}$ is finite. 
\enddefinition

\proclaim
{Lemma 4.8} 
Let $f_1,\dots,f_n\in\Cal O_{E_{X^l}}$. 
Then, for all but finitely many $c\in k_X^{l}$, 
$f_1- c,\dots,f_n-c\in \Cal O_{E_{X^l}}^{\times}$. 
\endproclaim

\demo{Proof}
We may and shall assume $n=1$ and set $f\defeq f_1$. (The assertion for general $n$ can be 
reduced to this special case immediately.) The function $f\in 
\Cal O_{E_{X^l}}$ descends to a finite extension of $k_X$. More precisely, there exists a 
finite subextension $k_0/k_X$ of $k_X^{l}/k_X$, such that 
$f\in \Cal O_{E_{X_{k_0}}}$, where 
$E_{X_{k_0}}\defeq E_X\times_{k_X}k_0$ and 
$\Cal O_{E_{X_{k_0}}}= \Cal O_{E_X}\otimes_{k_X}k_0$. Then $f$ defines a $k_0$-morphism 
$f: X_{k_0}\to\Bbb P^1_{k_0}$ of degree, say, $d$. One has 
$f^{-1}(
\Bbb P^1_{k_0}(k^{l})^{\cl}
) \subset 
X(k')^{\cl}
$, 
where $k'$ is the compositum of finite extensions 
of $k^{l}$ of degree $\leq d$, which is finite over $k^{l}$. 
On the other hand, $E_{X_{k_0}}\subset 
X_{k_0}(F)^{\cl}
$, where 
$F\defeq F_Xk_0$. By Proposition 2.13, we see that 
$S:=f^{-1}(
\Bbb P^1_{k_0}(k^{l})^{\cl}
)\cap 
E_{X_{k_0}}$ is finite. Now, for each $c\in k^{l}\setminus f(S)(k^{l})
\ (\subset \Bbb P^1_{k_0}(k^{l}))$, one has $f-c\in \Cal O_{E_{X^l}}^{\times}$, as desired. 
\qed\enddemo

\definition{Definition 4.9}
(i) For a subset $S$ of a field of characteristic $p>0$, we denote by $\langle S \rangle$ the 
abelian subgroup (or, equivalently, $\Bbb F_p$-vector subspace) generated by $S$. 

\noindent
(ii) $R_{X^l}\defeq \langle H_{X^l}^{\times}\rangle$, $R_{Y^l}\defeq \langle H_{Y^l}^{\times}\rangle$. 
\enddefinition

\proclaim{Corollary 4.10}
{\rm (i)} $(k_X^{l})^{\times}+\Cal O_{E_{X^l}}^{\times}=\Cal O_{E_{X^l}}$. 

\noindent
{\rm (ii)} $\langle\Cal O_{E_{X^l}}^{\times}\rangle=\Cal O_{E_{X^l}}$. In particular, 
$\langle\Cal O_{E_{X^l}}^{\times}\rangle$ coincides with the $k_X^{l}$-subalgebra 
generated by $\Cal O_{E_{X^l}}^{\times}$. 
\endproclaim

\demo{Proof}
(i) Clearly $k_X^{l}+\Cal O_{E_{X^l}}^{\times}\subset\Cal O_{E_{X^l}}$. 
To show the opposite, take any $f\in \Cal O_{E_{X^l}}$. Then, by Lemma 4.8, there exists $c\in 
(k_X^{l})^{\times}$ such that $f-c\in \Cal O_{E_{X^l}}^{\times}$. Thus, 
$f=c+(f-c)\in (k_X^{l})^{\times}+\Cal O_{E_{X^l}}^{\times}$, as desired. 

\noindent
(ii) The first assertion follows directly from (i). The second assertion follows from 
the first assertion, together with the fact that $\Cal O_{E_{X^l}}$ is a $k_X^{l}$-algebra. 
\qed
\enddemo

So far, we have only resorted to the assumptions 
that $\Sigma$ is $k_X$-large and $k_Y$-large and 
that at least one of $X$ and $Y$ is almost $\Sigma$-separated. From now, we will resort to 
the other assumption that $\alpha$ is pseudo-constants-additive. 

\proclaim
{Lemma 4.11} Let $f\in H_{X^l}^{\times}$, and assume that $1+f\in \Cal O_{E_{X^l}}^{\times}$. 
Then $1+f\in H_{X^l}^{\times}$.
\endproclaim

\demo{Proof} Replacing $k_X$, $k_Y$ by suitable 
finite subextensions of 
$k_X^{l}/k_X$ and $k_Y^{l}/k_Y$ corresponding to each other by $\alpha$, 
we may assume that $f\in H_X^{\times}$. 
Then the proof is similar to [Sa\"\i di-Tamagawa3], Lemma 3.16, 
using Proposition 2.4. 
More precisely, write $(\rho')^{-1}(f')=g'$. We have $((1+f)')^m\in {H_{X}'}$ 
and we may write $(\rho')^{-1}(((1+f)')^m)=h'$ with $h\in H_{Y}^{\times}$ 
for some positive integer $m$ divisible only by primes in $\Sigma'$. 
By evaluating the pseudo-function $h'$ at all points in 
$Y^{\cl}\setminus 
Y(F_Y)^{\cl}
$, and using Lemma 4.6 
and Proposition 2.4, 
we deduce that $h'=(\rho')^{-1}(((1+f)')^m)=((1+cg)')^m\in {H_{Y}'}\subset 
\Cal D_{Y\setminus E_{Y}}$ for some $c\in (k_Y)^{\times}\{\Sigma'\}$ 
(cf. loc. cit. for more details). Hence, 
$(\rho')^{-1}((1+f)')=(1+cg)'\in \Cal D_{Y\setminus E_{Y}}$ since $\Cal D_{Y\setminus E_{Y}}$ 
admits no nontrivial $\Sigma'$-primary torsion ($D_{Y\setminus E_{Y}}$ 
has no nontrivial 
torsion, and the kernel of the natural surjective 
map $\Cal D_{Y\setminus E_{Y}}\to D_{Y\setminus E_{Y}}$ is 
$\Sigma$-primary torsion).  As $(1+cg)'\in \Cal O_{E_Y}^{\times}/(k_{Y}^{\times}\{\Sigma'\})$, 
we have $(1+f)'\in  H_{X}'$ by definition, as desired.
\qed
\enddemo

\proclaim
{Lemma 4.12} 
{\rm (i)} $R_{X^l}
$ coincides with the 
$k_X^{l}$-subalgebra generated by $H_{X^l}^{\times}$. 

\noindent
{\rm (ii)} Let $f_1,\dots,f_n\in R_{X^l}$. 
Then, for all but finitely many $c\in k_X^{l}$, 
$f_1- c,\dots,f_n-c\in H_{X^l}^{\times}$. 

\noindent
{\rm (iii)} $(k_X^{l})^{\times}+H_{X^l}^{\times}=R_{X^l}$. 

\noindent
{\rm (iv)} $R_{X^l}\cap\Cal O_{E_{X^l}}^{\times}=R_{X^l}^{\times}=H_{X^l}^{\times}$. 
\endproclaim

\demo{Proof} 
(i) First, note that $k_X^{l}=(k_X^{l})^{\times}\cup \{0\}$ is contained in 
$R_{X^l}=\langle H_{X^l}^{\times}\rangle$, as $(k_X^{l})^{\times}\subset H_{X^l}^{\times}$. 
Thus, it suffices to prove that $R_{X^l}$ is stable under multiplication. But this just follows from 
the fact that $H_{X^l}^{\times}$ (which is a multiplicative subgroup) is stable under multiplication. 

\noindent
(ii) We may and shall assume $n=1$ and set $f\defeq f_1$. (The assertion for general $n$ can be 
reduced to this special case immediately.) As $f\in R_{X^l}=\langle H_{X^l}^{\times}\rangle$, 
$f$ can be written as $f=g_1+\dots+ g_m$ with $g_1,\dots,g_m\in H_{X^l}^{\times}$. We shall prove the 
assertion for $f$ by induction on $m$. 
The case $m=0$ (i.e., $f=0$) is easy: 
any $c\in k_X^{l}\setminus\{0\}$ satisfies the desired property. 
The case $m=1$ (i.e., $f=g_1\in H_{X^l}^{\times}$) follows immediately from Lemma 4.8 and Lemma 4.11. 
More precisely, by Lemma 4.8 (and the case $m=0$), for all but finitely many $c\in k_X^l$, one has 
$f-c\in \Cal O_{E_{X^l}}^{\times}$ and $-c\in H_{X^l}^{\times}$, hence 
$$f-c=\frac{f-c}{f}\cdot f=\left(\frac{-c}{f}+1\right) f\in H_{X^l}^{\times}$$
by Lemma 4.11. Now, assume $m>1$ and suppose that 
the assertion holds for $m-1$. Then it follows from Lemma 4.8 and the induction hypothesis that,  
for all but finitely many $c\in k_X^{l}$, one has 
$f-c\in\Cal O_{E_{X^l}}^{\times}$ and $(g_1+\dots+g_{m-1})-c\in H_{X^l}^{\times}$. Now, 
as 
$$f-c=\frac{f-c}{g_m}\cdot g_m=\left(\frac{(g_1+\dots+g_{m-1})-c}{g_m}+1\right)g_m,$$
one has $f-c\in H_{X^l}^{\times}$ by Lemma 4.11. 

\noindent
(iii) Clearly $(k_X^{l})^{\times}+H_{X^l}^{\times}\subset R_{X^l}$. 
To show the opposite, take any $f\in R_{X^l}$. Then, by (ii), there exists $c\in 
(k_X^{l})^{\times}$ such that $f-c\in H_{X^l}^{\times}$. Thus, 
$f=c+(f-c)\in (k_X^{l})^{\times}+H_{X^l}^{\times}$, as desired. 

\noindent
(iv) Clearly $R_{X^l}\cap\Cal O_{E_{X^l}}^{\times}\supset R_{X^l}^{\times}\supset H_{X^l}^{\times}$, 
hence it suffices to prove 
$R_{X^l}\cap\Cal O_{E_{X^l}}^{\times}\subset H_{X^l}^{\times}$. So, take any 
$f\in R_{X^l}\cap\Cal O_{E_{X^l}}^{\times}$. By (iii), there exist $c\in (k_X^{l})^{\times}$ 
and $g\in H_{X^l}^{\times}$ such that $f=c+g$. As 
$$f=\frac{f}{c}\cdot c=\left(1+\frac{g}{c}\right)c,$$
one has $f\in H_{X^l}^{\times}$ by Lemma 4.11. 
\qed\enddemo

\proclaim{Lemma 4.13}
$\Fr(\Cal O_{E_{X^l}})=\Fr(R_{X^l})$ and $\Cal O_{E_{X^l}}$ is the normalization of 
$R_{X^l}$. 
\endproclaim

\demo{Proof}
Write $K_{X^l}=\Fr(\Cal O_{E_{X^l}})\ (=K_Xk_X^{l})$ and 
$N_{X^l}=\Fr(R_{X^l})$. 

\medskip\noindent
{\it Step 1.} $\Cal O_{E_{X^l}}$ is the integral closure of 
$R_{X^l}$ in $K_{X^l}$. 

Indeed, as $\Cal O_{E_{X^l}}$ is the intersection of discrete valuation rings 
$\Cal O_{X^l,x}$ ($x\in E_{X^l}$), $\Cal O_{E_{X^l}}$ is integrally closed. 
On the other hand, as $\Cal O_{E_{X^l}}^{\times}/R_{X^l}^{\times}$ is torsion, 
each element of $\Cal O_{E_{X^l}}^{\times}$ is integral over $R_{X^l}$. 
As $\Cal O_{E_{X^l}}=\langle\Cal O_{E_{X^l}}^{\times}\rangle$, 
$\Cal O_{E_{X^l}}$ is integral over $R_{X^l}$, as desired. 

\medskip\noindent
{\it Step 2.} $\Cal O_{E_{X^l}}^{\times}/(k_X^{l})^{\times}$ and 
$R_{X^l}^{\times}/(k_X^{l})^{\times}$ are free $\Bbb Z$-modules of countably infinite rank. 

Indeed, each of these groups is injectively mapped into $\Div^0_{X^l\setminus E_{X^l}}$ 
with torsion cokernel. Now, the assertion follows from the 
fact that $\Div^0_{X^l\setminus E_{X^l}}$ itself is 
a free $\Bbb Z$-module of countably infinite rank. 
(Note that $(X^l)^{\cl}\setminus E_{X^l}$ is a countably infinite set.) 

\medskip\noindent
{\it Step 3.} $K_{X^l}/N_{X^l}$ is finite. 

Indeed, since $k_X^{l}\subset N_{X^l}\subset K_{X^l}$ and 
$K_{X^l}$ is a (regular) function field of one variable 
over $k_X^{l}$, it suffices to prove that $N_{X^l}\supsetneq k_X^{l}$. 
But this follows, for example, from Step 1 or Step 2. 

\medskip\noindent
{\it Step 4.} $K_{X^l}/N_{X^l}$ is separable. 

Indeed, otherwise, $N_{X^l}\subset (K_{X^l})^p$, since 
$K_{X^l}/N_{X^l}$ is a finite extension of (regular) function fields of one variable 
over a perfect field $k_X^{l}$. Then 
$R_{X^l}\subset \Cal O_{E_{X^l}}\cap (K_{X^l})^p=(\Cal O_{E_{X^l}})^p$ 
(where the last equality follows from the fact that $\Cal O_{E_{X^l}}$ 
is normal), and 
$R_{X^l}^{\times}\subset((\Cal O_{E_{X^l}})^p)^{\times}=
(\Cal O_{E_{X^l}}^{\times})^p$. Thus, one has 
$\Cal O_{E_{X^l}}^{\times}/R_{X^l}^{\times}\twoheadrightarrow 
\Cal O_{E_{X^l}}^{\times}/(\Cal O_{E_{X^l}}^{\times})^p$. 
But this is impossible, since the $p$-primary part of 
the torsion abelian group $\Cal O_{E_{X^l}}^{\times}/R_{X^l}^{\times}$ 
is embeddable into a finite direct sum of $\Bbb Q_p/\Bbb Z_p$ (hence 
is isomorphic to a direct sum of a finite number of 
copies of $\Bbb Q_p/\Bbb Z_p$ and a $p$-primary finite abelian group), 
while 
$\Cal O_{E_{X^l}}^{\times}/(\Cal O_{E_{X^l}}^{\times})^p$ is an 
$\Bbb F_p$-vector space of countably infinite dimension by Step 2. 

\medskip\noindent
{\it Step 5.} $K_{X^l}^{\times}/N_{X^l}^{\times}$ has finite torsion. 

Indeed, the homomorphism $N_{X^l}\to K_{X^l}$ of fields comes from 
a finite, generically \'etale $k_X^{l}$-morphism 
$X^l\to Z$ (where $Z$ is a proper, smooth, geometrically connected 
curve over $k_X^{l}$ with function field $N_{X^l}$) 
of degree $d\defeq [K_{X^l}:N_{X^l}]$. 
Then 
$$K_{X^l}^{\times}/N_{X^l}^{\times}
=(K_{X^l}^{\times}/(k_X^{l})^{\times})
/
(N_{X^l}^{\times}/(k_X^{l})^{\times})
=\Pri_{X^l}/\Pri_Z.$$
Considering the commutative diagram with two rows exact: 
$$
\CD
0@>>>\Pri_Z@>>>\Div_Z@>>>\Pic_Z@>>>0\\
@.@VVV@VVV@VVV@.\\
0@>>>\Pri_{X^l}@>>>\Div_{X^l}@>>>\Pic_{X^l}@>>>0
\endCD
$$
in which the vertical arrows are induced by the pull-back of divisors 
by the morphism $X^l\to Z$, 
one obtains an exact sequence 
$$0\to\Ker(\Pic_Z\to\Pic_{X^l})\to 
\Pri_{X^l}/\Pri_Z \to \Div_{X^l}/\Div_Z.$$
Now, on the one hand, by considering the norm map $\Pic_{X^l}\to\Pic_Z$, 
one sees that 
$$\Ker(\Pic_Z\to\Pic_{X^l})\subset \Pic_Z[d]=\Pic^0_Z[d]=J_Z(k_X^{l})[d],$$ 
hence that $\Ker(\Pic_Z\to\Pic_{X^l})$ is finite. On the other hand, 
by considering the definition of $\Div_Z\to\Div_{X^l}$, one sees that 
the torsion of $\Div_{X^l}/\Div_Z$ (all of which arises from the finitely many 
ramified points of the generically \'etale morphism $X^l\to Z$) is finite. 
Thus, the assertion follows. 

\medskip\noindent
{\it Step 6.} Let $\tilde R_{X^l}$ denote the normalization of $R_{X^l}$ in $N_{X^l}$. 
Then $(\Cal O_{E_{X^l}}^{\times}:\tilde R_{X^l}^{\times})<\infty$. 

Indeed, as $\Cal O_{E_{X^l}}^{\times}/\tilde R_{X^l}^{\times}\twoheadleftarrow 
\Cal O_{E_{X^l}}^{\times}/R_{X^l}^{\times}$, 
$\Cal O_{E_{X^l}}^{\times}/\tilde R_{X^l}^{\times}$ is torsion. On the other hand, 
as $\Cal O_{E_{X^l}}$ is integral over $\tilde R_{X^l}$ and 
$\tilde R_{X^l}$ is integrally closed, one has 
$\Cal O_{E_{X^l}}^{\times}/\tilde R_{X^l}^{\times}\hookrightarrow 
K_{X^l}^{\times}/N_{X^l}^{\times}$. Now, the assertion follows from Step 5. 

\medskip\noindent
{\it Step 7.} End of proof: $\Cal O_{E_{X^l}}=\tilde R_{X^l}$ and $K_{X^l}=N_{X^l}$. 

Indeed, by Step 6, we may write 
$\Cal O_{E_{X^l}}^{\times}
=\tilde R_{X^l}^{\times}f_1\cup \dots\cup \tilde R_{X^l}^{\times}f_r$ 
for some $f_1,\dots, f_r\in \Cal O_{E_{X^l}}^{\times}$. 
This, together with Corollary 4.10(i), implies 
that $\Cal O_{E_{X^l}}=(\tilde R_{X^l}f_1\cup \dots\cup \tilde R_{X^l}f_r)+k_X^{l}$, 
and that $\Cal O_{E_{X^l}}/k_X^{l}=\overline{\tilde R_{X^l} f_1}\cup 
\dots\cup\overline {\tilde R_{X^l}f_r}$ as a union of $k_X^{l}$-vector 
subspaces 
(where $\overline {\tilde R_{X^l}f_i}$ denotes 
the image of 
$\tilde R_{X^l}f_i$ in $\Cal O_{E_{X^l}}/k_X^{l}$). 
As $k_X^{l}$ is an infinite field, we must have 
$\Cal O_{E_{X^l}}/k_X^{l}=\overline {\tilde R_{X^l}f_i}$ for some $i$, 
hence $\Cal O_{E_{X^l}}=k_X^{l}+\tilde R_{X^l}f_i$. 

We claim that $f\defeq f_i\in N_{X^l}$. Indeed, otherwise, $1, f\in K_{X^l}$ are 
linearly independent over $N_{X^l}$. Namely, 
$N_{X^l}\oplus N_{X^l}f\hookrightarrow K_{X^l}$, 
hence, in particular, 
$\tilde R_{X^l}\oplus \tilde R_{X^l}f\hookrightarrow \Cal O_{E_{X^l}}$. 
As $k_X^{l}\subsetneq R_{X^l}\subset \tilde R_{X^l}$, one has
$k_X^{l}\oplus \tilde R_{X^l}f
\subsetneq \tilde R_{X^l}\oplus \tilde R_{X^l}f$. This implies that 
$k_X^{l}+\tilde R_{X^l}f\subsetneq 
R_{X^l}+ \tilde R_{X^l}f\subset \Cal O_{E_{X^l}}$, which is absurd. 

Now, $f\in N_{X^l}\cap \Cal O_{X^l}^{\times}=\tilde R_{X^l}^{\times}$, where the last equality 
follows from Step 1, together with the definition of $\tilde R_{X^l}$. 
Thus, 
$\Cal O_{E_{X^l}}=k_X^{l}+\tilde R_{X^l}f=
k_X^{l}+\tilde R_{X^l}=\tilde R_{X^l}$, as desired. 
In particular, 
$K_{X^l}=\Fr(\Cal O_{E_{X^l}})=\Fr(\tilde R_{X^l})=N_{X^l}$. 
\qed
\enddemo

Next, we will denote by $\Bbb P(R_{X^l})\defeq (R_{X^l}\setminus \{0\})/(k_X^{l})^{\times}$ 
(resp. $\Bbb P( R_{Y^l})\defeq  
( R_{Y^l}\setminus \{0\})/(k_Y^{l})^{\times}$) 
the projective space associated to the infinite-dimensional $k_X^{l}$-vector space
$R_{X^l}$ (resp. the $k_Y^{l}$-vector space $ R_{Y^l}$) 
and 
by $\Bbb L(R_{X^l})\ (\subset 2^{\Bbb P(R_{X^l})})$ 
(resp. $\Bbb L( R_{Y^l})\ (\subset 
2^{\Bbb P( R_{Y^l})})$) 
the set of lines on 
$\Bbb P(R_{X^l})$ (resp. $\Bbb P( R_{Y^l})$).  
We view $\Cal U_{X^l}\defeq \overline H_{X^l}=R_{X^l}^{\times}/(k_X^{l})^{\times}$ and $\Cal U_{Y^l}\defeq \overline H_{Y^l}=
 R_{Y^l}^{\times}/(k_Y^{l})^{\times}$ (cf. Lemma 4.12(iv))
as subsets of the projective spaces $\Bbb P(R_{X^l})$ and $\Bbb P( R_{Y^l})$, respectively.
Let $F_X/k_X$ and $F_Y/k_Y$ be the extensions in Definition 2.9. 
We define the following sets of subsets of $\Bbb P^1({k_X^{l}})$ and $\Bbb P^1({k_Y^{l}})$: 
$$\Cal S_{X^l}\defeq \{S\subset \Bbb P^1({k_X^{l}})\ \vert\  \sharp(S)<\infty\},\ 
\Cal S_{Y^l}\defeq \{S\subset \Bbb P^1({k_Y^{l}})\ \vert\  \sharp(S)<\infty\}.$$

\proclaim {Lemma 4.14} The sets $\Cal S_{X^l}$ and $\Cal S_{Y^l}$ are 
admissible 
(cf. Definition 5.4(ii) for the meaning of being 
admissible).
\endproclaim

\demo {Proof} 
%
%
This follows from the fact that $k_X^{l}$ and $k_Y^{l}$ 
are infinite. See Remark 5.5. 
\qed
\enddemo

\proclaim {Lemma 4.15} The subset 
$\Cal U_{X^l}\subset \Bbb P(R_{X^l})$ is $\Cal S_{X^l}$-ample, and the subset 
$\Cal U_{Y^l}\subset \Bbb P(R_{Y^l})$ is $\Cal S_{Y^l}$-ample (cf. Definition 5.6(iii)
for the meaning of being $\Cal S_{X^l}$-ample and $\Cal S_{Y^l}$-ample).
\endproclaim

\demo{Proof} We prove that $\Cal U_{X^l}$ is $\Cal S_{X^l}$-ample. 
The proof that $\Cal U_{Y^l}$ is $\Cal S_{Y^l}$-ample is similar.
We have to prove the following equality $\Bbb L(R_{X^l})_{\Cal U_{X^l}}=\Bbb L(R_{X^l})_{\Cal U_{X^l},\Cal S_{X^l}}$ (cf. Definition 5.6(iii)(2)). More precisely,
let $\ell \in \Bbb L(R_{X^l})$ be a line in the projective space  $\Bbb P(R_{X^l})$ such that $\ell _{\Cal U_{X^l}}\defeq
\ell \cap \Cal U_{X^l}\neq \varnothing$, then we have to prove
that $\ell \setminus \ell _{\Cal U_{X^l}}\in \Cal S_{X^l}$, i.e., 
$\sharp(\ell \setminus \ell _{\Cal U_{X^l}})<\infty$. We have $\ell =\Bbb P(V)$, where 
$V$ is a $2$-dimensional $k_X^{l}$-vector subspace of $R_{X^l}$. 
As $\ell _{\Cal U_{X^l}}\neq \varnothing$ we can take 
$f\in V\cap R_{X^l}^{\times}
$. Further, taking any $g\in V\setminus k_X^{l}f$, we may write 
$V=\{af+bg\ \vert\ a,b\in k_X^{l}\}$. 
Then $\ell =\{\overline f\}\cup \{\overline {cf+g},\ c\in k_X^{l}\}
=\{\overline f\}\cup \{\overline {(h-c)f},\ c\in k_X^{l}\}$, where $h\defeq \frac {g}{f}
\in R_{X^l}$. Now, by Lemma 4.12(ii), $h-c\in R_{X^l}^{\times}$ 
(hence $\overline{(h-c)f} \in \Cal U_{X^l}$) 
for all but finitely many $c\in k_X^{l}$, as desired. 
\qed
\enddemo

\proclaim {Lemma 4.16} The natural isomorphism $\bar \rho:\Cal U_{Y^l}\isom \Cal U_{X^l}$ 
induces a natural 
bijection $\bar \tau:\Bbb L( R_{Y^l})_{\Cal U_{Y^l}}\isom \Bbb L( R_{X^l})_{\Cal U_{X^l}}$ with the following property:
If $\ell \in \Bbb L( R_{Y^l})_{\Cal U_{Y^l}}$ then 
$\bar \tau (\ell )_{\Cal U_{X^l}}
=\bar \rho (\ell _{\Cal U_{Y^l}})$, where 
$\bar \tau (\ell )_{\Cal U_{X^l}}\defeq \bar \tau (\ell )\cap {\Cal U_{X^l}}$ and 
$\ell _{\Cal U_{Y^l}}\defeq \ell \cap {\Cal U_{Y^l}}$.
\endproclaim

\demo{Proof} We will define the map $\bar \tau$.
Let $\ell \in \Bbb L( R_{Y^l})_{\Cal U_{Y^l}}$. 
Then by Lemma 4.15, $\ell \setminus \ell _{\Cal U_{Y^l}}$ is finite, 
hence $\ell _{\Cal U_{Y^l}}$ is infinite. Take two distinct 
points $\overline f_1, \overline g_1\in \ell _{\Cal U_{Y^l}}$, 
and take any liftings $f_1,g_1\in R_{Y^l}^{\times}$ of 
$\overline f_1, \overline g_1$, respectively. 
Set $V_1=\{af_1+bg_1\ \vert\ a,b\in k_Y^{l}\}$, so that 
$V_1$ is a $2$-dimensional $k_Y^{l}$-vector subspace of $ R_{Y^l}$ 
and that $\ell =\Bbb P(V_1)\subset\Bbb P(R_{Y^l})$. Then we have 
$\ell =\{\overline g_1\}\cup\{\overline{(1+ch_1)f_1}\mid c\in k_Y^{l}\}$, 
where $h_1\defeq\frac{g_1}{f_1}\in R_{Y^l}^{\times}$, and 
$\ell _{\Cal U_{Y^l}}=\{\overline g_1\}\cup\{\overline{(1+ch_1)f_1}\mid c\in k_Y^{l}, 
(1+ch_1)\in R_{Y^l}^\times\}$. 
Take any $c\in (k_Y^{l})^{\times}$ such that 
$(1+ch_1)\in R_{Y^l}^\times$. Set 
$f_2'\defeq\rho'(f_1'), g_2'\defeq\rho'(g_1')$, where 
$\rho':H_{Y^l}'= R_{Y^l}^{\times}/((k_Y^{l})^{\times}\{\Sigma '\}) 
\isom 
H_{X^l}'= R_{X^l}^{\times}/((k_X^{l})^{\times}\{\Sigma '\})$ 
is the natural isomorphism (cf. discussion before Lemma 4.11), 
and take any liftings $f_2,g_2\in R_{X^l}^{\times}$ of 
$f_2', g_2'$, respectively. Set $h_2\defeq\frac{g_2}{f_2}\in R_{X^l}^{\times}$, 
so that $h_2'=\frac{g_2'}{f_2'}=\rho'(h_1')$. 
Further set 
$V_2=\{af_2+bg_2\ \vert\ a,b\in k_X^{l}\}$. 
Since $\bar\rho$ is bijective, one has $\overline f_2\neq\overline g_2$, hence 
$V_2$ is a $2$-dimensional $k_X^{l}$-vector subspace of 
$R_{X^l}$. Set $\bar\tau(\ell )\defeq \ell '\defeq \Bbb P(V_2)\subset \Bbb P(R_{X^l})$. As 
$\overline f_2\in \ell '$, one has $\ell '\in\Bbb L(R_{X^l})_{\Cal U_{X^l}}$. 

By evaluating the pseudo-function $\rho'((1+ch_1)')$ at all points 
$x\in (X^l)^{\cl}$ above $X^{\cl}\setminus  
X(F_X)^{\cl}
$, 
and applying Lemma 4.5, Lemma 4.6, and  Proposition 2.3 to suitable finite 
subextensions of $k_X^{l}/k_X$ and $k_Y^{l}/k_Y$, 
we deduce that $\rho'((1+ch_1)')=(\alpha +\beta h_2)'   
\in  R_{X^l}^{\times}/(k_{X^{l}}^{\times})\{\Sigma '\}$ for some 
$\alpha\in (k_X^{l})^{\times}\{\Sigma'\}$ and 
$\beta\in (k_X^{l})^{\times}$, hence $\bar \rho (\overline{1+ch_1})
=\overline {1+\gamma h_2}$ for $\gamma=\frac{\beta}{\alpha}\in (k_X^{l})^{\times}$. 
This, together with $\bar\rho(\overline f_1)
=\overline f_2, \bar\rho(\overline g_1)=\overline g_2$, 
implies that $\bar\rho(\ell _{\Cal U_{Y^l}})\subset \bar\tau(\ell )_{\Cal U_{X^l}}$. 
Further, $\bar\tau$ is bijective since 
it has an inverse $\bar\tau':\Bbb L( R_{X^l})_{\Cal U_{X^l}}\isom \Bbb L( R_{Y^l})_{\Cal U_{Y^l}}$ 
which is naturally deduced from $\alpha ^{-1}:\Pi_Y\isom \Pi_X$. By using $\bar\tau'$, 
we also conclude that $\bar\rho(\ell _{\Cal U_{Y^l}})=\ell '_{\Cal U_{X^l}}$, as desired. 
\qed
\enddemo

\proclaim
{Lemma 4.17} (Recovering the Additive Structure of $ R_{X^l}$ and $ R_{Y^l}$) The following hold.

\noindent
{\rm (i)} There exist natural isomorphisms
$\tilde \rho: \Bbb P( R_{Y^l})\isom \Bbb P( R_{X^l})$
and $\tilde \tau : \Bbb L( R_{Y^l})\isom \Bbb L( R_{X^l})$ 
which extend the isomorphisms $\bar \rho:\Cal U_{Y^l}\isom \Cal U_{X^l}$ 
and $\bar \tau:\Bbb L( R_{Y^l})_{\Cal U_{Y^l}}\isom \Bbb L( R_{X^l})_{\Cal U_{X^l}}$, respectively, 
such that for every $\ell\in 
\Bbb L( R_{Y^l})$, we have $\tilde \tau (\ell)=\tilde \rho (\ell)$ set-theoretically.

\noindent
{\rm (ii)} The bijection $\tilde \rho$ arises from a $\psi _0$-isomorphism 
$$\psi: (  R_{Y^l} ,+) \isom  ( R_{X^l},+),$$
where  $\psi _0:k_{Y}^{l} \isom k_{X}^{l}$ is a field isomorphism. Namely, 
$\psi$ is an isomorphism of abelian groups which is compatible with $\psi _0$ 
in the sense that $\psi(ax)=\psi_0(a)\psi(x)$ for $a\in k_{Y}^{l}$ and $x\in  R_{Y^l}$. 
Further, $\psi_0$ is uniquely determined and $\psi$ is uniquely determined 
up to scalar multiplication. 
\endproclaim

\demo {Proof} Assertion (i) follows formally from Lemma 4.16 and Theorem 5.7. Assertion (ii) follows from Theorem 5.7.\qed
\enddemo

\proclaim
{Lemma 4.18} (Recovering the Ring Structure of  $ R_{X^l}$ and $ R_{Y^l}$) 
The following hold.

\noindent
{\rm (i)} If we normalize the isomorphism 
$$\psi: (  R_{Y^l} ,+) \isom  ( R_{X^l},+),$$
in Lemma 4.17 by the condition $\psi(1)=1$ 
(which is possible as $\bar\rho(\bar 1)=\bar1$), it becomes a ring isomorphism such that 
the diagram
$$
\CD
R_{X^l}  @<\psi<<R_{Y^l}   \\
@AAA                    @AAA   \\
k_{X}^{l} @<\psi_0<< k_{Y}^{l}\\ 
\endCD
$$
commutes. 

\noindent
{\rm (ii)} $\psi$ induces a natural commutative diagram
$$
\CD
K_{X^l}
@<\psi<<   K_{Y^l}
\\
@AAA     @AAA  \\
K_X  @<\psi<<  K_Y \\
\endCD
$$
where the horizontal maps are field isomorphisms and the vertical maps are natural inclusions.
Further, $\psi: K_{Y^l}\isom K_{X^l}$ is Galois-equivariant with respect to the isomorphism 
$G_{k_X}\isom G_{k_Y}$ induced by $\alpha: \Pi_X\isom \Pi_Y$ (cf. Proposition 1.7(i)). 

\noindent
{\rm (iii)} $\psi$ induces a natural commutative diagram
$$
\CD
X^l
@>\psi>>   Y^l
\\
@VVV     @VVV  \\
X  @>\psi>>  Y \\
\endCD
$$
where the horizontal maps are scheme isomorphisms and the vertical maps are natural projections. 
Further, $\psi: X^l\isom Y^l$ is Galois-equivariant with respect to the isomorphism 
$G_{k_X}\isom G_{k_Y}$ induced by $\alpha: \Pi_X\isom \Pi_Y$ (cf. Proposition 1.7(i)). 
\endproclaim

\demo{Proof} (i) The proof is similar to the proof of [Sa\"\i di-Tamagawa3], Lemma 4.14. More precisely,
let $f\in  R_{Y^l}^{\times}$, then $\psi \circ \mu_f$ and $\mu _{\psi (f)}\circ \psi$ are $\psi_0$-isomorphisms
$( R_{Y^l} ,+)\isom ( R_{X^l} ,+)$, where $\mu_g$ denotes the $g$-multiplication map. Note that 
$\psi (f)\in  R_{X^l}^{\times}$, since $\overline {\psi (f)}
=\bar \rho (\overline f)\in   R_{X^l}^{\times}/(k_X^{l})^{\times}\subset
K_{X^l}^{\times}/(k_X^{l})^{\times}$. The isomorphisms $ R_{Y^l}^{\times}/(k_Y^{l})^{\times}
\isom  R_{X^l}^{\times} /(k_X^{l})^{\times}$ they induce coincide with each other:
$$\overline {\psi \circ \mu_f}=\bar \rho \circ \mu _{\bar f}
=\mu _{\bar \rho (\overline f)}\circ \bar \rho=\overline {\mu _{\psi (f)}\circ \psi},$$
where the second equality follows from the multiplicativity of $\bar \rho$. Further, we have 
$$\psi \circ \mu_f(1)=\psi (f)=\mu _{\psi (f)}(1)=\mu _{\psi (f)}\circ \psi (1).$$
Thus, the equality $\psi \circ \mu_f=\mu _{\psi (f)}\circ \psi$ follows from the uniqueness in Theorem 5.7. 
This equality means that $\psi(fg)=\psi(f)\psi(g)$ holds for $f\in R_{Y^l}^{\times}$ and $g\in R_{Y^l}$, 
which, together with the fact that $R_{Y^l}=\langle R_{Y^l}^\times \rangle$, 
implies the multiplicativity of $\psi$.

\noindent
(ii) By considering fields of fractions, we see that 
$\psi: R_{Y^l} \isom R_{X^l}$ naturally induces $\psi: K_{Y^l}\isom K_{X^l}$. 
Since the isomorphism $\bar \rho:\Cal U_{Y^l}\isom \Cal U_{X^l}$ is 
Galois-equivariant under the compatible actions of 
$G_{k_X}$ and $G_{k_Y}$ 
respectively, 
and since the isomorphism 
$\psi:  R_{Y^l} \isom   R_{X^l}$, $1\mapsto 1$ 
is uniquely determined by $\bar\rho$, it follows that 
$\psi:  R_{Y^l} \isom   R_{X^l}$ is 
Galois-equivariant under the compatible actions of 
$G_{k_X}$ and $G_{k_Y}$, 
respectively, 
hence so is $\psi:K_{Y^l}\isom K_{X^l}$. Now, 
the fact that the field isomorphism $\psi: K_{Y^l}\isom K_{X^l}$ maps $K_Y$ isomorphically to $K_X$ 
follows from this. 

\noindent
(iii) This follows immediately from (ii). 
\qed
\enddemo

Next, let $X'\to X$ be a finite \'etale covering corresponding to an open subgroup $\Pi_{X'}\subset \Pi_X$
and $Y'\to Y$ the corresponding \'etale covering via $\alpha$ (i.e., $\Pi_{Y'}\defeq \alpha (\Pi_{X'})$). 
(We refer to the case where $X'\to X$ is Galois (or, equivalently, $Y'\to Y$ is Galois) as 
a Galois case.) 
Then $\alpha$ induces, by restriction to $\Pi_{X'}$, an isomorphism 
$$\alpha:\Pi_{X'}\isom \Pi_{Y'}.$$
To apply the preceding arguments to this isomorphism, we need to show that 
the assumptions of Theorem 4.4 hold. 

\proclaim{Lemma 4.19}
{\rm (i)} $\Sigma$ is $k_{X'}$-large and $k_{Y'}$-large. 

\noindent
{\rm (ii)} $X'$ and $Y'$ are $\Sigma$-separated. More precisely, 
$F_{X'}$ and $F_{Y'}$ are finite extensions of $F_X$ and $F_Y$, 
respectively. In particular, a prime number is 
$(X,\Sigma)$-admissible 
(resp. $(Y,\Sigma)$-admissible) if and only if it is 
$(X',\Sigma)$-admissible 
(resp. $(Y',\Sigma)$-admissible). 

\noindent
{\rm (iii)} $\alpha: \Pi_{X'}\isom \Pi_{Y'}$ is pseudo-constants-additive. 
\endproclaim

\demo{Proof}
(i) This follows from the fact that $k_{X'}$ and $k_{Y'}$ are finite extensions of 
$k_X$ and $k_Y$, respectively. 

\noindent
(ii) This follows from Proposition 2.11(i). 

\noindent
(iii) By Theorem 4.2(iii), 
the isomorphism $\tau: (\bar k_{Y'}^{\times})^{\Sigma}\isom(\bar k_{X'}^{\times})^{\Sigma}$ 
induced by $\alpha: \Pi_{X'}\isom \Pi_{Y'}$ 
just coincides with 
the isomorphism $\tau: (\bar k_{Y}^{\times})^{\Sigma}\isom(\bar k_{X}^{\times})^{\Sigma}$ 
induced by $\alpha: \Pi_{X}\isom \Pi_{Y}$. Thus, the pseudo-constants-additivity 
of $\alpha: \Pi_{X'}\isom \Pi_{Y'}$ follows from that of 
$\alpha: \Pi_{X}\isom \Pi_{Y}$. 
\qed
\enddemo

Now, we may apply Lemma 4.18 to $\alpha: \Pi_{X'}\isom \Pi_{Y'}$ to obtain: 
(i) a ring isomorphism 
$\psi': R_{(Y')^l}\isom R_{(X')^l}$ 
compatible with a field isomorphism 
$\psi'_0:k_{Y'}^{l}\isom k_{X'}^{l}$; 
%
(ii) a field isomorphism 
$\psi': K_{(Y')^l}\isom K_{(X')^l}$ 
compatible with a field isomorphism 
$\psi': K_{Y'}\isom K_{X'}$ 
and Galois-equivariant with respect to the isomorphism 
$G_{k_{X'}}\isom G_{k_{Y'}}$ 
induced by 
$\alpha: \Pi_{X'}\isom\Pi_{Y'}$; and 
%
(iii) a scheme isomorphism 
$\psi': (X')^l\isom (Y')^l$ 
compatible with a scheme isomorphism 
$\psi': X'\isom Y'$ 
and Galois-equivariant with respect to the isomorphism 
$G_{k_{X'}}\isom G_{k_{Y'}}$ 
induced by 
$\alpha: \Pi_{X'}\isom\Pi_{Y'}$. Further, in the Galois case, these isomorphisms 
are Galois-equivariant with respect to the isomorphism $\alpha:\Pi_X\isom\Pi_Y$. 
Indeed, this Galois-equivariance can be proved just similarly as in 
the proof of Lemma 4.18(ii).

\proclaim{Lemma 4.20} 
{\rm (i)} The following diagram of rings is commutative: 
$$
\CD
R_{(X')^l}  @<\psi'<<R_{(Y')^l} \\
@AAA                    @AAA   \\
R_{{X}^l}  @<\psi<<R_{{Y}^l} 
\endCD
$$
In particular, the following diagram of fields is commutative: 
$$
\CD
k_{X'}^{l}  @<\psi'_0<<k_{Y'}^{l} \\
@AAA                    @AAA   \\
k_X^{l}  @<\psi_0<<k_Y^{l} 
\endCD
$$

\noindent
{\rm (ii)}
The following diagram of fields is commutative: 
$$
\CD
K_{(X')^l}  @<\psi'<<K_{(Y')^l} \\
@AAA                    @AAA   \\
K_{{X}^l}  @<\psi<<K_{{Y}^l} 
\endCD
$$
In particular, the following diagram of fields is commutative: 
$$
\CD
K_{X'}  @<\psi'<<K_{Y'} \\
@AAA                    @AAA   \\
K_{X}  @<\psi<<K_{Y} 
\endCD
$$

\noindent
{\rm (iii)} 
The following diagram of schemes is commutative: 
$$
\CD
(X')^l  @>\psi'>>(Y')^l \\
@VVV                    @VVV   \\
{X}^l  @>\psi>>{Y}^l 
\endCD
$$
In particular, the following diagram of schemes is commutative: 
$$
\CD
X'  @>\psi'>>Y' \\
@VVV                    @VVV   \\
X  @>\psi>>Y 
\endCD
$$
\endproclaim

\demo{Proof}
(i) To prove the commutativity of the first diagram, we may and shall 
assume that we are in the Galois case. 
Let 
$\iota_X$ and $\iota_Y$ denote the natural inclusions 
$R_{X^l}\to R_{(X')^l}$ and $R_{Y^l}\to R_{(Y')^l}$, 
respectively. Applying Theorem 4.2(iii) to various 
finite covers of $X$ and $Y$, we see that the diagram in question 
induces a commutative diagram 
$$
\CD
R_{(X')^l}^{\times}/(k_{X'}^{l})^{\times}  @<<<R_{(Y')^l}^{\times}/(k_{Y'}^{l})^{\times} \\
@AAA                    @AAA   \\
R_{X^l}^{\times}/(k_{X}^{l})^{\times}  @<<<R_{Y^l}^{\times}/(k_{Y}^{l})^{\times} 
\endCD
$$
In particular, we have 
$$\psi'\circ\iota_Y(R_{{Y}^l}^{\times})\cdot(k_{X'}^{l})^{\times}
=\iota_X\circ\psi(R_{Y^l}^{\times})\cdot(k_{X'}^{l})^{\times}.$$
Consider the $\Pi_{X^l}$-fixed parts of both sides of this equality. 
Then, since $\psi': 
R_{(Y')^l}  \isom R_{(X')^l}$ is Galois-equivariant with respect to 
$\alpha: \Pi_X\isom\Pi_Y$  and since $\alpha(\Pi_{X^l})=\Pi_{Y^l}$, 
we obtain 
$$\psi'\circ\iota_Y(R_{{Y}^l}^{\times})
=\iota_X\circ\psi(R_{Y^l}^{\times}),$$
hence 
$$
R\defeq \psi'\circ\iota_Y(R_{{Y}^l})
=
\langle\psi'\circ\iota_Y(R_{{Y}^l}^{\times})\rangle
=
\langle\iota_X\circ\psi(R_{Y^l}^{\times})\rangle
=
\iota_X\circ\psi(R_{Y^l})
.$$
Accordingly, each of $\psi'\circ\iota_Y$ and $\iota_X\circ\psi$ induces 
an isomorphism $R_{Y^l}\isom R$ that maps $1$ to $1$. Now, the desired equality 
$\psi'\circ\iota_X=\iota_Y\circ\psi$ follows from the uniqueness assertion 
in Theorem 5.7. 

The commutativity of the second diagram follows from that of the first diagram. 

\noindent
(ii) The commutativity of the first diagram follows from that of the first diagram in (i). 
The commutativity of the second diagram follows from that of the first diagram. 

\noindent 
(iii) This follows immediate from (ii). 
\qed
\enddemo 
 
\proclaim{Corollary 4.21} The isomorphism $\alpha$ induces a natural isomorphism
$\Tilde X \isom \Tilde Y$
which fits into a commutative diagram
$$
\CD
\Tilde X @>{\sim}>> \Tilde Y \\
@VVV   @VVV \\
X @>{\sim}>> Y \\
\endCD
$$
where the vertical maps are the pro-\'etale coverings corresponding to $\Pi_X$ and $\Pi_Y$, respectively.
\endproclaim

\demo{Proof} This follows from Lemma 4.20 by passing to open subgroups of $\Pi_X$ and $\Pi_Y$ which correspond 
to each other via $\alpha$, together with the Galois-equivariance of various isomorphisms in the Galois case. 
\qed
\enddemo

This finishes the proof of Theorem 4.4. \qed

\proclaim {Theorem 4.22} 
(A Refined Version of  the Grothendieck Conjecture for Proper Hyperbolic Curves over Finite Fields)
Let $X$, $Y$ be proper hyperbolic curves  over finite fields $k_X$, $k_Y$ of characteristic $p_X$, $p_Y$, respectively. Let $\Sigma_X, \Sigma_Y \subset \Primes$
be sets of prime numbers containing at least one prime number different from 
$p_X$, $p_Y$, respectively, 
and set $\Sigma _X'\overset \text {def}\to=\Primes \setminus  \Sigma_X$,
$\Sigma _Y'\overset \text {def}\to=\Primes \setminus  \Sigma_Y$. 
Assume that 
$\Sigma_X$ is $J_X$-large and that $\Sigma_Y$ is $J_Y$-large. 
Write $\Pi_X$, $\Pi_Y$
for the geometrically pro-$\Sigma_X$ \'etale fundamental group of $X$ and the geometrically pro-$\Sigma_Y$ \'etale fundamental group of $Y$, respectively. Let
$$\alpha:\Pi_X\isom \Pi_Y$$
be an isomorphism of profinite groups. Then $\alpha$ arises from a uniquely determined commutative diagram of schemes:
$$
\CD
\Tilde X @>{\sim}>> \Tilde Y \\
@VVV   @VVV \\
X @>{\sim}>> Y \\
\endCD
$$
in which the horizontal arrows are isomorphisms and the vertical arrows are the 
profinite \'etale coverings corresponding to the groups $\Pi_X$, $\Pi_Y$, 
respectively. 
\endproclaim

\demo{Proof} 
By Proposition 1.7, we have 
$\Sigma\defeq\Sigma_X=\Sigma_Y$. By Lemma 2.7, 
$\Sigma$ is $k_X$-large and $k_Y$-large. 
By Proposition 2.5(i), we have 
$p\defeq p_X=p_Y$. By Proposition 2.12, 
$X$ and $Y$ are almost $\Sigma$-separated. 
Now, by Theorem 4.4, it suffices to prove that $\alpha$ is pseudo-constants-additive. 

Let the notations be as in the proof of Theorem 4.4. 
To prove that $\alpha$ is pseudo-constants-additive, 
set 
$$\tilde T_{X^l}\defeq\{f\in \Cal O_{E_{X^l}} \mid 
\text{At least one pole of $f$ is $k_X^{l}$-rational.}\},$$ 
$$\tilde T_{X^l}^{\times}
\defeq\{f\in \Cal O_{E_{X^l}}^{\times} \mid 
f, f^{-1}\in \tilde T_{X^l}\},$$ 
and define 
$T_{X^l}^{\times}$ 
to be the intersection of $\tilde T_{X^l}$ with $H_{X^l}^{\times}$. 
Further, 
set 
$\overline{T_{X^l}^{\times}}\defeq T_{X^l}^{\times}/(k_X^{l})^{\times}$. 
We define 
$\tilde T_{Y^l}$, 
$\tilde T_{Y^l}^{\times}$, 
$T_{Y^l}^{\times}$ 
and 
$\overline{T_{Y^l}^{\times}}$ 
similarly. 
Since the divisors of functions are preserved under $\bar\rho$ (cf. Lemma 4.5(i)) 
and since the degrees of extensions of residue fields are preserved under 
$\phi: X^{\cl}\setminus E_X \isom Y^{\cl}\setminus E_Y$ 
(cf. discussion before Corollary 3.3), 
we have 
$\bar\rho(\overline {T_{Y^l}^{\times}})=\overline {T_{X^l}^{\times}}$. 

Observe that $\PGL_2(k_X^l)$ acts on $K_{X^l}\setminus k_X^l=\Bbb P^1(K_{X^l})\setminus\Bbb P^1(k_X^l)$ via 
linear fractional transformation: 
$$A=
\pmatrix
a&b\\
c&d
\endpmatrix
\text{mod } (k_X^{l})^{\times}\in \PGL_2(k_X^l),\  
h\in K_{X^l}\setminus k_X^l 
\implies 
A\cdot h\defeq \frac{ah+b}{ch+d}.$$
Claim: For any $h\in K_{X^l}\setminus k_X^l $, there exists $A\in \PGL_2(k_X^l)$ such that 
$f\defeq A\cdot h$ satisfies $f, f+1\in H_{X^l}^{\times}$. 
(In particular, $T_{X^l}^{\times}\neq\varnothing$.) 
A similar statement holds for $Y$. 

\medskip
Indeed, for simplicity, set $P_X\defeq\PGL_2(k_X^l)$ and write $P_Xh$ for 
the $P_X$-orbit of $h$. 
Take a finite subextension $k_0/k_X$ of $k_X^l/k_X$ such that 
$h\in K_Xk_0$. Then $h$ can be regarded as a finite $k_0$-morphism 
$X_{k_0}\to\Bbb P^1_{k_0}$. On the one hand, by the Weil estimate, we have 
$\sharp (
X_{k_0}(k_0^l)^{\cl}
)=\infty$. On the other hand, 
by Proposition 2.13, 
$\sharp(h(
X_{k_0}(k_0^l)^{\cl}
)\cap h(E_X\times_{k_X} k_0))<\infty$. 
Thus, there exist infinitely many points $x \in X(k_X^l)$ 
whose image under $h$ is not contained in $h(E_X\times_{k_X} k_0)\cup\{\infty\}$. 
Take such an $x$ and set $d\defeq h(x)\in k_X^l=\Bbb P^1(k_X^l)\setminus\{\infty\}$. 
Then the set of zeros of $h-d\in K_X^l$ intersects trivially 
with $E_{X^l}$ and includes at least one $k_X^l$-rational point. 
Thus, $h_1\defeq \frac{1}{h-d}\in\tilde T_{X^l}\cap P_Xh$. 
Next, it follows from the above argument again 
(applied to $h_1$ instead of $h$) that 
for all but finitely many $c\in k_X^l$, 
the set of zeros of $h_1-c\in K_X^l$ intersects trivially 
with $E_{X^l}$ and includes at least one $k_X^l$-rational point, 
hence 
$h_1-c\in \tilde T_{X^l}^{\times}\cap P_Xh$. 
Further, by Remark 4.7, $(\Cal O_{E_{X^l}}^{\times}: H_{X^l}^{\times})<\infty$, 
hence there exists an infinite subset $C$ of $k_X^l$ such that 
for all $c\in C$, $h_1-c\in \tilde T_{X^l}^{\times}\cap P_Xh$ and 
$h_1-c$ belongs to the same coset of 
$\Cal O_{E_{X^l}}^{\times}/H_{X^l}^{\times}$. 
Now, take mutually distinct three elements $a,b,c\in C$. 
Then 
$$f\defeq\frac{b-c}{a-b}\cdot\frac{h_1-a}{h_1-c}
\in \tilde T_{X^l}^{\times}\cap H_{X^l}^\times\cap P_Xh\subset T_{X^l}^{\times}\cap P_Xh$$ 
and 
$$f+1\defeq\frac{a-c}{a-b}\cdot\frac{h_1-b}{h_1-c}
\in \tilde T_{X^l}^{\times}\cap H_{X^l}^\times\cap P_Xh\subset T_{X^l}^{\times}\cap P_Xh,$$ 
as desired. 

As the degree of a function depends only on its $P_X$-orbit, the above claim 
particularly implies that 
the set of 
degrees of functions in $K_{X^l}\setminus k_X^l$ coincides with 
that of $H_{X^l}^\times$: 
$\deg(K_{X^l}\setminus k_X^l)=\deg(H_{X^l}^\times)$, 
and, similarly, 
that 
the set of 
degrees of functions in $K_{Y^l}\setminus k_Y^l$ coincides with 
that of $H_{Y^l}^\times$: 
$\deg(K_{Y^l}\setminus k_Y^l)=\deg(H_{Y^l}^\times)$. In particular, 
it follows from Lemma 4.5(i) that 
the gonality $\gamma_{X^l}$ of $X^l$ coincides with the gonality 
$\gamma_{Y^l}$ of $Y^l$. 

Now, take any $h\in K_{Y^l}\setminus k_Y^l$ attaining the gonality: 
$\deg(h)=\gamma_{Y^l}$. Then, by the above claim, there exists 
$f\in P_Yh$ (where $P_Y\defeq \PGL_2(k_Y^l)$) such that 
$f, f+1\in H_{Y^l}^{\times}$. 
As $f\in P_Yh$, we have $\deg(f)=\deg(h)=\gamma_{Y^l}$. 

Set $g'\defeq \rho'(f')\in H_{X^l}'$, 
$g_1'\defeq\rho'((f+1)')\in H_{X^l}'$ and 
take any lifts $g,g_1\in H_{X^l}^\times$ of 
$g',g_1'$, respectively. 
Then, by Lemma 4.5(i), $\deg(g)=\deg(f)=\gamma_{Y^l}=\gamma_{X^l}$. 

As the pole divisors of $f$ and $f+1$ coincide, 
the pole divisors of $g$ and $g_1$ coincide by Lemma 4.5(i). 
Also, as $f$ admits at least one $k_Y^l$-rational pole, 
$g$ admits at least one $k_X^l$-rational pole, say, $x$. 
Now, 
by considering 
the leading terms of Laurent expansions of $g$ and $g_1$ at the 
$k_X^l$-rational pole $x$, we see that 
there exists $\beta\in (k_X^{l})^{\times}$,  
such that the pole divisor $D_\beta$ of $g_1-\beta g$ is 
strictly smaller than the pole divisor $D$ of $g$. 
(That is to say, the divisor $D-D_{\beta}$ is effective and non-zero.) 
As $\deg(g)=\gamma_{X^l}$, this implies that $g_1-\beta g$ is constant, hence 
we may write $g_1=\alpha+\beta g$ with $\alpha\in k_X^l$. 
By evaluating this equation at $\phi^{-1}(y_1)$,  where $y_1\in (Y^l)^{\cl}$ is a zero of 
$f$, we see $\alpha'=1'$, i.e., $\alpha\in (k_X^{l})^{\times}\{\Sigma'\}$. Similarly, 
by evaluating the equation $\frac{g_1}{g}=\frac{\alpha}{g}+\beta$ at $\phi^{-1}(y_2)$, where 
$y_2\in (Y^l)^{\cl}$ is a pole of $g$, we see $\beta'=1'$, i.e., 
$\beta\in (k_X^{l})^{\times}\{\Sigma'\}$. 
Now, the proof of the assertion follows just 
similarly to the last paragraph of 
the proof of [Sa\"\i di-Tamagawa3], Lemma 4.9. 
\qed
\enddemo

As a consequence of Theorem 4.22 we can deduce the following refined version of  the Grothendieck conjecture 
for (not necessarily proper) hyperbolic curves over finite fields.

\proclaim {Theorem 4.23} (A Refined Version of  the Grothendieck Conjecture for 
(Not Necessarily Proper) Hyperbolic Curves over Finite Fields)
Let $U$, $V$ be (not necessarily proper) hyperbolic curves
over finite fields $k_U$, $k_V$ of characteristic $p_U$, $p_V$, respectively. Let 
$\Sigma_U, \Sigma_V\subset \Primes$ 
be sets of prime numbers 
and 
set $\Sigma _U'\overset \text {def}\to=\Primes \setminus \Sigma_U$,
$\Sigma _V'\overset \text {def}\to=\Primes \setminus \Sigma_V$. Write $\Pi_U$, $\Pi_V$,
for the geometrically pro-$\Sigma_U$ tame fundamental group of $U$ and 
the geometrically pro-$\Sigma_V$ tame fundamental group of $V$, respectively.  
Let
$$\alpha:\Pi_U\isom \Pi_V$$
be an isomorphism of profinite groups. Assume that 
there exist open subgroups $\Pi_{U'}\subset \Pi_U$,  $\Pi_{V'}\subset \Pi_V$, 
which correspond to each other via 
$\alpha$, i.e., $\Pi_{V'}\defeq \alpha (\Pi_{U'})$,
corresponding to \'etale coverings $U'\to U$, $V'\to V$, 
such 
that the smooth compactifications $X'$ of $U'$ and $Y'$ of $V'$ are hyperbolic, 
that $\Sigma_U$ is $J_{X'}$-large 
and 
that $\Sigma_V$ is $J_{Y'}$-large. 
Then $\alpha$ arises from a uniquely determined commutative diagram of schemes:
$$
\CD
\Tilde U @>{\sim}>> \Tilde V  \\
@VVV   @VVV \\
U @>{\sim}>> V \\
\endCD
$$
in which the horizontal arrows are isomorphisms and the vertical arrows are the 
profinite \'etale coverings corresponding to the groups $\Pi_U$, $\Pi_V$, 
respectively. 
\endproclaim

The rest of this section is devoted to the proof of Theorem 4.23. 
Let $X$, $Y$, $X'$, $Y'$ be the smooth compactifications of $U$, $V$, $U'$, $V'$, respectively. 
We consider any open subgroup $\Pi_{U''}\subset \Pi_U$ corresponding to an 
\'etale covering $U''\to U$ 
such 
that $\Pi_{U''}$ is normal in $\Pi_U$, 
and that $\Pi_{U''}$ is contained in $\Pi_{U'}$. We refer to such a subgroup 
$\Pi_{U''}\subset\Pi_U$ (resp. a covering $U''\to U$)  
as a nice subgroup of $\Pi_U$ (resp. a nice covering of $U$). 
Set $\Pi_{V''}\defeq\alpha(\Pi_{U''})$, which corresponds to an \'etale covering 
$V''\to V$. 
Let $X''$, $Y''$ be the smooth compactification of $U''$, $V''$, respectively. 
(Note that $X''$, $Y''$ are hyperbolic, as they dominate $X'$, $Y'$, respectively.) 
By Theorem 1.8(i), $\alpha: \Pi_{U''}\isom \Pi_{V''}$ induces 
$\alpha: \Pi_{X''}\isom \Pi_{Y''}$. 

Note that by (the proof of) Theorem 4.22, 
applied to $\alpha: \Pi_{X'}\isom \Pi_{Y'}$ 
induced by $\alpha: \Pi_{U'}\isom\Pi_{V'}$, we have 
$\Sigma\defeq\Sigma_U=\Sigma_V$; $p\defeq p_U=p_V$; 
$\Sigma$ is $k_{X'}$-large and $k_{Y'}$-large; 
$X'$ and $Y'$ are almost $\Sigma$-separated; and 
$\alpha: \Pi_{X'}\isom\Pi_{Y'}$ is pseudo-constants-additive. 

\proclaim{Lemma 4.24} 
{\rm (i)} $\Sigma$ is $k_{X''}$-large and $k_{Y''}$-large. 

\noindent
{\rm (ii)} $X''$ and $Y''$ are almost $\Sigma$-separated. 

\noindent
{\rm (iii)} $\alpha: \Pi_{X''}\isom\Pi_{Y''}$ is pseudo-constants-additive. 
\endproclaim

\demo{Proof}
(i) This follows from the fact that $k_{X''}$ (resp. $k_{Y''}$) is 
a finite extension of $k_{X'}$ (resp. $k_{Y'}$). 

\noindent
(ii) This follows from Proposition 2.11, together with the fact that 
$X''\to X'$, $Y''\to Y'$ are tame-Galois. 

\noindent
(iii) This follows from the fact that the diagram 
$$
\CD
(\bar k_{X''}^{\times})^{\Sigma}@<\tau<< (\bar k_{Y''}^{\times})^{\Sigma} \\
@AAA   @AAA \\
(\bar k_{X'}^{\times})^{\Sigma} @<\tau << (\bar k_{Y'}^{\times})^{\Sigma}
\endCD
$$
is commutative, where the horizontal arrows are isomorphisms induced 
by $\alpha: \Pi_{X''}\isom \Pi_{Y''}$ and $\alpha: \Pi_{X'}\isom \Pi_{Y'}$, 
and the vertical arrows are 
isomorphisms induced (via Kummer theory) by the natural homomorphisms  
$\Pi_{X''}\to\Pi_{X'}$ and 
$\Pi_{Y''}\to\Pi_{Y'}$ 
(induced by 
$\Pi_{U''}\hookrightarrow\Pi_{U'}$ and 
$\Pi_{V''}\hookrightarrow\Pi_{V'}$, respectively). 
Here, as in the proof of Theorem 4.2, 
the commutativity of this diagram follows basically from 
the functoriality of Kummer theory, together with the commutativity of the diagram 
$$
\CD
M_{X''} @<<<  M_{Y''} \\
@AAA     @AAA  \\
M_{X'}  @<<< M_{Y'}   \\
\endCD
$$
where the horizontal arrows are isomorphisms induced by $\alpha$ 
(or, more precisely, by $\alpha^{-1}$, via the functoriality of $H^2$) 
and the vertical arrows are 
isomorphisms defined geometrically via the identification of $M_{X'}$ and $M_{X''}$ 
(resp. $M_{Y'}$ and $M_{Y''}$) 
with the Tate module of $\Bbb G_m$ over $\bar k_{X'}=\bar k_{X''}$ (resp. 
$\bar k_{Y'}=\bar k_{Y''}$). The commutativity of this last diagram follows from the fact that 
the above (geometrically defined) isomorphism $M_{X'}\isom M_{X`'}$ 
(resp. $M_{Y'}\isom M_{Y''}$) is identified with the composite of 
the $(X''_{\bar k_{X''}}:X'_{\bar k_{X'}})$-multiplication map $M_{X'}\to M_{X'}$ 
(resp. the $(Y''_{\bar k_{Y''}}:Y'_{\bar k_{Y'}})$-multiplication map $M_{Y'}\to M_{Y'}$) 
and the inverse of the natural isomorphism 
$M_{X''}\isom (X''_{\bar k_{X''}}:X'_{\bar k_{X'}})M_{X'}\subset M_{X'}$ (resp. 
$M_{Y''}\isom (Y''_{\bar k_{Y''}}:Y'_{\bar k_{Y'}})M_{Y'}\subset M_{Y'}$) induced by 
the inclusion $M_{X''}\to M_{X'}$ (resp. 
$M_{Y''}\to M_{Y'}$) 
arising from the functoriality of $H^2$, 
and the fact that 
$(X''_{\bar k_{X''}}:X'_{\bar k_{X'}})
=(\Delta_{U'}:\Delta_{U''})
=(\Delta_{V'}:\Delta_{V''})
=(Y''_{\bar k_{Y''}}:Y'_{\bar k_{Y'}})$. 
\qed
\enddemo

Thus, by Theorem 4.4, we obtain an isomorphism 
$X''\isom Y''$. Further, this isomorphism is Galois-equivariant with respect 
to $\alpha:\Pi_X\isom \Pi_Y$. 
Indeed, this Galois-equivariance can be proved just similarly as in 
the proof of Lemma 4.18(ii). 

Further, let $U''_1\to U$, $U''_2\to U$ be nice coverings of $U$ such that 
$\Pi_{U''_2}\subset \Pi_{U''_1}$, $V''_2\to V$, $V''_1\to V$ corresponding 
nice coverings of $V$ (via $\alpha: \Pi_U\isom \Pi_V$), and 
$X''_1, X''_2,Y''_1,Y''_2$ the smooth compactifications of 
$U''_1,U''_2,V''_1,V''_2$, respectively. Then the following diagram 
$$
\CD
X''_2 @>>>  Y''_2 \\
@VVV     @VVV  \\
X''_1  @>>> Y''_1   \\
\endCD
$$
is commutative, where the horizontal arrows are (Galois-equivariant) isomorphisms 
induced by $\alpha: \Pi_{X''_2}\isom\Pi_{Y''_2}$ and $\alpha: \Pi_{X''_1}\isom\Pi_{Y''_1}$, and 
the vertical arrows are natural (finite) morphisms induced by 
finite \'etale coverings $U''_2\to U''_1$ and $V''_2\to V''_1$. Indeed, the proof of this commutativity 
is similar to the proof of Lemma 4.20. 
More precisely, recall the proof of Theorem 4.4. The isomorphism $X''_i\isom Y''_i$ ($i\in\{1,2\}$) 
is induced by the isomorphism $R_{(Y''_i)^l}\isom R_{(X''_i)^l}$
obtained by applying Theorem 5.7 to 
$$\matrix
\Bbb P(R_{(X''_i)^l}) &&\Bbb P(R_{(Y''_i)^l}) \\
&&\\
\cup&&\cup\\
&&\\
R_{(X''_i)^l}^{\times}/(k_{X''_i}^{l})^\times &\isomleft&
R_{(Y''_i)^l}^{\times}/(k_{Y''_i}^{l})^\times
\endmatrix
$$
where 
$$
\matrix
\langle H_{(X''_i)^l}^\times\rangle&=&R_{(X''_i)^l}&\subset& \Cal O_{E_{(X''_i)^l}}\\
&&&&\\
\cup&&\cup&&\cup\\
&&&&\\
H_{(X''_i)^l}^\times&=&R_{(X''_i)^l}^{\times}&\subset& \Cal O_{E_{(X''_i)^l}}^\times
\endmatrix
$$
and 
$$
\matrix
\langle H_{(Y''_i)^l}^\times\rangle&=&R_{(Y''_i)^l}&\subset& \Cal O_{E_{(Y''_i)^l}}\\
&&&&\\
\cup&&\cup&&\cup\\
&&&&\\
H_{(Y''_i)^l}^\times&=&R_{(Y''_i)^l}^{\times}&\subset& \Cal O_{E_{(Y''_i)^l}}^\times
\endmatrix
$$

But the problem here (which does not occur in the proof of Lemma 4.20) is that 
in general the natural inclusions 
$K_{(X''_1)^l}\hookrightarrow K_{(X''_2)^l}$, 
$K_{(Y''_1)^l}\hookrightarrow K_{(Y''_2)^l}$ 
induced by the (ramified) coverings 
$f: X''_2\to X''_1$, 
$g: Y''_2\to Y''_1$, respectively, 
may not induce inclusions 
$\Cal O_{E_{(X''_1)^l}}\hookrightarrow \Cal O_{E_{(X''_2)^l}}$, 
$\Cal O_{E_{(Y''_1)^l}}\hookrightarrow \Cal O_{E_{(Y''_2)^l}}$, 
respectively. 
This is because it is unclear if 
$f^{-1}(E''_{X''_1})=E''_{X''_2}$, 
$g^{-1}(E''_{Y''_1})=E''_{Y''_2}$ hold. 

Here, the remedy is to resort to Proposition 2.11(ii) instead of Proposition 2.11(i). 
So, for each nice covering $U''\to U$ (resp. $V''\to V$), define 
$\Cal E_{X''}$ (resp. $\Cal E_{Y''}$) to be 
the inverse image of $E_{X'}\cup (X'\setminus U')$
(resp. $E_{Y'}\cup (Y'\setminus V')$) in $X''$ (resp. $Y''$). 
Then we have $E_{X''}\subset \Cal E_{X''}$, $E_{Y''}\subset\Cal E_{Y''}$. 
In particular, let $l$ be a $(X',\Sigma)$-admissible prime number 
(then $l$ is automatically $(Y',\Sigma)$-admissible). Then 
$l$ is $(X'',\Sigma)$-admissible and $(Y'',\Sigma)$-admissible 
for all nice coverings $U''\to U$ and $V''\to V$. Now, 
replacing $E_{X''}$, $E_{Y''}$ by $\Cal E_{X''}$, $\Cal E_{Y''}$ in the various definitions, 
we obtain 
$$
\matrix
\langle \Cal H_{(X'')^l}^\times\rangle&=&\Cal R_{(X'')^l}&\subset& \Cal O_{\Cal E_{(X'')^l}}\\
&&&&\\
\cup&&\cup&&\cup\\
&&&&\\
\Cal H_{(X'')^l}^\times&=&\Cal R_{(X'')^l}^{\times}&\subset&\Cal O_{\Cal E_{(X'')^l}}^\times
\endmatrix
$$
and 
$$
\matrix
\langle \Cal H_{(Y'')^l}^\times\rangle&=&\Cal R_{(Y'')^l}&\subset& \Cal O_{\Cal E_{(Y'')^l}}\\
&&&&\\
\cup&&\cup&&\cup\\
&&&&\\
\Cal H_{(Y'')^l}^\times&=&\Cal R_{(Y'')^l}^{\times}&\subset&\Cal O_{\Cal E_{(Y'')^l}}^\times
\endmatrix
$$
As in Lemma 4.13, we have 
$\Fr(\Cal R_{(X'')^l})=K_{(X'')^l}$ and 
$\Fr(\Cal R_{(Y'')^l})=K_{(Y'')^l}$. 
Then we may apply Theorem 5.7 to 
$$\matrix
\Bbb P(\Cal R_{(X'')^l}) &&\Bbb P(\Cal R_{(Y'')^l}) \\
&&\\
\cup&&\cup\\
&&\\
\Cal R_{(X'')^l}^{\times}/(k_{X''}^{l})^\times &\isomleft&
\Cal R_{(Y'')^l}^{\times}/(k_{Y''_i}^{l})^\times 
\endmatrix
$$
to obtain $\Cal R_{(Y'')^l}\isom \Cal R_{(X'')^l}$. By the uniqueness assertion of Theorem 5.7, 
the diagram 
$$
\CD
R_{(X'')^l} @<<< R_{(Y'')^l} \\
@AAA     @AAA  \\
\Cal R_{(X'')^l} @<<< \Cal R_{(Y'')^l} \\
\endCD
$$
commutes, where the vertical arrows are natural inclusions. In particular, 
$R_{(Y'')^l}\isom R_{(X'')^l}$ and 
$\Cal R_{(Y'')^l}\isom \Cal R_{(X'')^l}$ 
induce the same isomorphisms $K_{Y''}\isom K_{X''}$ and $X''\isom Y''$. 

Now, as in the proof of Lemma 4.20, we can prove that the diagrams 
$$
\CD
\Cal R_{(X''_2)^l} @<<< \Cal R_{(Y''_2)^l} \\
@AAA     @AAA  \\
\Cal R_{(X''_1)^l} @<<< \Cal R_{(Y''_1)^l} \\
\endCD
$$
$$
\CD
K_{X''_2} @<<< K_{Y''_2} \\
@AAA     @AAA  \\
K_{X''_1} @<<< K_{Y''_1} \\
\endCD
$$
and 
$$
\CD
X''_2 @>>> Y''_2 \\
@VVV     @VVV  \\
X''_1 @>>> Y''_1 \\
\endCD
$$
commute, as desired. 

Now, passing to nice open subgroups of $\Pi_U$ and $\Pi_V$ which correspond 
to each other via $\alpha$, we obtain an isomorphism 
$$\Tilde X_U \isom \Tilde Y_V $$ which is Galois-equivariant with respect to $\alpha:\Pi_U\isom\Pi_V$, 
where $\Tilde X_U $, $\Tilde Y_V $ are the integral closures of 
$X$, $Y$ in (the function fields of) $\Tilde U$, $\Tilde V$, respectively. 
(Note that in general $\Tilde X_U $, $\Tilde Y_V $ do not coincide with 
$\Tilde X$, $\Tilde Y$.) 
Further, dividing both sides of this isomorphism by the actions of 
$\Pi_U$ and $\Pi_V$, it follows that the isomorphism $\Tilde X_U \isom \Tilde Y_V $ 
fits into a commutative diagram  
$$
\CD
\Tilde X_U  @>{\sim}>> \Tilde Y_V  \\
@VVV   @VVV \\
X @>{\sim}>> Y \\
\endCD
$$
Finally, removing the ramification loci from this last diagram, 
we obtain the desired commutative diagram 
$$
\CD
\Tilde U @>{\sim}>> \Tilde V \\
@VVV   @VVV \\
U @>{\sim}>> V \\
\endCD
$$
in which the horizontal arrows are isomorphisms and the vertical arrows are the 
profinite \'etale coverings corresponding to the groups $\Pi_U$, $\Pi_V$, 
respectively.

This finishes the proof of Theorem 4.23. \qed

\subhead
\S 5. 
On the fundamental theorem of projective geometry
\endsubhead
Throughout this section all fields are assumed to be commutative. 
For a field $k$ and a vector space $V$ over $k$, 
define $\Bbb P(V)$ to be the projective space associated to $V$: 
$$\Bbb P(V)\defeq (V\setminus\{0\})/k^{\times},$$
and define $\Bbb L(V)\ (\subset 2^{\Bbb P(V)})$ to be the set of 
lines on $\Bbb P(V)$. 
Thus, $\Bbb P(V)$ (resp. $\Bbb L(V)$) is 
also identified with the set of $1$-dimensional (resp. $2$-dimensional) 
$k$-vector subspaces of $V$. For each $x\in V\setminus\{0\}$, 
denote the point of $\Bbb P(V)$ corresponding to $x$ by $\bar x$. 
We say that a set of points in $\Bbb P(V)$ is collinear, if 
there exists a line in $\Bbb L(V)$ which contains all of them. 
We say that a set of lines in $\Bbb L(V)$ is concurrent, 
if there exists a point in $\Bbb P(V)$ which is contained in all of them. 
We set $\Bbb P^n(k)\defeq \Bbb P(k^n)$ for each $n\geq 0$. 

It is well-known and easily proved 
that the projective space $(\Bbb P(V),\Bbb L(V))$ satisfies 
the following proposition. 

\proclaim{Proposition 5.1 (cf.~[EDM], 343)} 
Let $k$ be a field and $V$ a vector space over $k$. 

\noindent
{\rm (i)} (Axioms of projective geometry) The following (I)-(III) hold.  

\noindent
{\rm (I)} If $p,q\in\Bbb P(V)$, $p\neq q$, then there exists a unique 
$\ell\in\Bbb L(V)$, such that 
$p,q\in\Bbb \ell$. We denote this line $\ell$ by $p\vee q$. 

\noindent
{\rm (II)} If $p_0,p_1,p_2,q_1,q_2\in\Bbb P(V)$, 
$p_0,p_1,p_2$ are not collinear, 
$q_1\neq q_2$, 
$p_0,p_1,q_1$ are collinear and 
$p_0,p_2,q_2$ are collinear, 
then 
$p_1\vee p_2, q_1\vee q_2$ are concurrent. 

\noindent
{\rm (III)} If $\ell\in\Bbb L(V)$, then $\sharp(\ell)\geq 3$. 

\noindent
{\rm (ii)} (Desargues' theorem) 
If $p_1,p_2,p_3,q_1,q_2,q_3\in\Bbb P(V)$, 
$p_1,p_2,p_3$ are not collinear, 
$q_1,q_2,q_3$ are not collinear, 
$p_1\neq q_1$, 
$p_2\neq q_2$ and 
$p_3\neq q_3$, 
then: 
``$p_1\vee q_1,p_2\vee q_2,p_3\vee q_3$ are concurrent''
if and only if 
``$(p_2\vee p_3)\cap(q_2\vee q_3)$, 
$(p_3\vee p_1)\cap(q_3\vee q_1)$, 
$(p_1\vee p_2)\cap(q_1\vee q_2)$ 
are collinear''. \qed
\endproclaim 


Now, roughly speaking, the fundamental theorem of projective geometry asserts 
that the information carried by the pair $(k,V)$ is equivalent to 
that carried by the pair $(\Bbb P(V), \Bbb L(V))$. 

{}To be more precise, from now on, let $k_i$ be a field and $V_i$ a 
vector space over $k_i$, for $i=1,2$. 

\definition
{Definition 5.2}
(i) An (A semilinear) isomorphism $(k_1,V_1)\isom (k_2,V_2)$ 
is a pair $(\mu,\lambda)$, 
where $\mu$ is an isomorphism $k_1\isom k_2$ of fields and 
$\lambda$ is an isomorphism $V_1\isom V_2$ of abelian groups, 
such that, 
for each $a\in k_1$ and $x\in V_1$, one has 
$\lambda(ax)=\mu(a)\lambda(x)$. 
(In fact, when $V_i\neq 0$, $\mu$ is determined uniquely by $\lambda$, 
hence we may say that $\lambda$ is an (a semilinear) isomorphism.) 

\noindent
(ii) A collineation $(\Bbb P(V_1), \Bbb L(V_1))\isom (\Bbb P(V_2), \Bbb L(V_2))$ 
is a pair $(\sigma,\tau)$, 
where $\sigma$ is a bijection $\Bbb P(V_1)\isom\Bbb P(V_2)$ 
and 
$\tau$ is a bijection $\Bbb L(V_1)\isom\Bbb L(V_2)$, 
such that, for each $\ell\in\Bbb L(V_1)$, one has 
$\tau(\ell)=\sigma(\ell)(\defeq \{\sigma(p)\mid p\in\ell\})$. 
(In fact, $\tau$ is determined uniquely by $\sigma$, hence 
we may say that $\sigma$ is a collineation.) 
\enddefinition

\proclaim{Theorem 5.3 (Fundamental Theorem of Projective Geometry, cf.~[Artin])}

\noindent
{\rm (i)} Each isomorphism $(\mu,\lambda): (k_1,V_1)\isom (k_2,V_2)$ 
naturally induces a collineation 
$(\sigma, \tau): (\Bbb P(V_1), \Bbb L(V_1))\isom (\Bbb P(V_2), \Bbb L(V_2))$ 
by setting $\sigma(\bar x)\defeq\overline{\lambda(x)}$ for 
$x\in V_1\setminus\{0\}$. 

\noindent
{\rm (ii)} Assume that $\dim_{k_i}(V_i)\geq 3$ for $i=1,2$. Then, 
for each collineation 
$(\sigma, \tau): (\Bbb P(V_1), \Bbb L(V_1))\isom (\Bbb P(V_2), \Bbb L(V_2))$, 
there exists an isomorphism $(\mu,\lambda): (k_1,V_1)\isom (k_2,V_2)$ 
that induces $(\sigma, \tau): (\Bbb P(V_1), \Bbb L(V_1))\isom (\Bbb P(V_2), \Bbb L(V_2))$ 
(in the sense of (i)). Further, such an isomorphism $(\mu,\lambda)$ is unique 
up to scalar multiplication. More precisely, if $(\mu,\lambda),(\mu',\lambda')$ 
are such isomorphisms (that induce the same collineation $(\sigma, \tau)$), then 
there exists an (in fact, a unique) element $a\in k_1^{\times}$ such that 
$\mu'=\mu$ and $\lambda'(-)=\lambda(a\cdot-)$. \qed
\endproclaim

The aim of this section is to give a refined version of the fundamental theorem 
of projective geometry, where certain ``partial'' collineations defined over 
``sufficiently large'' subsets of projective spaces are considered. 
To formulate it precisely, let us first define what are ``sufficiently large'' 
subsets of projective spaces. 

\definition{Definition 5.4} 
Let $k$ be a field. Let $\Cal S$ be a set of subsets of $\Bbb P^1(k)$. 

\noindent
(i) We say that $\Cal S$ is $\PGL_2$-stable, if, for any $S\in\Cal S$ and 
any $\sigma\in\PGL_2(k)$, one has $\sigma(S)\in\Cal S$. 

\noindent
(ii) Let $m,n$ be integers $\geq 0$. Then we say that 
$\Cal S$ is $(m,n)$-admissible, if $\Cal S$ is 
$\PGL_2$-stable and, 
for 
any $0\leq m'\leq m$, $0\leq n'\leq n$, 
any $S_1,\dots,S_{m'}\in\Cal S$, 
and 
any $p_1,\dots,p_{n'}\in\Bbb P^1(k)$, 
one has 
$$S_1\cup\dots\cup S_{m'}\cup\{p_1,\dots,p_{n'}\}\subsetneq\Bbb P^1(k).$$
(Thus, if $m_1\geq m_2\geq 0$, $n_1\geq n_2\geq 0$, then
 $\Cal S$: $(m_1,n_1)$-admissible $\implies$ $\Cal S$: $(m_2,n_2)$-admissible.) 
)

\noindent
(iii) We say that $\Cal S$ is admissible, if $\Cal S$ is $(m,n)$-admissible 
for all integers $m,n\geq 0$. 

\noindent
(iv) Assume that $\Cal S$ is $\PGL_2$-stable. Then, for each $1$-dimensional 
projective space $\ell$ over $k$, we set 
$$\Cal S_{\ell}\defeq \{S\subset \ell\mid\text{$\alpha(S)\in \Cal S$ for some 
$\alpha:\ell\underset{k}\to{\isom} \Bbb P^1(k)$}\}.$$ 
(By assumption, we see that ``for some $\alpha$'' in this definition may be 
replaced by ``for all $\alpha$'' and that $\Cal S= \Cal S_{\Bbb P^1(k)}$.) 
\enddefinition

\definition{Remark 5.5} 
We have: 
$$
\text{$\Cal S$: $(m,n)$-admissible} 
\Longleftarrow 
\forall S\in \Cal S, \sharp(\Bbb P^1(k))>m\sharp(S)+n 
\Longleftrightarrow 
\sharp (k)> m\sharp (S)+n-1.
$$
In particular, if $\Cal S=\{\varnothing\}$, we have:
$$
\text{$\Cal S$: $(m,n)$-admissible} 
\Longleftrightarrow 
\Bbb \sharp(P^1(k))>n 
\Longleftrightarrow 
\sharp (k)\ge n.
$$
\enddefinition

\definition{Definition 5.6}
Let $k$ be a field, $V$ a vector space over $k$, and $U$ a subset of $\Bbb P(V)$. 

\noindent
(o) For each line $\ell\in\Bbb L(V)$, we write 
$\ell_U\defeq \ell\cap U$ and 
$\ell_{U^c}\defeq \ell\cap(\Bbb P(V)\setminus U)=\ell\setminus\ell_U$. 

\noindent
(i) We define $\Bbb L(V)_{U}\subset \Bbb L(V)$ by: 
$$\Bbb L(V)_{U}\defeq\{\ell\in\Bbb L(V)\mid \ell_U\neq\varnothing\}.$$

\noindent
(ii) Let $\Cal S$ be a $\PGL_2$-stable set of subsets of $\Bbb P^1(k)$. 
Then we define $\Bbb L(V)_{U,\Cal S}\subset \Bbb L(V)$ by: 
$$\Bbb L(V)_{U,\Cal S}\defeq\{\ell\in\Bbb L(V)\mid \ell_{U^c}
\in \Cal S_{\ell}\}.$$

\noindent
(iii) Let $\Cal S$ be a $\PGL_2$-stable set of subsets of $\Bbb P^1(k)$. 
We say that $U$ is $\Cal S$-ample, if the following conditions (1)(2) hold. 

\noindent
(1) $U\neq\varnothing$. 

\noindent
(2) $\Bbb L(V)_{U}\subset \Bbb L(V)_{U,\Cal S}$. Equivalently, 
for each $\ell\in\Bbb L(V)$, either $\ell_U=\varnothing$ 
or $\ell_{U^c} \in \Cal S_{\ell}$. 
(When $\Cal S$ is $(1,0)$-admissible, one automatically has 
$\Bbb L(V)_{U,\Cal S}\subset \Bbb L(V)_{U}$, and the 
above condition (2) is then equivalent to: 
$\Bbb L(V)_{U}= \Bbb L(V)_{U,\Cal S}$.) 
\enddefinition 

Now, return to the situation of the fundamental theorem of projective geometry. 
Namely, let $k_i$ be a field and $V_i$ a vector space over $k_i$, for $i=1,2$. 
The main result in this section is the following refinement of the fundamental 
theorem of projective geometry. 

\proclaim{Theorem 5.7}
%
Assume that $\dim_{k_i}(V_i)\geq 3$ for $i=1,2$. 
Let $U_i$ be a subset of $\Bbb P(V_i)$ for $i=1,2$, and assume that 
$U_i$ is $\Cal S_i$-ample for some 
$(3,2)$-admissible 
set $\Cal S_i$ of subsets of $\Bbb P^1(k_i)$ for $i=1,2$. 
Let $\sigma: U_1\isom U_2$ and $\tau: \Bbb L(V_1)_{U_1}\isom \Bbb L(V_2)_{U_2}$ 
be bijections such that 
for each $\ell\in\Bbb L(V_1)_{U_1}$, one has 
$\tau(\ell)_{U_2}= \sigma(\ell_{U_1})$. 
Then, each such 
$(\sigma, \tau): (U_1, \Bbb L(V_1)_{U_1})\isom (U_2, \Bbb L(V_2)_{U_2})$ 
uniquely extends to a collineation 
$(\tilde\sigma, \tilde\tau): (\Bbb P(V_1), \Bbb L(V_1))\isom (\Bbb P(V_2), \Bbb L(V_2))$. 
In particular, 
there exists an isomorphism $(\mu,\lambda): (k_1,V_1)\isom (k_2,V_2)$ 
that induces $(\sigma, \tau): (U_1, \Bbb L(V_1)_{U_1})\isom 
(U_2, \Bbb L(V_2)_{U_2})$, and 
such an isomorphism $(\mu,\lambda)$ is unique 
up to scalar multiplication. 
\endproclaim

\demo{Proof}
{\it Step 0.} The second assertion follows from the first assertion and 
Theorem 5.3. So, let us concentrate on the proof of the first assertion. 

\medskip\noindent
{\it Step 1.} Claim: If $p\in U_1$, then $\tau: \Bbb L(V_1)_{U_1}\isom\Bbb L(V_2)_{U_2}$ induces 
a bijection of subsets $\Bbb L(V_1)_{\{p\}}\isom\Bbb L(V_2)_{\{\sigma(p)\}}$. 

Indeed, for each $\ell\in\Bbb L(V_1)_{U_1}$, one has 
$$
\align
\ell\in\Bbb L(V_1)_{\{p\}} &\iff p\in\ell \\
&\iff p\in\ell_{U_1} \\
&\iff \sigma(p)\in\sigma(\ell_{U_1})=\tau(\ell)_{U_2} \\
&\iff \sigma(p)\in\tau(\ell) \\
&\iff \tau(\ell)\in \Bbb L(V_2)_{\{\sigma(p)\}}, 
\endalign
$$
as desired. 

\medskip\noindent
{\it Step 2.}  Claim: If $p\in \Bbb P(V_1)$, 
then there exists a unique point $p'\in \Bbb P(V_2)$, 
such that $\tau: \Bbb L(V_1)_{U_1}\isom\Bbb L(V_2)_{U_2}$ induces 
a bijection of subsets 
$\Bbb L(V_1)_{\{p\}}\cap \Bbb L(V_1)_{U_1}
\isom\Bbb L(V_2)_{\{p'\}}\cap \Bbb L(V_2)_{U_2}$. 
If, moreover, $p\in U_1$, then $p'=\sigma(p)$. 

\medskip\noindent
{\it Step 2-1.} Claim: 
Two lines $\ell_1,\ell_2\in\Bbb L(V_1)_{U_1}$ are concurrent, if and only if 
so are $\tau(\ell_1),\tau(\ell_2)\in\Bbb L(V_2)_{U_2}$. 

Indeed, it suffices to prove the `only if' part, since 
the `if' part is obtained by applying the `only if' part 
to $\sigma^{-1}:U_2\isom U_1$. 
Now, the `only if' part is clear if $\ell_1=\ell_2$. So, assume $\ell_1\neq\ell_2$. 
Then, as $\ell_1,\ell_2$ 
are concurrent, there is a (unique) point $p\in\Bbb P(V)$ such that 
$\ell_1\cap\ell_2=\{p\}$. 
If $p\in U_1$, then $\sigma(p)\in\tau(\ell_1), \tau(\ell_2)$ 
by Step 1, hence $\tau(\ell_1), \tau(\ell_2)$ are concurrent, 
as desired. So, we may and shall assume that $p\not\in U_1$. 
For each $i=1,2$, choose 
$p_i\in (\ell_i)_{U_1}
=\ell_i\setminus((\ell_i)_{U_1^c}
)\nemp$ 
($(1,0)$-admissibility). 
As $p_1\in\ell_1\setminus\ell_2$ and $p_2\in\ell_2\setminus\ell_1$, one has 
$p_1\neq p_2$. So, 
set $m\defeq p_1\vee p_2$. 
As $p_1\in U_1$, one has $m\in\Bbb L(V_1)_{U_1}$. Next, take 
$q\in m_{U_1}\setminus\{p_1,p_2\}=m\setminus(m_{U_1^c}\cup\{p_1,p_2\})\nemp$ 
($(1,2)$-admissibility). 
Consider the projection (or perspective mapping) $\alpha: \ell_1\isom\ell_2$ with respect to the 
center $q$. More precisely, $\alpha$ is defined by $\{\alpha(x)\}=(x\vee q)\cap\ell_2$ 
for each $x\in\ell_1$. (In particular, $\alpha(p_1)=p_2$ and $\alpha(p)=p$.) 
Take $q_1\in((\ell_1)_{U_1}\cap \alpha^{-1}((\ell_2)_{U_1}))
\setminus\{
p_1\}
=\ell_1\setminus((\ell_1)_{U_1^c}\cup\alpha^{-1}((\ell_2)_{U_1^c})
\cup\{
p_1\})\nemp$ 
($(2,1)$-admissibility), 
and set $q_2\defeq\alpha(q_1)$ and $n=q_1\vee q$ ($=q_2\vee q=q_1\vee q_2$). (Thus, $q_2\in 
((\ell_2)_{U_1}\cap \alpha((\ell_1)_{U_1}))
\setminus\{
p_2\}$, and, as 
$q_1\in U_1$, one has $n\in\Bbb L(V_1)_{U_1}$. 
) 
Now, one has 
$q, p_1,q_1,p_2,q_2\in U_1$, 
$q, p_1, q_1$ are not collinear, 
$p_2\neq q_2$, 
$q, p_1,p_2\in m$ and 
$q, q_1,q_2\in n$. 
Accordingly, one has 
$\sigma(q),\sigma(p_1),\sigma(q_1),\sigma(p_2),\sigma(q_2)\in U_2$, 
$\sigma(q), \sigma(p_1), \sigma(q_1)$ not collinear, 
$\sigma(p_2)\neq \sigma(q_2)$, 
and, by Step 1, 
$\sigma(p_1),\sigma(p_2), \sigma(q)\in \tau(m)$ and 
$\sigma(q_1),\sigma(q_2),\sigma(q)\in \tau(n)$. 
Further, as $p_1,q_1\in \ell_1$ ($p_1\neq q_1$) and 
$p_2,q_2\in \ell_2$ ($p_2\neq q_2$), one has 
$\sigma(p_1), \sigma(q_1)\in\tau(\ell_1)$
(with $\sigma(p_1)\neq\sigma(q_1)$) and 
$\sigma(p_2), \sigma(q_2)\in\tau(\ell_2)$ 
(with $\sigma(p_2)\neq\sigma(q_2)$) 
by Step 1. 
Now, by Proposition 5.1(i)(II), this implies that 
$\tau(\ell_1)=\sigma(p_1)\vee\sigma(q_1)$ and 
$\tau(\ell_2)=\sigma(p_2)\vee\sigma(q_1)$ are concurrent, 
as desired. 

\medskip\noindent
{\it Step 2-2.} Claim: Three lines $\ell_1,\ell_2,\ell_3\in\Bbb L(V_1)_{U_1}$ are concurrent, 
if and only if 
so are $\tau(\ell_1),\tau(\ell_2), 
\tau(\ell_3)\in\Bbb L(V_2)_{U_2}$. 

Indeed, it suffices to prove the `only if' part, since 
the `if' part is obtained by applying the `only if' part 
to $\sigma^{-1}:U_2\isom U_1$. 
Now, the `only if' part 
follows from (the `only if' part of) Step 2-1 if $\ell_i=\ell_j$ 
for some $i\neq j$. 
So, assume that $\ell_1,\ell_2,\ell_3$ are mutually distinct and 
set $\ell_1\cap\ell_2\cap\ell_3=\{p\}$. 
If $p\in U_1$, then $\sigma(p)\in\tau(\ell_1), \tau(\ell_2), \tau(\ell_3)$ 
by Step 1, hence $\tau(\ell_1), \tau(\ell_2), \tau(\ell_3)$ are concurrent, 
as desired. So, we may and shall assume that $p\not\in U_1$. 
Take $p_1\in(\ell_1)_{U_1}\nemp$ 
and 
$p_2\in(\ell_2)_{U_1}\nemp$. 
Take $p_3\in(\ell_3)_{U_1}\setminus((p_1\vee p_2)\cap\ell_3)=
\ell_3\setminus((\ell_3)_{U_1^c}\cup((p_1\vee p_2)\cap\ell_3))\nemp$ 
($(1,1)$-admissibility). 
Take $q_1\in (\ell_1)_{U_1}\setminus\{p_1\}=
\ell_1\setminus((\ell_1)_{U_1^c}\cup\{p_1\})\nemp$ 
($(1,1)$-admissibility). 
Let $\alpha: \ell_2\isom (p_1\vee p_2)$ be the projection 
with respect to the center $q_1$: 
$\{\alpha(x)\}=(x\vee q_1)\cap(p_1\vee p_2)$. 
Take $q_2\in 
(\ell_2)_{U_1}\cap \alpha^{-1}((p_1\vee p_2)_{U_1})\setminus\{p_2\}=
\ell_1\setminus((\ell_1)_{U_1^c}\cup \alpha^{-1}((p_1\vee p_2)_{U_1^c})
\cup\{p_2\})\nemp$ 
($(2,1)$-admissibility), 
and set 
$r_{12}\defeq \alpha(q_2)
\in (p_1\vee p_2)_{U_1}\setminus\{p_1,p_2\}$. 
Let $\beta: \ell_3\isom (p_2\vee p_3)$ be the projection 
with respect to the center $q_2$: 
$\{\beta(x)\}=(x\vee q_2)\cap(p_2\vee p_3)$ 
and 
$\gamma: \ell_3\isom (p_3\vee p_1)$ the projection 
with respect to the center $q_1$: 
$\{\gamma(x)\}=(x\vee q_1)\cap(p_3\vee p_1)$. 
Take $q_3\in
((\ell_3)_{U_1}\cap \beta^{-1}((p_2\vee p_3)_{U_1})
\cap\gamma^{-1}((p_3\vee p_1)_{U_1}))\setminus
(\{p_3\}\cup((q_1\vee q_2)\cap\ell_3))
=\ell_3\setminus
((\ell_3)_{U_1^c}\cup \beta^{-1}((p_2\vee p_3)_{U_1^c})
\cup\gamma^{-1}((p_3\vee p_1)_{U_1^c})\cup\{p_3\}
\cup((q_1\vee q_2)\cap\ell_3))\nemp$ 
(here, we use the $(3,2)$-admissibility assumption fully), 
and set 
$r_{23}\defeq \beta(q_3)
\in (p_2\vee p_3)_{U_1}\setminus\{p_2,p_3\}$ 
and 
$r_{31}\defeq \gamma(q_3)
\in (p_3\vee p_1)_{U_1}\setminus\{p_3,p_1\}$. 
Now, one has 
$p_1,p_2,p_3,q_1,q_2,q_3,r_{12}, r_{23}, r_{31}\in U_1$ 
with $p_1,p_2,p_3$ not collinear and $q_1,q_2,q_3$ not collinear, 
$p_1,q_1\in \ell_1$ with $p_1\neq q_1$, 
$p_2,q_2\in \ell_2$ with $p_2\neq q_2$, 
$p_3,q_3\in \ell_3$ with $p_3\neq q_3$, 
$\{r_{12}\}=(p_1\vee p_2)\cap(q_1\vee q_2)$, 
$\{r_{23}\}=(p_2\vee p_3)\cap(q_2\vee q_3)$, 
$\{r_{31}\}=(p_3\vee p_1)\cap(q_3\vee q_1)$, and 
$\ell_1=p_1\vee q_1$,  
$\ell_2=p_2\vee q_2$, 
$\ell_3=p_3\vee q_3$ are concurrent. 
Then, by (the `$\implies$' part of) Proposition 5.1(ii), $r_{12}, r_{23}, r_{31}$ are 
collinear. 
Accordingly (by Step 1), one has 
$\sigma(p_1),\sigma(p_2),\sigma(p_3),\sigma(q_1),\sigma(q_2),\sigma(q_3),
\sigma(r_{12}), \sigma(r_{23}), \sigma(r_{31})\in U_2$ 
with $\sigma(p_1),\sigma(p_2),\sigma(p_3)$ not collinear and $\sigma(q_1),\sigma(q_2),\sigma(q_3)$ 
not collinear, 
$\sigma(p_1),\sigma(q_1)\in \tau(\ell_1)$ with $\sigma(p_1)\neq \sigma(q_1)$, 
$\sigma(p_2),\sigma(q_2)\in \tau(\ell_2)$ with $\sigma(p_2)\neq \sigma(q_2)$,  
$\sigma(p_3),\sigma(q_3)\in \tau(\ell_3)$ with $\sigma(p_3)\neq \sigma(q_3)$, 
$\{\sigma(r_{12})\}=(\sigma(p_1)\vee \sigma(p_2))\cap(\sigma(q_1)\vee \sigma(q_2))$, 
$\{\sigma(r_{23})\}=(\sigma(p_2)\vee \sigma(p_3))\cap(\sigma(q_2)\vee \sigma(q_3))$,
$\{\sigma(r_{31})\}=(\sigma(p_3)\vee \sigma(p_1))\cap(\sigma(q_3)\vee \sigma(q_1))$, and 
$\sigma(r_{12}), \sigma(r_{23}), \sigma(r_{31})$ are collinear. Now, 
by (the `$\impliedby$' part of) Proposition 5.1(ii), 
$\tau(\ell_1)=\sigma(p_1)\vee\sigma(q_1)$, 
$\tau(\ell_2)=\sigma(p_2)\vee\sigma(q_2)$, 
$\tau(\ell_3)=\sigma(p_3)\vee\sigma(q_3)$ 
are concurrent, as desired. 

\medskip\noindent
{\it Step 2-3.} Claim: If $p\in\Bbb P(V_i)$ for $i=1,2$ , one has 
$$\bigcap_{\ell\in \Bbb L(V_i)_{\{p\}}\cap \Bbb L(V_i)_{U_i}}\ell =\{p\}.$$ 

Indeed, we may assume that $i=1$. First, 
``$\supset$'' is clear. So, to prove ``$\subset$'', it suffices to 
show that the left-hand side is of cardinality at most one. Then, since two distinct 
lines intersect at at most one point, it suffices to prove that 
there are at least two elements  
$\ell\in \Bbb L(V_1)_{\{p\}}\cap \Bbb L(V_1)_{U_1}$. 
Consider the two cases separately: (i) $p\in U_1$; and (ii) $p\not\in U_1$. 
In case (i), as $\dim(\Bbb P(V_1))\geq 2$, hence there exist $q,r\in \Bbb P(V_1)$ 
such that $p,q,r$ are not collinear. Then $p\vee q$ and $p\vee r$ 
are two distinct lines that belong to $\Bbb L(V_1)_{\{p\}}\cap \Bbb L(V_1)_{U_1}$. 
In case (ii), as $U_1\neq\varnothing$, take $q\in U_1$. (Thus, $q\neq p$.) 
As $\dim(\Bbb P(V_1))\geq 2$, hence there exists $r\in \Bbb P(V_1)$ 
such that $p,q,r$ are not collinear. Observe 
$q\vee r\in \Bbb L(V_1)_{U_1}$, and take 
$s\in (q\vee r)_{U_1}\setminus\{q,r\}
=(q\vee r)\setminus((q\vee r)_{U_1^c}\cup\{q,r\})\nemp$ 
($(1,2)$-admissibility). 
Now, $p\vee q$ and $p\vee s$ 
are two distinct lines that belong to $\Bbb L(V_1)_{\{p\}}\cap \Bbb L(V_1)_{U_1}$. 

\medskip\noindent
{\it Step 2-4.} Claim: If $p\in\Bbb P(V_1)$, 
$$\bigcap_{\ell\in \Bbb L(V_1)_{\{p\}}\cap \Bbb L(V_1)_{U_1}}\tau(\ell)$$ 
is a subset of $\Bbb P(V_2)$ of cardinality one. (Denote it by $\{p'\}$.) 

Indeed, for each pair $\ell,m\in \Bbb L(V_1)_{\{p\}}\cap \Bbb L(V_1)_{U_1}$ with 
$\ell\neq m$, $\tau(\ell)\cap\tau(m)$ is of cardinality one by Step 2-1. So, 
set $\tau(\ell)\cap\tau(m)=\{p'_{\ell,m}\}$. In fact, 
the point $p'_{\ell,m}$ does not depend on the pair $\ell,m$. Indeed, 
let $\ell',m'\in \Bbb L(V_1)_{\{p\}}\cap \Bbb L(V_1)_{U_1}$ be any pair with 
$\ell'\neq m'$. 
If $\sharp\{\ell,m,\ell',m'\}=2$, 
then it is clear that $p'_{\ell,m}=p'_{\ell',m'}$. 
If $\sharp\{\ell,m,\ell',m'\}=3$, 
then it follows from Step 2-2 that $\tau(\ell), \tau(m), \tau(\ell'), \tau(m')$ 
are concurrent, hence $p'_{\ell,m}=p'_{\ell',m'}$. 
If $\sharp\{\ell,m,\ell',m'\}=4$, 
then, again by Step 2-2, one has 
$$p'_{\ell,m}=p'_{\ell,m'}=p'_{\ell',m'}.$$
Now, write $p'=p'_{\ell,m}$ for some (or, equivalently, all) pair 
$\ell,m\in \Bbb L(V_1)_{\{p\}}\cap \Bbb L(V_1)_{U_1}$ with 
$\ell\neq m$. Here, note that, 
as shown in the proof of Step 2-3, 
$\Bbb L(V_1)_{\{p\}}\cap \Bbb L(V_1)_{U_1}$ is of cardinality at least two, 
hence at least one such pair $\ell,m$ exists. Then, by definition, 
one has 
$$\bigcap_{\ell\in \Bbb L(V_1)_{\{p\}}\cap \Bbb L(V_1)_{U_1}}\tau(\ell)=\{p'\},$$ 
as desired. 

\medskip\noindent
{\it Step 2-5.} End of Step 2. 

Let $p\in\Bbb L(V_1)$ and define $p'\in\Bbb L(V_2)$ as in Step 2-4. 
Then one has 
$$\tau(\Bbb L(V_1)_{\{p\}}\cap \Bbb L(V_1)_{U_1})
\subset \Bbb L(V_2)_{\{p'\}}\cap \Bbb L(V_2)_{U_2}.$$ 
Applying this to 
$\sigma^{-1}:U_2\isom U_1$ and $p'\in\Bbb L(V_2)$, we obtain 
$$
\tau^{-1}(\Bbb L(V_2)_{\{p'\}}\cap \Bbb L(V_2)_{U_2}
)
\subset 
\Bbb L(V_1)_{\{p''\}}\cap \Bbb L(V_1)_{U_1}
$$  
for some unique $p''\in\Bbb L(V_1)$. Combining these containment 
relations, we conclude 
$$\Bbb L(V_1)_{\{p\}}\cap \Bbb L(V_1)_{U_1}
\subset 
\Bbb L(V_1)_{\{p''\}}\cap \Bbb L(V_1)_{U_1},$$ 
which, together with Step 2-3, implies that 
$\{p\}\supset\{p''\}$, hence 
$p=p''$, and that 
$$\tau(\Bbb L(V_1)_{\{p\}}\cap \Bbb L(V_1)_{U_1})= \Bbb L(V_2)_{\{p'\}}\cap \Bbb L(V_2)_{U_2},$$
as desired. 

The uniqueness of $p'$ is clear by Step 2-3. This uniqueness, together with Step 1, implies 
$p'=\sigma(p)$ for $p\in U_1$. 

\medskip\noindent
{\it Step 3.} We define $\tilde\sigma: \Bbb P(V_1)\to\Bbb P(V_2)$ to be the 
map that sends $p\in \Bbb P(V_1)$ to $p'\in \Bbb P(V_2)$ defined in Step 2. 
(Thus, in particular, $\tilde\sigma(p)=\sigma(p)$ if $p\in U_1$.) 

\medskip\noindent
Claim: $\tilde\sigma$ is a collineation. (More precisely, there exists a (unique) bijection  
$\tilde\tau: \Bbb L(V_1)\isom \Bbb L(V_2)$, such that $(\tilde\sigma,\tilde\tau)$ 
is a collineation.) 

\medskip\noindent
{\it Step 3-1.} Claim: $\tilde\sigma: \Bbb P(V_1)\to\Bbb P(V_2)$ is a bijection. 

Indeed, starting with $\sigma^{-1}:U_2\isom U_1$ instead of $\sigma$, we obtain 
$\widetilde{\sigma^{-1}}: \Bbb P(V_2)\to\Bbb P(V_1)$, which turns out to 
be the inverse of $\tilde\sigma$ from the uniqueness assertion of Step 2. 

\medskip\noindent
{\it Step 3-2.} Claim: If $\ell\in\Bbb L(V_1)_{U_1}$, then $\tilde\sigma(\ell)=\tau(\ell)\ 
(\in \Bbb L(V_2)_{U_2})$. 

Indeed, let $\ell\in\Bbb L(V_1)_{U_1}$ and $p\in\ell$. Then, just 
by the definition of $\tilde\sigma$, we have $\tilde\sigma(p)\in\tau(\ell)$. 
Namely, we obtain $\tilde\sigma(\ell)\subset \tau(\ell)$. Applying this to 
$\sigma^{-1}: U_2\isom U_1$ and $\tau(\ell)\in \Bbb L(V_2)_{U_2}$ (and 
noting that $\widetilde{\sigma^{-1}}=\tilde\sigma^{-1}$ as shown in 
Step 3-1), we obtain $\tilde\sigma^{-1}(\tau(\ell))\subset
\tau^{-1}(\tau(\ell))=\ell$, or, equivalently, 
$\tau(\ell)\subset \tilde\sigma(\ell)$. Combining these, we obtain 
$\tilde\sigma(\ell)= \tau(\ell)$, as desired. 

\medskip\noindent
{\it Step 3-3.} Claim: Three points $p_1,p_2,p_3\in\Bbb P(V_1)$ are collinear, 
if and only if so are $\tilde\sigma(p_1), \tilde\sigma(p_2), \tilde\sigma(p_3)$. 

Indeed, it suffices to prove the `only if' part, since 
the `if' part is obtained by applying the `only if' part 
to $\sigma^{-1}:U_2\isom U_1$. 
If $\sharp\{p_1,p_2,p_3\}\leq 2$, the assertion is clear. So, we may assume 
that $p_1,p_2,p_3$ are mutually distinct. In particular, the line $\ell$ 
containing $p_1,p_2,p_3$ (whose existence is ensured by the collinearity of 
$p_1,p_2,p_3$) is unique. Next, if $\ell\in\Bbb L(V_1)_{U_1}$, then 
the assertion follows immediately from Step 3-2. So, we may assume that 
$\ell\not\in \Bbb L(V_1)_{U_1}$, i.e., $\ell\cap U_1=\varnothing$. 

Take $q_1\in U_1\nemp$. As $\ell\cap U_1=\varnothing$, we have $q_1\not\in \ell$, 
hence there exists a unique plane $P\subset\Bbb P(V_1)$ containing both $\ell$ 
and $q_1$. We shall construct various points in $P\cap U_1$. 
Let $\alpha: (q_1\vee p_3)\isom (q_1\vee p_2)$ be the projection 
with respect to the center $p_1$:  
$\{\alpha(x)\}=(x\vee p_1)\cap(q_1\vee p_2)$. 
Take $q_2\in 
(q_1\vee p_3)_{U_1}\cap\alpha^{-1}((q_1\vee p_2)_{U_1})\setminus\{q_1\}
=(q_1\vee p_3)\setminus ((q_1\vee p_3)_{U_1^c}\cup 
\alpha^{-1}((q_1\vee p_2)_{U_1^c})\cup\{q_1\})\nemp$ 
($(2,1)$-admissibility), 
and set 
$q_3\defeq\alpha(q_2)\in 
(q_1\vee p_2)_{U_1}\cap \alpha((q_1\vee p_3)_{U_1})\setminus\{q_1\}$. 
Next, take $r_1\in P\cap U_1\setminus((q_1\vee p_3)\cup(q_1\vee p_2)) 
\supset (q_1\vee p_1)_{U_1}\setminus\{q_1\}
=(q_1\vee p_1)\setminus ((q_1\vee p_1)_{U_1^c}\cup\{q_1\})\nemp$ 
($(1,1)$-admissibility). 
Let $\beta: (r_1\vee p_3)\isom (r_1\vee p_2)$ be the projection 
with respect to the center $p_1$:  
$\{\beta(x)\}=(x\vee p_1)\cap(r_1\vee p_2)$. 
Let $\gamma: (r_1\vee p_3)\isom (r_1\vee q_1)$ be the projection 
with respect to the center $q_2$:  
$\{\gamma(x)\}=(x\vee q_2)\cap(r_1\vee q_1)$. 
Take $r_2\in 
(r_1\vee p_3)_{U_1}\cap\beta^{-1}((r_1\vee p_2)_{U_1})\cap\gamma^{-1}((r_1\vee q_1)_{U_1})
\setminus
\{r_1\}
=(r_1\vee p_3)\setminus ((r_1\vee p_3)_{U_1^c}\cup 
\beta^{-1}((r_1\vee p_2)_{U_1^c})
\cup\gamma^{-1}((r_1\vee q_1)_{U_1^c})
\cup
\{r_1\}
)
\nemp$ 
($(3,1)$-admissibility), 
and set $r_3\defeq\beta(r_2)\in 
(r_1\vee p_2)_{U_1}
\setminus\{r_1\}$. 

First, $q_1,q_2,q_3$ are not collinear. Indeed, otherwise, $q_1,p_2,p_3$ must 
also be collinear, which contradicts the choice of $q_1$. 
Second, $r_1,r_2,r_3$ are not collinear. Indeed, otherwise, $r_1,p_3,p_2$ must 
also be collinear, which contradicts the choice of $r_1$. ($r_1\in U_1$ and 
$p_3\vee p_2=\ell\subset U_1^c$.)
Third, $q_1\neq r_1$. Indeed, this follows from the definition of $r_1$. 
Fourth, $q_2\neq r_2$. Indeed, otherwise, $q_2=r_2\in r_1\vee p_3$, hence 
$r_1\in q_2\vee p_3=q_1\vee p_3$, which contradicts the choice of $r_1$. 
Fifth, $q_3\neq r_3$. Indeed, otherwise, $q_3=r_3$, hence 
$r_1\in r_1\vee p_2=r_3\vee p_2=q_3\vee p_2=q_1\vee p_2$, 
which contradicts the choice of $r_1$. 

Thus, one may apply (the `$\impliedby$' part of) Proposition 5.1(ii) to 
$q_1,q_2,q_3,r_1,r_2,r_3$. Then 
there exists $s\in\Bbb P(V)$, such that 
$s, q_1,r_1$ are collinear, $s,q_2,r_2$ are collinear, and $s, q_3, r_3$
are collinear. 
By definition, $s=\gamma(r_2)\in U_1$. 

By Step 3-2, 
$\sigma(q_1),\sigma(q_2),\sigma(q_3)$ are not collinear, 
$\sigma(r_1),\sigma(r_2),\sigma(r_3)$ are not collinear, 
$\sigma(s),\sigma(q_1),\sigma(r_1)$ are collinear, 
$\sigma(s),\sigma(q_2),\sigma(r_2)$ are collinear, 
and 
$\sigma(s), \sigma(q_3), \sigma(r_3)$ are collinear. 
Also, as $\sigma$ is a bijection, one has
$\sigma(q_1)\neq\sigma(r_1)$, 
$\sigma(q_2)\neq\sigma(r_2)$ and 
$\sigma(q_3)\neq\sigma(r_3)$.  
Thus, applying (the `$\implies$' part of) Proposition 5.1(ii) to 
$\sigma(q_1)$,$\sigma(q_2)$,$\sigma(q_3)$,$\sigma(r_1)$,$\sigma(r_2)$,$\sigma(r_3)$, 
we conclude that 
$\tilde\sigma(p_1)=(\sigma(q_2)\vee \sigma(q_3))\cap (\sigma(r_2)\vee \sigma(r_3))$, 
$\tilde\sigma(p_2)=(\sigma(q_3)\vee \sigma(q_1))\cap (\sigma(r_3)\vee \sigma(r_1))$, 
$\tilde\sigma(p_3)=(\sigma(q_1)\vee \sigma(q_2))\cap (\sigma(r_1)\vee \sigma(r_2))$ 
are collinear, as desired. 

\medskip\noindent
{\it Step 3-4.} Claim: If $\ell\in\Bbb L(V_1)$, then $\tilde\sigma(\ell)\in \Bbb L(V_2)$. 

Indeed, for each pair $p,q\in\ell$ with $p\neq q$, set 
$\ell'_{p,q}\defeq\tilde\sigma(p)\vee\tilde\sigma(q)\in\Bbb L(V_2)$. 
In fact, the line $\ell'_{p,q}$ does not depend on the pair $p,q$. Indeed, 
let $p',q'\in\ell$ be any pair with $p'\neq q'$. 
If $\sharp\{p,q,p',q'\}=2$, 
then it is clear that $\ell'_{p,q}=\ell'_{p',q'}$. 
If $\sharp\{p,q,p',q'\}=3$, 
then it follows from Step 3-3 that $
\tilde\sigma(p), \tilde\sigma(q), \tilde\sigma(p'), \tilde\sigma(q')$ 
are collinear, hence $\ell'_{p,q}=\ell'_{p',q'}$. 
If $\sharp\{p,q,p',q'\}=4$, 
then, again by Step 3-3, one has 
$$\ell'_{p,q}=\ell'_{p,q'}=\ell'_{p',q'}.$$
Now, write $\ell'=\ell'_{p,q}$ for some (or, equivalently, all) pair 
$p,q\in \ell$ with 
$p\neq q$. Here, note that, 
by Proposition 5.1(i)(III), 
at least one such pair $p,q$ exists. Then, by definition, 
$\tilde\sigma(\ell)\subset\ell'$ or, equivalently, 
$\sigma(\ell)\subset\tilde\sigma^{-1}(\ell')$. 

Applying this to 
$\sigma^{-1}: U_2\isom U_1$ and $\ell'\in\Bbb L(V_2)$, 
we obtain $\tilde\sigma^{-1}(\ell')=\widetilde{\sigma^{-1}}(\ell')\subset\ell''$ 
for some $\ell'\in\Bbb L(V_1)$. Combining these containment relations, 
we conclude $\ell\subset\ell''$, which implies $\ell=\ell''$ and 
$\tilde\sigma(\ell)=\ell'$, as desired. 

\medskip\noindent
{\it Step 3-5.} End of Step 3. 

By Step 3-4, 
we may define $\tilde\tau: \Bbb L(V_1)\to\Bbb L(V_2)$ to be the 
map that sends $\ell\in \Bbb L(V_1)$ to $\tilde\sigma(\ell)\in \Bbb L(V_2)$. 
(Note that 
$\tilde\tau: \Bbb L(V_1)\to\Bbb L(V_2)$ is an extension of 
$\tau: \Bbb L(V_1)_{U_1}\to\Bbb L(V_2)_{U_2}$, by Step 3-2.)
Applying this to 
$\sigma^{-1}: U_2\isom U_1$, 
we may also define 
$\tilde\tau': \Bbb L(V_2)\to\Bbb L(V_1)$ to be the map 
that sends $\ell'\in \Bbb L(V_2)$ to $\tilde\sigma^{-1}(\ell')
=\widetilde{\sigma^{-1}}(\ell')\in \Bbb L(V_1)$. By definition, 
it is immediate to prove that $\tilde\tau'$ is the inverse map of 
$\tilde\tau$, and, in particular, that $\tilde\tau$ is a bijection. 
Now, by the very definition of $\tilde\tau$, 
$(\tilde\sigma,\tilde\tau): 
(\Bbb P(V_1), \Bbb L(V_1))\isom (\Bbb P(V_2), \Bbb L(V_2))$ 
is a collineation, as desired. 

\medskip\noindent
{\it Step 4.} Uniqueness. 
If 
$(\tilde\sigma, \tilde\tau): (\Bbb P(V_1), \Bbb L(V_1))\isom (\Bbb P(V_2), \Bbb L(V_2))$ 
is a collineation, then, for each $p\in\Bbb P(V_1)$, the bijection 
$\tilde\tau:\Bbb L(V_1)\isom\Bbb L(V_2)$ induces a bijection 
of subsets $\Bbb L(V_1)_{\{p\}}\isom\Bbb L(V_2)_{\{\tilde\sigma(p)\}}$. 
So, if $(\tilde\sigma, \tilde\tau)$ extends $(\sigma, \tau): 
(U_1, \Bbb L(V_1)_{U_1})\isom (U_2, \Bbb L(V_2)_{U_2})$, 
then $\tilde\tau$ (or, equivalently, $\tau$) induces a bijection of subsets 
$\Bbb L(V_1)_{\{p\}}\cap \Bbb L(V_1)_{U_1}
\isom\Bbb L(V_2)_{\{p'\}}\cap \Bbb L(V_2)_{U_2}$. Thus, the uniqueness assertion 
in Theorem 5.7 follows from the uniqueness assertion in Step 2. (Recall the fact 
that $\tilde\tau$ is determined by $\tilde\sigma$ uniquely.) 
\qed\enddemo

$$\text{References.}$$
\noindent
[Artin] E. Artin, Geometric algebra, Interscience Publishers, Inc., 
1957. 

\noindent
[EDM] Encyclopedic dictionary of mathematics, Vol. I--IV (Translated from Japanese), Second edition, 
Kiyosi It\^o ed., MIT Press, 
1987.

\noindent
[Grothendieck]\ Grothendieck, A., Brief an G. Faltings, (German), with an English 
translation on pp. 285-293, London Math. Soc. Lecture Note Ser., 242, Geometric Galois 
actions, 1, 49--58, Cambridge Univ. Press, Cambridge, (1997).

\noindent
[Mochizuki1] Mochizuki, S., Absolute anabelian cuspidalizations of proper hyperbolic curves, 
J. Math. Kyoto Univ. 47 (2007), 
451--539. 

\noindent
[Mochizuki2] Mochizuki, S., Topics in absolute anabelian geometry I: generalities, J. Math. Sci. Univ. Tokyo 
19 (2012), 
139--242.

\noindent
[Sa\"\i di-Tamagawa1]\ Sa\"\i di, M., and Tamagawa, A., A prime-to-$p$ version of the Grothendieck 
anabelian conjecture for hyperbolic curves in characteristic $p>0$, Publ. RIMS, Kyoto Univ. 45 (2009), 135--186.

\noindent
[Sa\"\i di-Tamagawa2]\ Sa\"\i di, M., and Tamagawa, A., 
On the anabelian geometry of hyperbolic curves over finite fields, 
in Algebraic Number Theory and Related Topics 2007, RIMS Kokyuroku Bessatsu B12, 
RIMS, Kyoto Univ., 2009, 67--89. 

\noindent
[Sa\"\i di-Tamagawa3]\ Sa\"\i di, M., and Tamagawa, A., A refined version of Grothendieck's birational 
anabelian conjecture for curves over finite fields, preprint, submitted. 

\noindent
[Tamagawa] Tamagawa, A., 
The Grothendieck conjecture for affine curves, 
Compositio Math. 109 (1997), 
135--194.

\bigskip

\noindent
Mohamed Sa\"\i di

\noindent
College of Engineering, Mathematics, and Physical Sciences

\noindent
University of Exeter

\noindent
Harrison Building

\noindent
North Park Road

\noindent
EXETER EX4 4QF 

\noindent
United Kingdom

\noindent
M.Saidi\@exeter.ac.uk

\bigskip
\noindent
Akio Tamagawa

\noindent
Research Institute for Mathematical Sciences

\noindent
Kyoto University

\noindent
KYOTO 606-8502

\noindent
Japan

\noindent
tamagawa\@kurims.kyoto-u.ac.jp
\enddocument